\documentclass[12pt]{amsart}
\usepackage{color,currfile,pdfpages}
\usepackage{epstopdf}
\usepackage{amsmath,amssymb,amsthm}
\usepackage{ifthen,verbatim,here}
\usepackage{mathrsfs}
\usepackage{yhmath}
\usepackage{booktabs}
\usepackage{underscore}
\usepackage{longtable}
\usepackage{url}
\usepackage[english]{babel}
\usepackage{verbatim,here}
\usepackage[T1]{fontenc}
\usepackage{floatflt,graphicx,graphics,subfig}
\usepackage{a4wide}
\numberwithin{equation}{section}
\nonstopmode
\theoremstyle{plain}

\newtheorem{lem}[equation]{Lemma}

\newtheorem{thm}[equation]{Theorem}
\theoremstyle{definition}
{\qed\bigskip}

\newcounter{alphabet}

\newcommand{\be}{\begin{eqnarray}}
	\newcommand{\ee}{\end{eqnarray}}
\newcommand{\ba}{\begin{array}}
	\newcommand{\ea}{\end{array}}
\newcommand{\ben}{\begin{eqnarray*}}
	\newcommand{\een}{\end{eqnarray*}}

\newcommand{\B}{\mathbb{B}}

\DeclareMathOperator{\capa}{\mathrm{cap}}

\let\th\relax
\DeclareMathOperator{\th}{\mathrm{th}}
\DeclareMathOperator{\sh}{\mathrm{sh}}

\DeclareMathOperator{\cth}{\mathrm{cth}}

\newcommand{\cT}{\mathcal{T}}
\newcommand{\eE}{\varepsilon}
\newcommand{\osc}{\mathrm{osc}}

\newcommand{\arsh}{\,\mathrm{arsh}\,}
\newcommand{\arth}{\,\mathrm{arth}\,}

\let\Re\relax
\DeclareMathOperator{\Re}{\mathrm{Re}}
\let\Im\relax
\DeclareMathOperator{\Im}{\mathrm{Im}}

\newcommand {\M} {\mathsf{M}}

\renewcommand{\i}{\mathrm{i}}
\newcommand{\bs}{{\bf s}}

\newcommand{\bI}{{\bf I}}
\newcommand{\bM}{{\bf M}}
\newcommand{\bN}{{\bf N}}
\renewcommand{\Im}{{ \rm Im}\,}
\renewcommand{\Re}{{ \rm Re}\,}

\newcounter{minutes}
\divide\time by 60
\newcounter{hours}
\multiply\time by 60
\addtocounter{minutes}{-\time}

\begin{document}
	
	\title[Mobile disks in hyperbolic space]
	{
		Mobile disks in hyperbolic space and minimization of conformal capacity 
	}
	
	\def\thefootnote{}
	\footnotetext{
		\texttt{\tiny File:~\jobname .tex,
			printed: \number\year-\number\month-\number\day,
			\thehours.\ifnum\theminutes<10{0}\fi\theminutes}
	}
	\makeatletter\def\thefootnote{\@arabic\c@footnote}\makeatother

\author[H. Hakula]{Harri Hakula}
\address{Aalto University, Department of Mathematics and Systems Analysis, P.O. Box 11100, FI-00076 Aalto, FINLAND}
\email{harri.hakula@aalto.fi}
\author[M. Nasser]{Mohamed M. S. Nasser}
\address{Department of Mathematics, Statistics, and Physics, Wichita State University, Wichita, KS 67260-0033, USA}
\email{mms.nasser@wichita.edu}
\author[M. Vuorinen]{Matti Vuorinen}
\address{Department of Mathematics and Statistics, University of Turku, FI-20014 Turku, Finland}
\email{vuorinen@utu.fi}
	
\begin{abstract}
Our focus is to study constellations of disjoint disks in the hyperbolic
space, the unit disk equipped with the hyperbolic metric. Each constellation
corresponds to a set $E$ which is the union of  $m>2$ disks with
hyperbolic radii $r_j>0, j=1,...,m$. The centers of the disks are not
fixed and hence individual disks of the constellation are allowed to
move under the constraints that they do not overlap and their 
hyperbolic radii remain invariant.  Our main objective is to find
computational lower bounds for the conformal capacity of a given
constellation. The capacity depends on the centers and radii in a very
complicated way even in the simplest cases when $m=3$ or $m=4$.  In the
absence of analytic methods our work is based on numerical simulations
using two different numerical methods, the boundary integral equation
method and the $hp$-FEM method, resp. Our simulations combine capacity
computation with minimization methods and produce extremal cases where
the disks of the constellation are grouped next to each other. This
resembles the behavior of animal colonies minimizing
heat flow in arctic areas.   
\end{abstract}
	
\maketitle

\section{Introduction}
Many extremal problems of physics, exact sciences, and mathematics have 
solutions
which exhibit varying degree of symmetry. A typical situation is to minimize or
maximize
a set functional of a planar set under the constraint that some other functional
is constant. The classical isoperimetric problem \cite{PS} is an example. Here 
	one maximizes the area of a planar set given its perimeter and the extremal 
domain is the disk.  G. P\'{o}lya and G. Szeg\"{o} \cite{PS} initiated a
systematic study of a large class of isoperimetric type problems of mathematical
physics for domain functionals such as moment of inertia, principal frequency,
torsional rigidity, and, in particular, capacities of condensers. 
Certain geometric transformations, known under the general name
``symmetrization'' have the property that they decrease the value of domain 
functionals and thus can give hints about the extremal configuration of
isoperimetric problems \cite{B,Du}. We study here new types of transformations
which decrease the value of conformal capacity. 

\begin{figure}
	\centering
	\subfloat[]{
		\includegraphics[width=0.6\textwidth]{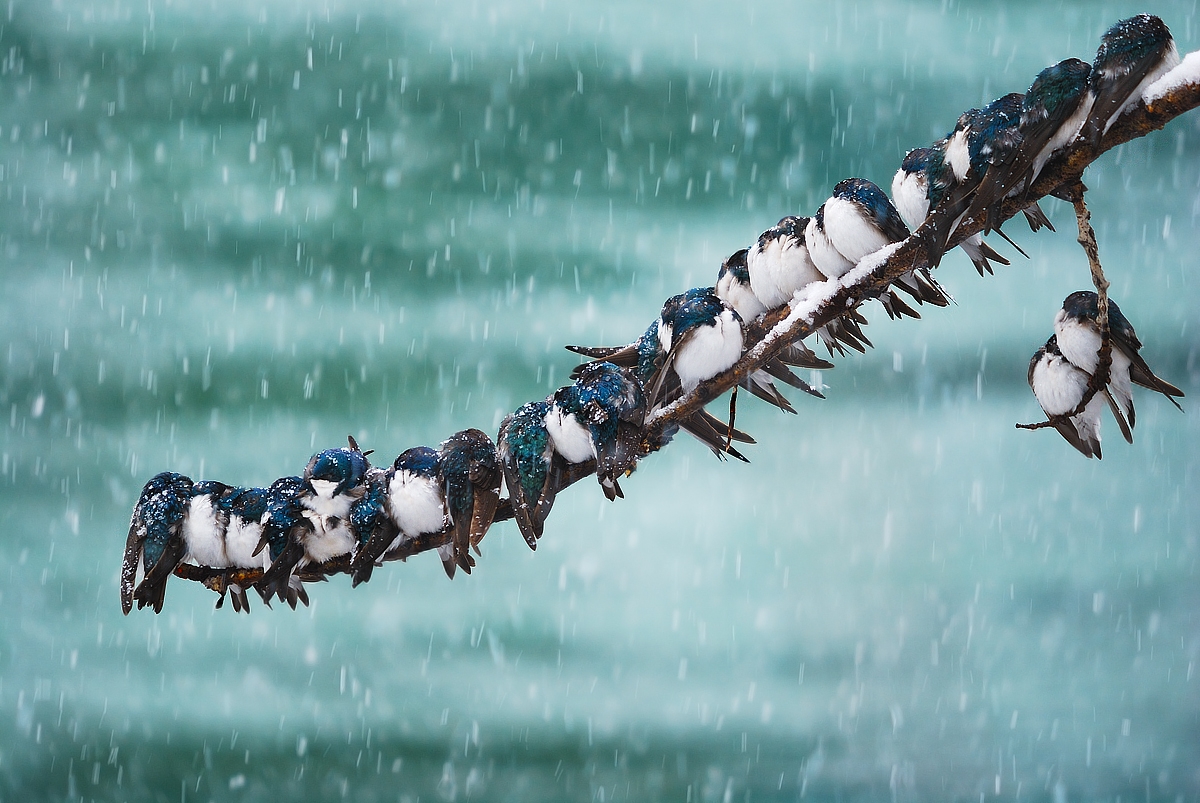}
	}\\
	\subfloat[]{
		\includegraphics[width=0.35\textwidth]{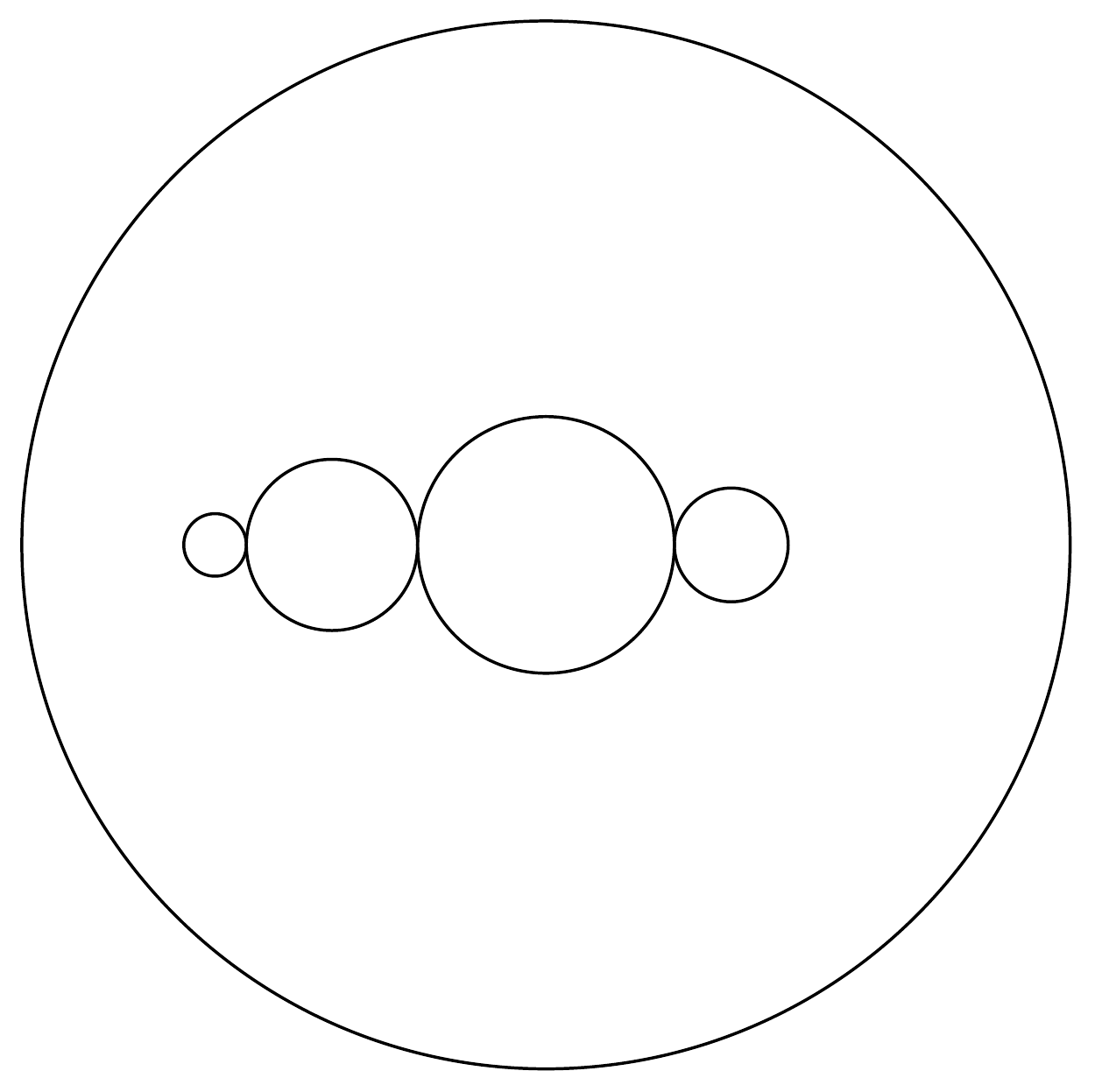}
	}
	\subfloat[]{
		\includegraphics[width=0.35\textwidth]{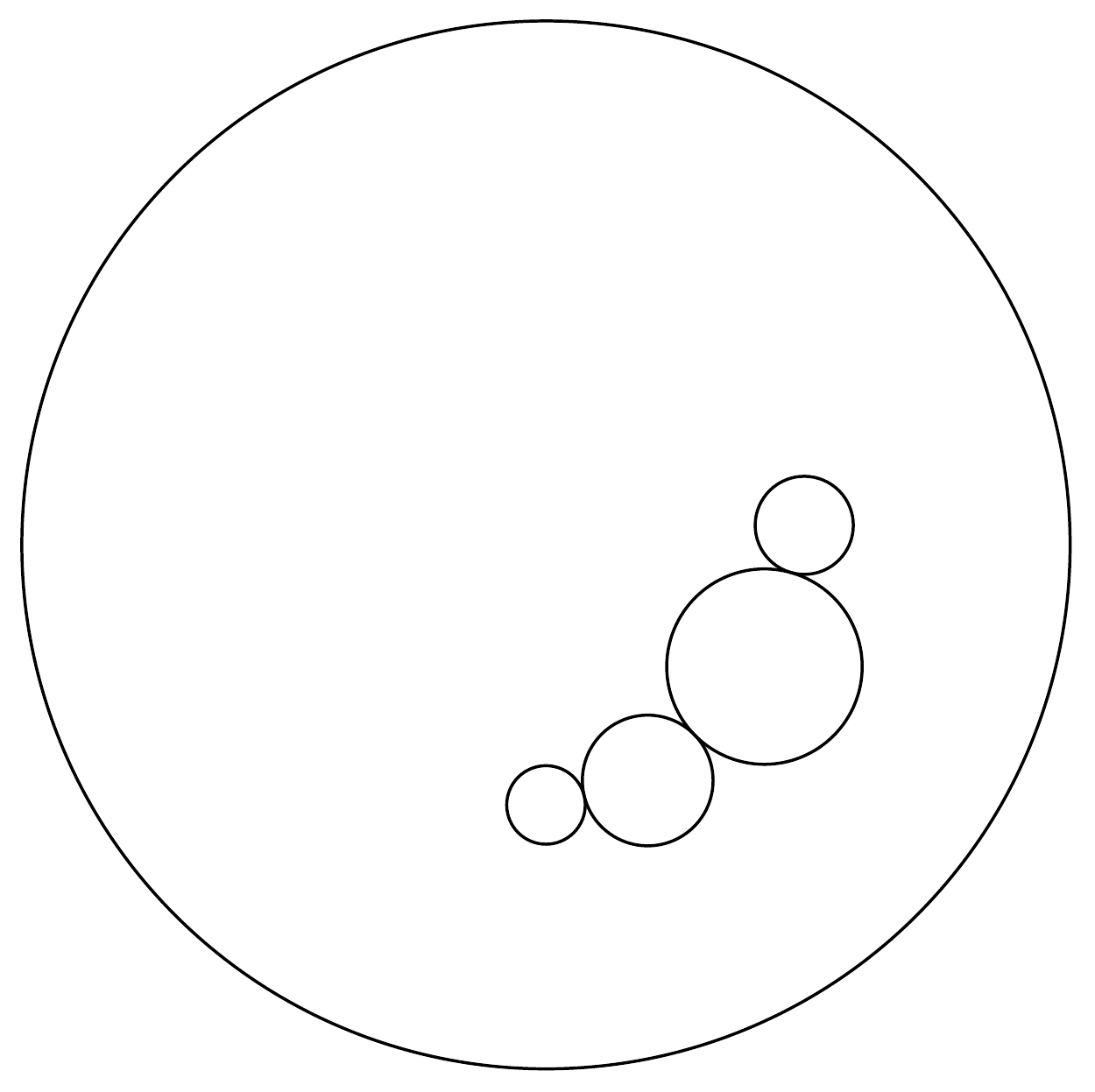}
	}
	
	\caption{Examples of constrained optimisation. (a) Tree swallows huddle on a branch during a spring snowstorm \cite{Wil}.
		(b) Minimal capacity configuration for four hyperbolic disks on a diameter.
		(c) Minimal capacity configuration for four hyperbolic disks on a hyperbolic circle.
		In (b) and (c) the hyperbolic disks are inside the unit disk equipped with the hyperbolic metric.}\label{fig:appetizer}
\end{figure}

%
In a very interesting recent paper,
A. Solynin \cite{S3} describes capacity problems, motivated by
the behavior of herds of arctic animals which keep close together to 
minimize the total loss of heat of the herd or to defend against predators
(see figures in \cite{S3}). 
Such a herd behavior seems to suggest the heuristic idea that 
``minimization of herd's outer perimeter'' minimizes the loss of heat 
or danger from predators. This kind of extremal problem can be 
classified as special type of isoperimetric problem.
As an illustration of the connection between the kind of transformations
we are interested in and the observed behavior in nature, see Figure~\ref{fig:appetizer}.

In a recent paper \cite{NV2}, isoperimetric inequalities in hyperbolic 
geometry  were applied to estimate the conformal capacity
of condensers of the form $(\mathbb{B}^2, E)$ where $E $ is a union of 
finitely many disjoint closed disks $E_j, j=1,...,m,$ in the unit 
disk $\mathbb{B}^2$. Thus $E$ is a constellation of disks.
Gehring's lower bound~\cite{G} (see also \cite{NV2}) is given by condensers of the form 
$(\mathbb{B}^2, E^*)$ where $E^* $ is a disk with the hyperbolic area
equal to that of $ \cup_{j=1}^m E_j.$
Further recent investigations of condenser capacity
in the framework of hyperbolic geometry include 
\cite{nv1, NRV1, NRV2}, where pointers to earlier work can be found. It
should be noticed that due to the conformal invariance of the conformal
capacity, the hyperbolic geometry provides the natural setup for this 
study.

We continue here this work and our goal is to analyse extremal cases of 
the aforementioned capacity and how the capacity depends on the geometry
of
the disk constellation. 
The constraint that the disks do not overlap leads to problems of combinatorial
geometry. Some examples of such geometric problems, related to this work and
the herd behavior mentioned above, are Descartes' problem of four circles with each circle tangent to three
circles, Apollonian circle packing, and Soddy's ``complex kiss precise'' problem for
configurations of mutually tangent circles \cite{LMW}. 
Combinatorial geometry extremal problems motivated by  biochemistry 
research and drug  development  are described in \cite{LB}. 
A very interesting discussion of many topics of combinatorial geometry
including packing problems is given in the encyclopedic work of 
M. Berger \cite{ber}.
The three dimensional case is much more
difficult than the planar case and it is the subject
of the extensive review paper \cite{KKLS} where topics range from optimal packing
of spheres to constrained motion of small spheres on the surface of the unit
sphere. For an extensive survey of potential theoretic extremal problems
see \cite{bhs}.

Analysing the extremal cases of the lower bound for
\[
{\rm cap} (\mathbb{B}^2, \cup_{j=1}^m E_j)
\]
for a constellation of disjoint hyperbolic disks $E_j$ seems to be very difficult even in the
simplest cases $m=3,4.$
Therefore we consider this problem in special cases such as 
the case when the circle
centers are at the same distance from origin or analyse constrained motion
of one circle along three other fixed circles (see Figure~\ref{fig:appetizer}). 
Simulations indicate that several constellations yield local minima of the capacity.
Throughout, the hyperbolic 
geometry provides the natural geometric  framework for our study, because of
the conformal invariance of the capacity.
We use two numerical methods for computing the capacity, the $hp$-FEM and
the boundary integral equation (BIE) method.
The numerical results lead to a number of conjectures and improved bounds.
Indeed, the existing lower bound for constellations considered here is improved of the order of 10\%
for disks of unit hyperbolic radius. 
Moreover, the asymptotic nature of the
theoretical lower bound as the hyperbolic radii $r_j \to \infty$ 
is easily understood in the context of
hyperbolic geometry.

In modern physics, in particular in condensed matter physics,
there has been a lot of interest in geometric settings with
negative curvature  \cite{Kol,Len}, that is, exactly our natural setup.
The purpose of this paper is also to show how computations can be formulated and
performed in both Euclidean and hyperbolic geometries, even with the possibility
of moving from one to another. This is highlighted in the last section
where the optimal configurations in hyperbolic geometry are found
by successive transformations to a Euclidean coordinate system employed in the optimization routines. For information about potential theory and its applications, see~\cite{bhs,PS,R,Tsuji}.

The contents are organized into sections as follows. Section 2 contains
the key facts about hyperbolic geometry, including the transformation formulae
from Euclidean disks to Poincar\'e disks and back.
Section 3 covers the preliminary notations of conformal capacity, 
collected from various sources, e.g. from \cite{Be,Du,GMP,GoRe,HKM,HKV}.
These are the cornerstones of the geometric setup of the computations 
in the sequel. Section 3 also provides an overview of the $hp$-FEM  \cite{HRV1,hno} 
adjusted to the
present computational tasks, our second computational work horse,
the BIE method \cite{Nas-ETNA,nv1}, and
the interior-point method used in optimization.
The numerical experiments are discussed in Sections 4 and 5.
In Section 4 the selected configurations have been designed a priori, with the
goal of forming an understanding of the identifiable geometric features
of the minimal capacity configurations.
In Section 5 that understanding is challenged by searching for the minimal
capacity configurations using numerical optimization starting with
random initial configurations. 
Finally, the conclusions are drawn in Section 6.

\section{From Euclidean Disk to Poincare and Back}

\begin{figure}
	\centering
	\includegraphics[width=0.35\textwidth]{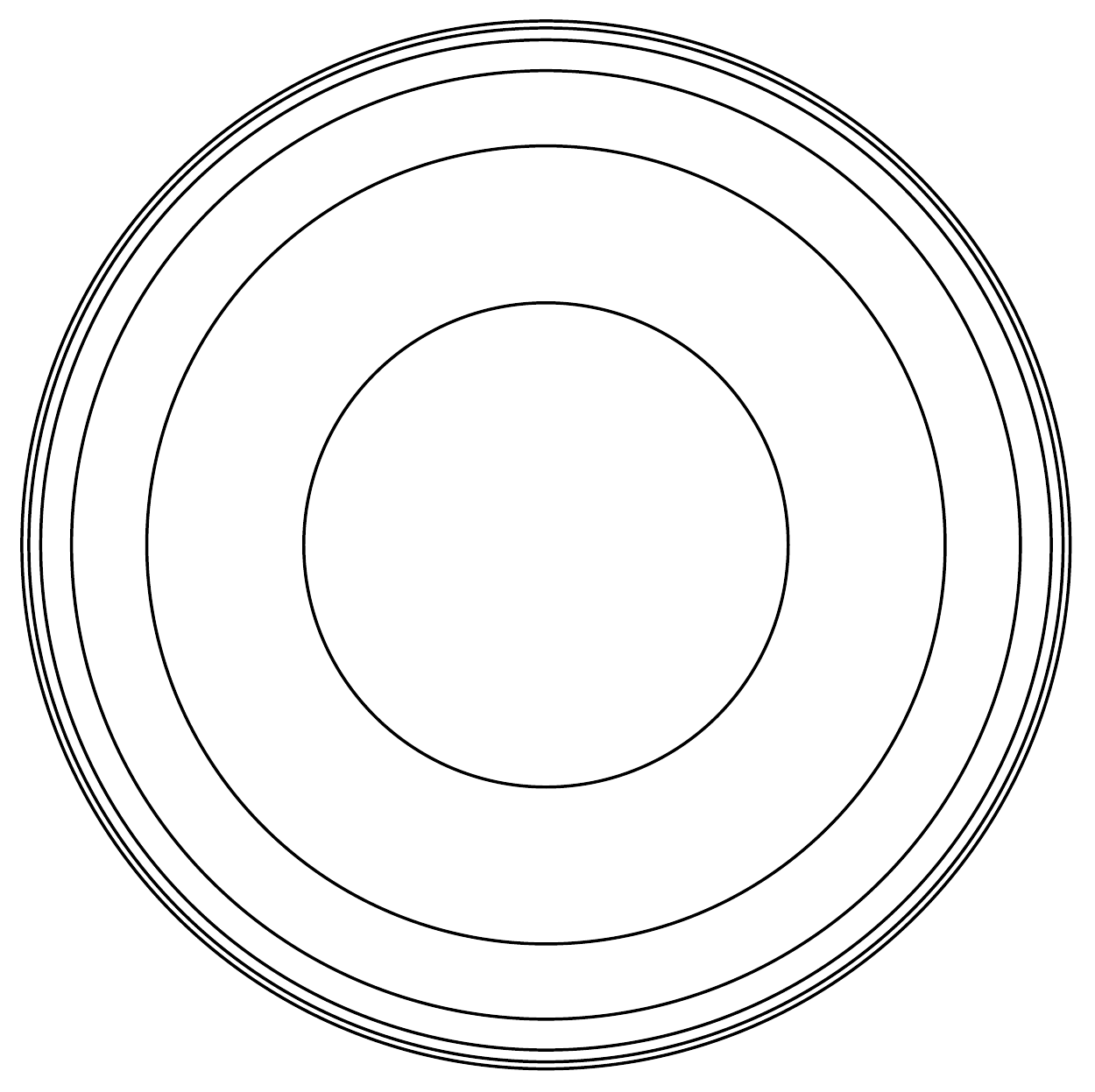}
	\hfill
	\includegraphics[width=0.35\textwidth]{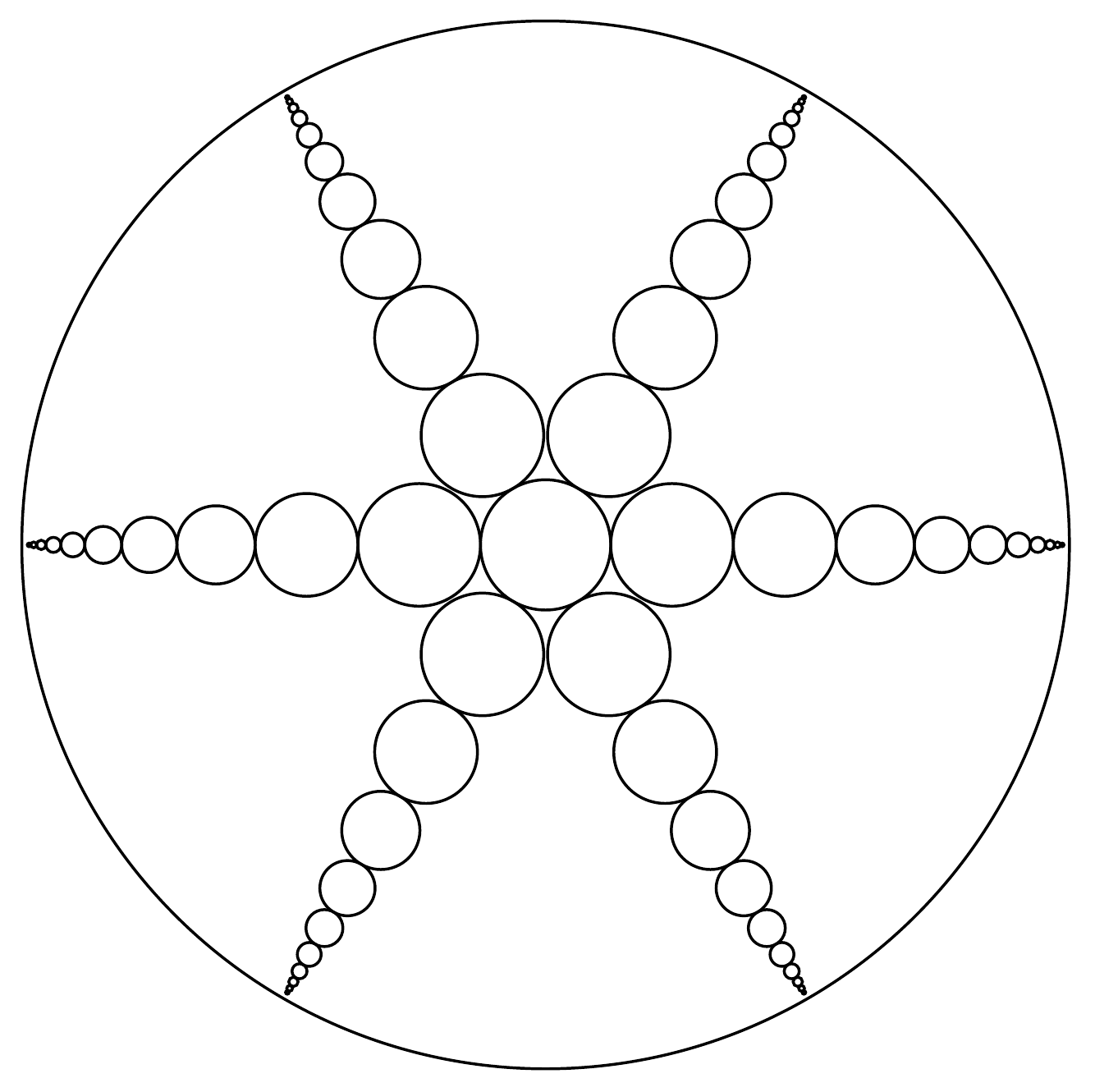}
	\caption{Visualisations on Poincar\'e disk. Left: Images of hyperbolic disks with hyperbolic  radii $= 1,2,3,4,5.$ Right: Hyperbolic disks on three diameters all with equal radii. Notice the lens-shaped regions containing the disks on each diameter.}\label{fig:hypgeom}
\end{figure}

In this section the central transformation formulae collected from various sources
are presented. In Figure~\ref{fig:hypgeom} different properties of
geometry on the Poincar\'e disk have been illustrated. In particular, the 
facts that for all $\epsilon>0$, $M>0$ there are hyperbolic disks with radii $M$ but Euclidean diameter $<\epsilon$ and hyperbolic
disks with equal radii have different Euclidean radii depending
on their location are important for our discussion below.

For a point $x\in\mathbb{R}^n$ and a radius $r>0$,  define an open Euclidean ball 
$B^n(x,r)=\{y\in\mathbb{R}^n\text{ }|\text{ }|x-y|<r\}$ and its boundary sphere 
$S^{n-1}(x,r)=\{y\in\mathbb{R}^n\text{ }|\text{ }|x-y|=r\}$. For the unit ball and 
sphere, we use the simplified notations $\mathbb{B}^n=B^n(0,1)$ and 
$S^{n-1}=S^{n-1}(0,1)$.
The segment joining two points $x, y\in\mathbb{R}^n$ is denoted $[x,y].$

Define the hyperbolic metric in the Poincar\'e unit disk 
$\mathbb{B}^2$ as in  \cite{Be}, \cite[(2.8) p. 15]{BM}
\begin{align}\label{myrho}
	\text{sh}^2\frac{\rho_{\mathbb{B}^2}(x,y)}{2}&=
	\frac{|x-y|^2}{(1-|x|^2)(1-|y|^2)},\quad x,y\in \mathbb{B}^2.
\end{align}
We use the notation $\sh$ and $\arsh$ for the hyperbolic sine and its inverse, respectively, and similarly, $\th$ and $\arth$ for the hyperbolic tangent and its inverse. 
The hyperbolic midpoint of $x,y\in\mathbb{B}^2$ is given by~\cite{WVZ}
\begin{equation*}
	m_H(x,y) = \frac{y \left(1-| x| ^2\right)+x \left(1-| y| ^2\right)}{1-|x|^2|y|^2 + A[x,y] \sqrt{(1-|x|^2)(1-|y|^2)}}
\end{equation*}
where $A[x,y]= \sqrt{|x-y|^2 +(1-|x|^2)(1-|y|^2)}$.
We use the notation
\[ B_{\rho}(x,M) = \{z \in \mathbb{B}^2: \rho_{\mathbb{B}^2}(x,z)<M\}\]
for the hyperbolic disk centered at $x\in \mathbb{B}^2$ with radius $M>0\,.$
It is a basic fact that they are Euclidean disks with
the center and radius given
by \cite[p.56, (4.20)]{HKV}
\begin{equation}\label{hypdisk}
	\begin{cases}
		B_\rho(x,M)=B^2(y,r)\;,&\\
		\noalign{\vskip5pt}
		{\displaystyle y=\frac{x(1-t^2)}{1-|x|^2t^2}\;,\;\;
			r=\frac{(1-|x|^2)t}{1-|x|^2t^2}\;,\;\;t={\th} ( M/2)\;,}&
	\end{cases}
\end{equation}
Note the special case $x= 0$,
\begin{equation} \label{hypDat0}
	B_\rho(0,M)=B^2(0,{\th} ( M/2)) \,.
\end{equation}
Conversely, the Euclidean disks can be considered as hyperbolic ones by
\cite{WVZ}
\begin{equation}\label{eucdisk}
	\begin{cases}
		B^2(y,r) = B_\rho(x,M)\;,&\\
		\noalign{\vskip5pt}
		{\displaystyle 
			x=t\,y/|y|\;,\;\;
			M=\rho_{\mathbb{B}^2}(x,z)\;,\;\;t={m_H}\, (|y|-r, |y|+r)\;,}&
	\end{cases}
\end{equation}

\begin{lem}[{\cite[Thm 7.2.2, p. 132]{Be}}] \label{perimBe}
	The area of a hyperbolic disc of radius $r$ is 
	$4\pi\sh^2(r\slash2)$ and  the length of a hyperbolic 
	circle of radius $r$ is $2\pi\sh(r)$.
\end{lem}

\section{Conformal Capacity and Numerical Methods}\label{prel}

A \emph{condenser} is a pair $(G,E)$, where $G \subset \mathbb{B}^2$ 
is a domain and $E$ is a compact non-empty subset of $G$. 
The \emph{ conformal capacity} of this condenser is defined as  
\cite{Du,GMP,GoRe,HKV,HKM}
\begin{align}\label{def_condensercap}
	{\rm cap}(G,E)=\inf_{u\in A}\int_{G}|\nabla u|^2 dm,
\end{align}
where $A$ is the class of $C^\infty_0(G)$ functions $u: G\to[0,\infty)$ 
with  $u(x) \ge 1$ for all $x \in E$ and $dm$ is the $2$-dimensional Lebesgue measure.
In this paper we assume that $G=\mathbb{B}^2$ is the unit disk and $E=\cup_{j=1}^{m}E_j$ where $E_1,\ldots,E_m$ are $m$ closed disjoint disks in the unit disk. Hence $\Omega=G\backslash E$ is a multiply connected circular domain of connectivity $m+1$. In this case, the infimum is attained by a function $u$ which is harmonic in $\Omega$ and satisfies the boundary conditions $u=0$ on $\partial G$ and $u=1$ on $\partial E$~\cite{Du}. The capacity can be expressed in terms of this extremal function as
\begin{equation}\label{eq:cap}
	{\rm cap}(G,E)=\iint \limits_{\Omega}|\nabla u|^2 dm.
\end{equation}

The conformal capacity of a condenser is one 
of the key notions of potential theory of elliptic partial 
differential equations \cite{HKM, GoRe} and it has numerous 
applications to geometric function theory, both in the plane and 
in higher dimensions, \cite{Du, GMP, GoRe, HKV, HKM}. 
Numerous variants of the definition \eqref{def_condensercap} of 
capacity are given in \cite{GMP, GoRe}. First, the family $A$ may be 
replaced by several other families by \cite[Lemma 5.21, p. 161]{GMP}. 
Furthermore,
\begin{align}
	\capa(G,E)=\M(\Delta(E,\partial G;G)),    
\end{align}
where $\Delta(E,\partial G;G)$ is the family of all curves joining $E$ 
with the boundary $\partial G$ in the domain $G$ and $\M$ stands 
for the modulus of a curve family \cite[Thm 5.23, p. 164]{GMP}. 
For the basic facts about capacities and moduli, the reader is 
referred to \cite{GMP, GoRe, HKV, HKM}.

\subsection{{Numerical Methods}}
In this section the numerical methods used in the numerical experiments
are briefly described. The capacities are computed using the $hp$-version
of the finite element method (FEM) and the boundary integral equation with
the generalized Neumann kernel method (BIE).
The minimization problems are computed using the interior-point method
as implemented in MATLAB and Mathematica.

Since the Dirichlet problem \eqref{def_condensercap}
is one of the primary numerical model problems, 
any standard solution technique can be viewed as having been validated.
Verification of the results is discussed in connection with one of the numerical
experiments below.

\begin{figure}
	\centering
	\subfloat{
		\includegraphics[width=0.29\textwidth]{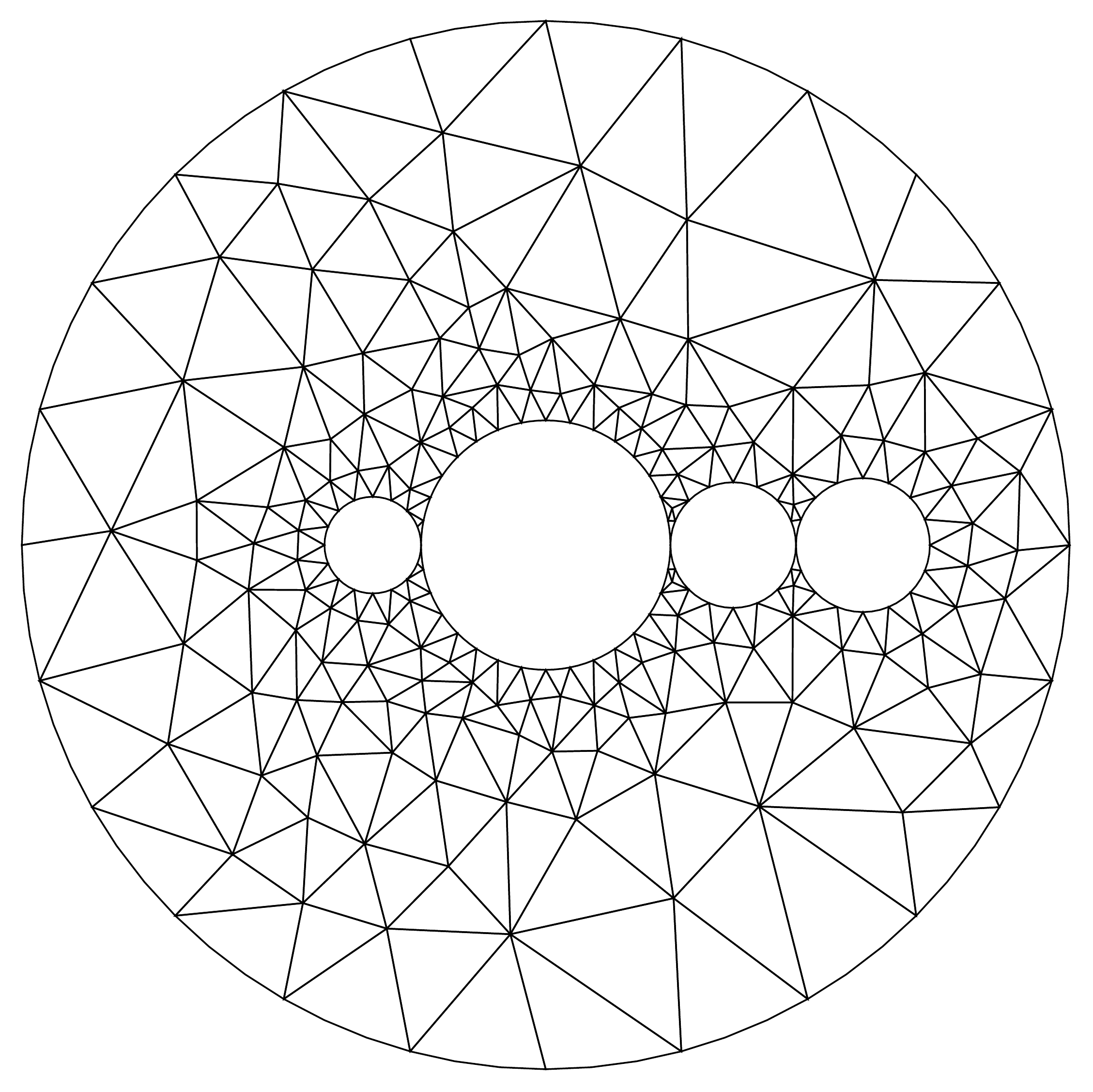}
	}
	\subfloat{
		\includegraphics[width=0.29\textwidth]{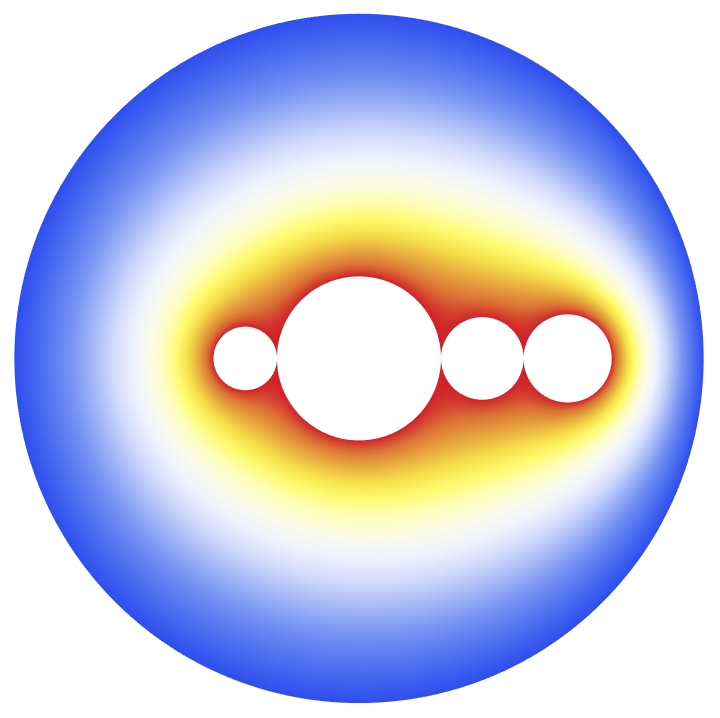}
	}
	\subfloat{
		\includegraphics[width=0.34\textwidth]{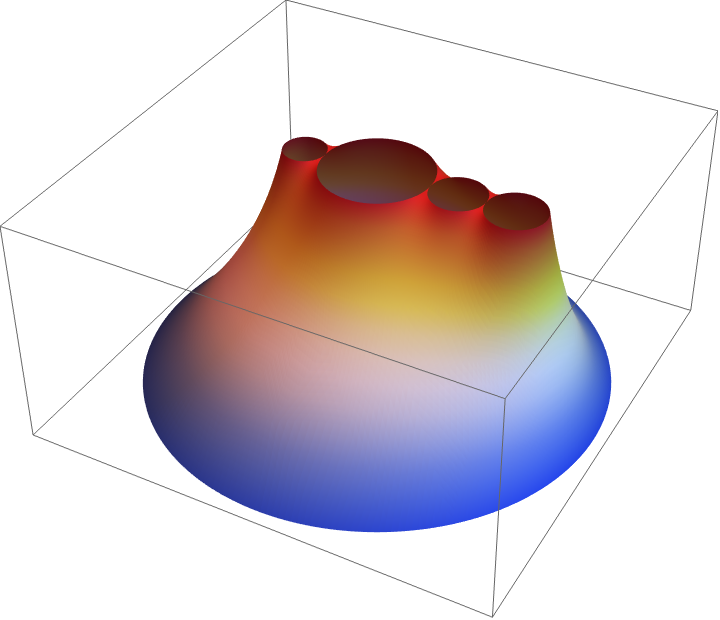}
	}
	\caption{Discretization and optimization. For the given set of four hyperbolic disks with centers constrained on a diameter, the configuration shown here minimizes the capacity.
		Left: Configuration and $hp$-FEM mesh.
		Center: Potential in 2D.
		Right: Potential in 3D.
	}\label{fig:methods}
\end{figure}

\subsubsection{{$hp$-FEM}}\label{sec:hpfem}
What is of particular interest in the context of this paper is that the
$hp$-FEM allows for large curved elements without
significant loss of accuracy. Since the number of elements can be kept
relatively low given that additional refinement can always be
added via elementwise polynomial degree, variation in the boundary
can be addressed directly at the level of the boundary
representation in some exact parametric form. 
This is illustrated in Figure~\ref{fig:methods}.

The following theorem due to Babu{\v{s}}ka and Guo \cite{BaGuo}
sets the limit to the rate of convergence. Notice that construction of
the appropriate spaces is technical. For rigorous treatment of the theory involved
see Schwab \cite{schwab} and references therein. 
\begin{thm} \label{propermesh}
	Let $\Omega \subset \mathbb{R}^2$ be a polygon, $v$ the FEM-solution 
	of \eqref{def_condensercap}, and
	let the weak solution $u_0$ be in a suitable countably normed space where
	the derivatives of arbitrarily high order are controlled.
	Then
	\[
	\inf_v \|u_0 - v\|_{H^1(\Omega)} \leq C\,\exp(-b \sqrt[3]{N}),
	\]
	where $C$ and $b$ are independent of $N$, the number of degrees of freedom. 
	Here $v$ is computed on a proper geometric mesh, where the order of an individual
	element is set to be its element graph distance to the nearest singularity.
	(The result also holds for meshes with constant polynomial degree.)
\end{thm}

Consider the abstract problem setting with $u$ defined on the standard
piecewise polynomial finite element space on some discretization
$\cT$ of the computational domain $\Omega$. Assuming that the exact
solution $u \in H_0^1(D)$ has finite energy, we arrive at the
approximation problem: Find $\hat{u} \in V$ such that
\begin{equation}\label{eq:approximation}
	a(\hat{u},v)=l(v)\ (= a(u,v))\quad (\forall v \in V),
\end{equation}where
$a(\,\cdot\,,\,\cdot\,)$ and $l(\,\cdot\,)$, are the bilinear form
and the load potential, respectively. Additional degrees of
freedom can be introduced by enriching the space $V$. This is
accomplished via introduction of an auxiliary subspace or ``error
space'' $W \subset H_0^1(D)$ such that $V \cap W = \{0\}$. We can
then define the error problem: Find $\varepsilon \in W$ such that
\begin{equation}\label{eq:error}
	a(\varepsilon,v)=l(v)- a(\hat{u},v) (= a(u-\hat{u},v))\quad (\forall v \in W).
\end{equation}
This can be interpreted as a projection of the residual to the auxiliary space.

The main
result on this kind of estimators for the Dirichlet problem \eqref{def_condensercap}
is the following
theorem.
\begin{thm}[\cite{hno}]\label{KeyErrorThm}
	There is a constant $K$ depending only on the dimension $d$,
	polynomial degree $p$, continuity and coercivity constants $C$ and
	$c$, and the shape-regularity of the triangulation $\cT$ such that
	\begin{align*}
		\frac{c}{C}\,\|\eE\|_1\leq\|u-\hat{u}\|_{1}\leq K
		\left(\|\eE\|_{1}+\osc(R,r,\cT)\right),
	\end{align*}
	where the residual oscillation depends on the volumetric and face
	residuals $R$ and $r$, and the triangulation $\cT$.
\end{thm}

\subsubsection{{BIE method}}\label{sec:bie}

We review a BIE method from~\cite{nv1} for computing the capacity $\capa(\B,E)$. The method is based on the BIE with the generalized Neumann kernel. 
The domains considered in this paper are circular domains, i.e., domains whose boundary components are circles. The external boundary is the unit circle, denoted by $C_0$, is parametrized by $\eta_0(t)=e^{\i t}$ for $t\in J_0=[0,2\pi]$. The inner circles $C_j$ are parametrized by $\eta_j(t)=z_j+r_je^{-\i t}$, $t\in J_j=[0,2\pi]$, for $j=1,2,\ldots,m$, where $z_j$ is the center of the circle $C_j$ and $r_j$ is its radius.
Let $J$ be the disjoint union of the $m+1$ intervals $J_j=[0,2\pi]$, $j=0,1,\ldots,m$. We define a parametrization of the whole boundary $C=\cup_{j=0}^m C_j$ on $J$ by (see~\cite{Nas-ETNA} for the details)
\[
\eta(t)=\left\{
\begin{array}{cc} 
	\eta_0(t), & t\in J_0, \\ 
	\eta_1(t), & t\in J_1, \\
	\vdots \\
	\eta_m(t), & t\in J_m. \\ 
\end{array}
\right.
\]

With the parametrization $\eta(t)$ of the whole boundary $C$, we define a complex function $A$ by
\begin{equation}\label{eq:A}
	A(t) = \eta(t)-\alpha,
\end{equation}
where $\alpha$ is a given point in the domain $G$. The generalized Neumann kernel $N(s,t)$ is 
defined for $(s,t)\in J\times J$ by
\begin{equation}\label{eq:N}
	N(s,t) :=
	\frac{1}{\pi}\Im\left(\frac{A(s)}{A(t)}\frac{\eta'(t)}{\eta(t)-\eta(s)}\right).
\end{equation}
We define also the following kernel 
\begin{equation}\label{eq:M}
	M(s,t) :=
	\frac{1}{\pi}\Re\left(\frac{A(s)}{A(t)}\frac{\eta'(t)}{\eta(t)-\eta(s)}\right),
	\quad (s,t)\in J\times J.
\end{equation}
The kernel $N(s,t)$ is continuous and the kernel $M(s,t)$ is singular where the 
singular part involves the cotangent function. Hence, the integral operator $\bN$ 
with the kernel $N(s,t)$ is compact and the integral operator $\bM$ with 
the kernel $M(s,t)$ is singular. Further details can be found in~\cite{Weg-Nas}.

For each $k=1,2,\ldots,m$, let the function $\gamma_k$ be defined by
\begin{equation}\label{eq:gam-k}
	\gamma_k(t)=\log|\eta(t)-z_k|,
\end{equation}
let $\mu_k$ be the unique solution of the BIE
\begin{equation}\label{eq:ie}
	\mu_k-\bN\mu_k=-\bM\gamma_k,
\end{equation}
and let the piecewise constant function $h_k=(h_{0,k},h_{1,k},\ldots,h_{m,k})$ be given by
\begin{equation}\label{eq:h}
	h_k=[\bM\mu_k-(\bI-\bN)\gamma_k]/2.
\end{equation}
For each $k=1,2,\ldots,m$, the solution $\mu_k$ of the BIE~\eqref{eq:ie} and the piecewise constant function $h_k$ in~\eqref{eq:h} will be computed using the MATLAB {\tt fbie} from~\cite{Nas-ETNA}.
In the function {\tt fbie}, the integral equation~\eqref{eq:ie} is solved using the Nystr\"om method with the trapezoidal rule. Solving the integral equation is then reduced to solving an $(m+1)n\times(m+1)n$ linear system which is solved by the MATLAB function {\tt gmres}. The matrix-vector product in {\tt gmres} is computed by the MATLAB function {\tt zfmm2dpart} from the MATLAB toolbox FMMLIB2D~\cite{Gre-Gim12}.
To use the MATLAB function {\tt fbie}, we define a vector $\bs=[s_1,\ldots,s_n]$ where $s_k=2(k-1)\pi/n$, $k=1,\ldots,n$, and $n$ is a given even positive integer. Then we compute the $(m+1)n\times1$ discretization vectors {\tt et} and {\tt etp} of the parametrization $\eta(t)$ of the boundary $C$ and its derivative $\eta'(t)$ by
\[
{\tt et} =[\eta_0(\bs),\eta_1(\bs),\ldots,\eta_m(\bs)]^T, \quad {\tt etp} =[\eta'_0(\bs),\eta'_1(\bs),\ldots,\eta'_m(\bs)]^T.
\]
We also discretize the functions $A(t)$ and $\gamma_k(t)$ by ${\tt A}={\tt et}-\alpha$ and ${\tt gamk}=\gamma_k({\tt et})$, $k=1,\ldots,m$. Then we compute $(m+1)n\times1$ approximate discretizations {\tt muk} and {\tt hk} of the functions $\mu_k(t)$ and $h_k(t)$ by calling
\[
{\tt [muk,hk]=fbie(et,etp,A,gamk,n,5,[\;],1e-14,100)},
\]
i.e., the tolerance of the FMM is $0.5\times10^{-15}$, the GMRES is used without restart, the tolerance of the GMRES method is $10^{-14}$ and the maximal number of GMRES iterations is $100$. 

By computing the $(m+1)n\times1$ vector {\tt hk}, we obtain approximate discretizations of the piecewise constant function $h_k=(h_{0,k},h_{1,k},\ldots,h_{m,k})$ in~\eqref{eq:h}. Note that, for $k=1,\ldots,m$, the constant $h_{j,k}$ is the value of the function $h_k$ on the boundary component $\Gamma_j$. We approximate the values of the real constants $h_{j,k}$ by taking arithmetic means
\[
h_{j,k}=\frac{1}{n}\sum_{i=1+jn}^{(j+1)n}{\tt hk}_i,\quad j=0,1,\ldots,m, \quad k=1,\ldots,m.
\]
The values of the $m$ real constants $a_1,\ldots,a_{m}$ are then approximated by solving the $(m+1)\times(m+1)$ linear system~\cite{nv1}
\begin{equation}\label{eq:sys-method}
	\left[\begin{array}{ccccc}
		h_{0,1}    &h_{0,2}    &\cdots &h_{0,m}      &1       \\
		h_{1,1}    &h_{1,2}    &\cdots &h_{1,m}      &1       \\
		\vdots     &\vdots     &\ddots &\vdots       &\vdots  \\
		h_{m,1}    &h_{m,2}    &\cdots &h_{m,m}      &1       \\
	\end{array}\right]
	\left[\begin{array}{c}
		a_1    \\a_2    \\ \vdots \\ a_{m} \\  c 
	\end{array}\right]
	= \left[\begin{array}{c}
		0 \\  1 \\  \vdots \\ 1  
	\end{array}\right].
\end{equation}
Since $m+1$ is the number of boundary components of the domain $\Omega=G\setminus E$, we can assume that $m$ is small and solve the linear system~\eqref{eq:sys-method} using the Gauss elimination method. By solving the linear system, the capacity $\capa(\B,E)$ will be computed by~\cite[Eq.~(3.9)]{nv1} 
\begin{equation}
	\capa(\B,E)=2\pi\sum_{k=1}^m a_k.
\end{equation}

In this paper, the boundary components of the domain $\Omega$ are circles. Thus, the integrands in~\eqref{eq:ie} and~\eqref{eq:h} will be $2\pi$-periodic functions, and can be extended holomorphically to some parallel strip $|\Im t|<\sigma$ in the complex plane. Hence, the trapezoidal rule will then converge exponentially with $O(e^{-\sigma n})$~\cite{Tre1} when it is used to discretize the integrals in~\eqref{eq:ie} and~\eqref{eq:h}. The numerical solution of the integral equation will converge with a similar rate of convergence~\cite[p.~322]{Atk97} (see Figure~\ref{fig:fourmesherror} (right) below).

\subsubsection{{Nonlinear Optimization: Interior-Point Method}}\label{sec:ipm}
The two methods outlined above are combined with a numerical optimization routine
in the last set of numerical experiments below.
The task is to find an optimal configuration for a set of 
hyperbolic disks $E$ with fixed radii.
We use the interior-point method
as implemented in Mathematica (\texttt{FindMinimum}, \cite{Wol}) and Matlab (\texttt{fmincon},  \cite{Matlab}).

In the most general case the problem is defined as in \eqref{eq:min},
where the only constraint is a geometric one, that is, the
disks are not allowed to overlap. Here, the radii are fixed and the
optimization concerns only the locations of the disks.
\begin{alignat}{3}
	\min_E              &\quad&  {\rm cap}(G,E)  &&             & \nonumber \\
	\text{subject to: } &\quad&  E_i \cap E_j        && = \emptyset &\quad \forall\ i,j=1,\ldots,m, i \neq j \label{eq:min}\\
	&\quad&  E_j                 && \subset G   &\quad \forall\ j=1,\ldots,m. \nonumber 
\end{alignat}
This nonlinear optimization problem can be solved using the interior-point method. 
This solution would be a local minimum.
The standard textbook reference is Nocedal and Wright \cite{Noc}.

Notice, that the objective function is indeed the capacity of the
constellation. Often optimization problems with geometric constraints
are related to packing and fitting problems. The task here is orders of magnitude
more demanding since, at every point evaluation one solution
of the capacity problem has to be computed, and as the disks move the
constraints change as well. The number of evaluations is greater than the number of iteration
steps, since the gradients and Hessians must be approximated numerically.
It should be noted that the success of the optimization depends on the high accuracy of
the capacity solver, since otherwise the approximate derivatives are not sufficiently accurate.

In the context of this work, there have been no attempts to devise a special
method that would incorporate some of the insights gathered during this study.
Instead, the numerical optimization is used to challenge those insights and
therefore the optimizations have been computed with minimal input information.

\section{Minimizing Capacity: Constrained Configurations}\label{ted}
As mentioned above, even with a small number of disks the combinatorial explosion of
the number of configurations is evident. Therefore, we restrict ourselves to a
series of experiments each with increasing complexity building toward an
understanding of the fundamental geometric principles behind the minimal configurations.
In each case we consider a set of hyperbolic disks $E_j$ with radii $r_j$, where
some geometric constraint is placed on all or some of the disks in the constellation.

An initial observation is that due to conformal invariance of the capacity, its
numerical value remains invariant under a M\"obius transformation of the unit
disk onto itself. Therefore we may assume that the disk with the largest
radius $r_1$ is centered at the origin.

Further, consider a disk $B_{\rho}(z_2, r_2)$ with center $z_2$ on 
the segment $(0,1)$. The disk lies in the lens-shaped region
\[
W=B^2(\i\tau,\sqrt{1+\tau^2})\cap B^2(-\i\tau,\sqrt{1+\tau^2}),\quad \tau>0,
\]
with $\rho_{{\mathbb B}^2}(0,\i v)= r_2$ where $v=\sqrt{1+\tau^2}-\tau$ and is
tangent to  both boundary arcs of $W$ and $\pm 1 \in \partial W,$ see
Figure~\ref{fig:hypgeom} (right). Every disk lies within its own associated lens-shaped domain. 

\subsection{Disks with collinear centers}\label{sec:collineardisks} 
Consider a set of $m$ hyperbolic disks $E_j$ with radii $r_j$ 
and centers on the diameter $(-1,1)$ with $\sum_{j=1}^m 2 r_j=d_1=\rho_{\mathbb{B}^2}(-0.6,0.6)$.
We choose the hyperbolic centers of these disks so that the hyperbolic distance between them is 
$d\ge 0$ where $d=0$ corresponds to the case when they touch each other.
The goal is to establish upper and lower bounds for ${\rm cap}(\B, \cup E_j)$.
Since the hyperbolic radius of a hyperbolic disk is invariant under M\"obius transform, in view of~\eqref{hypdisk}, we have ${\rm cap}(\B,E_j)=2\pi/\log(1/\th(r_j/2))$ for all $E_j$.

\begin{figure}
	\centering
	\subfloat{
		\includegraphics[width=0.275\textwidth]{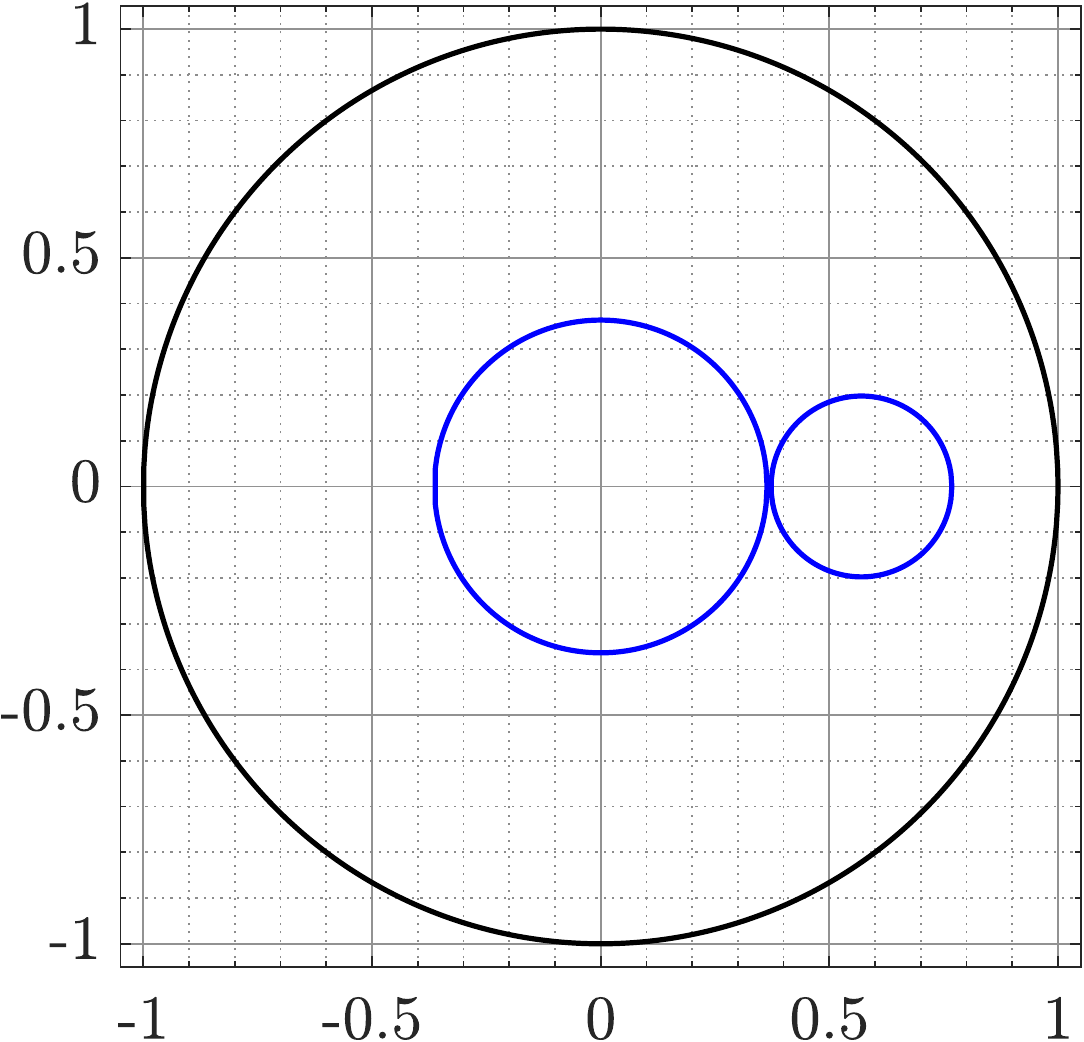}
		\includegraphics[width=0.275\textwidth]{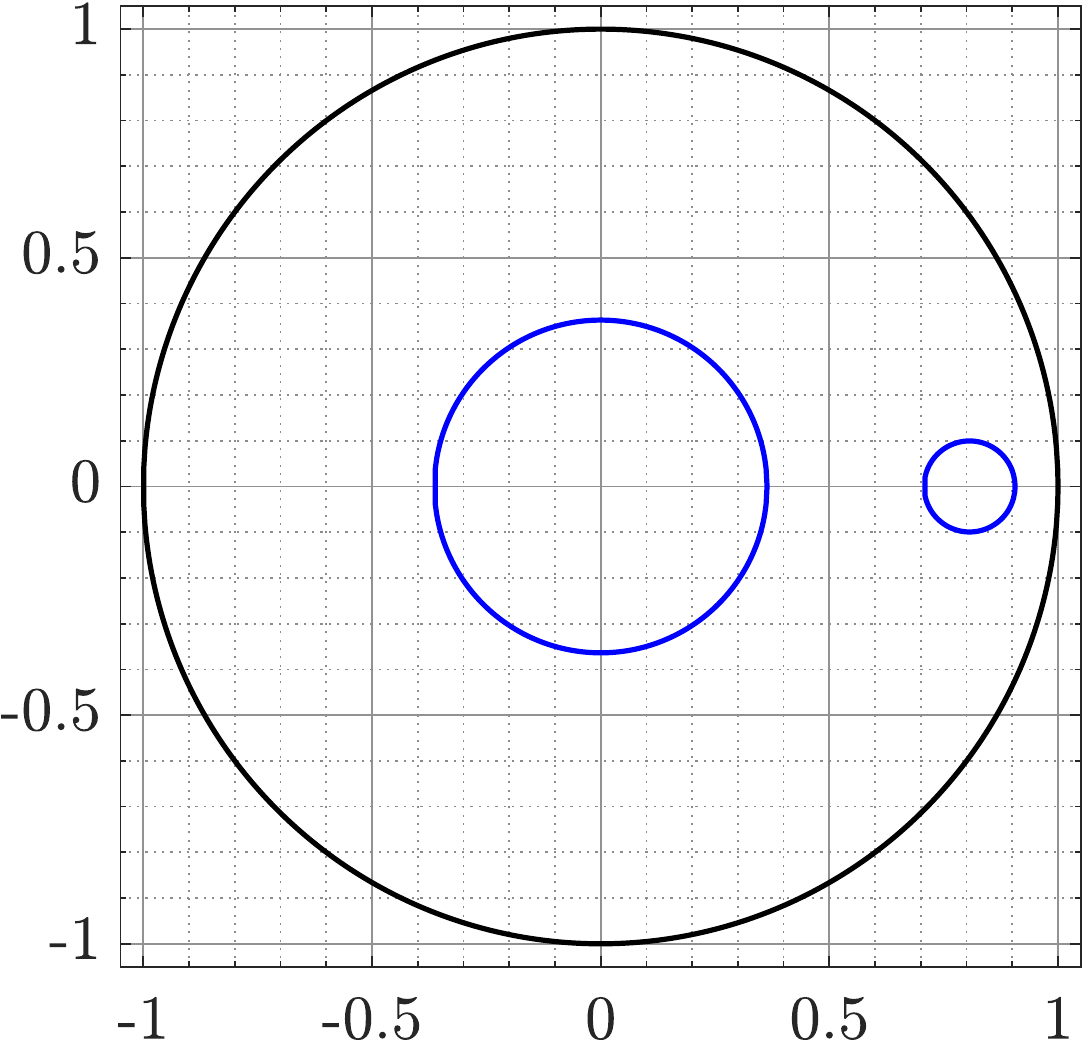}
		\includegraphics[width=0.2525\textwidth]{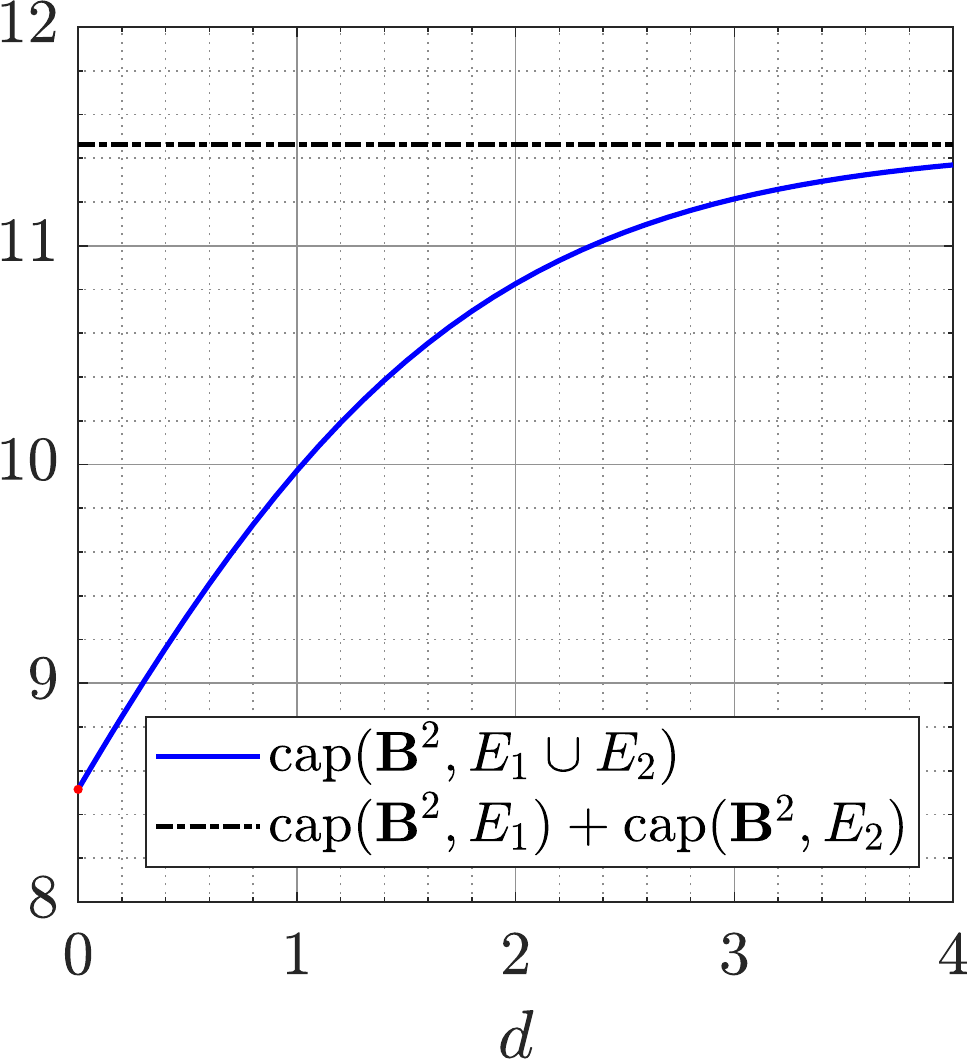}
	}\\
	\subfloat{
		\includegraphics[width=0.275\textwidth]{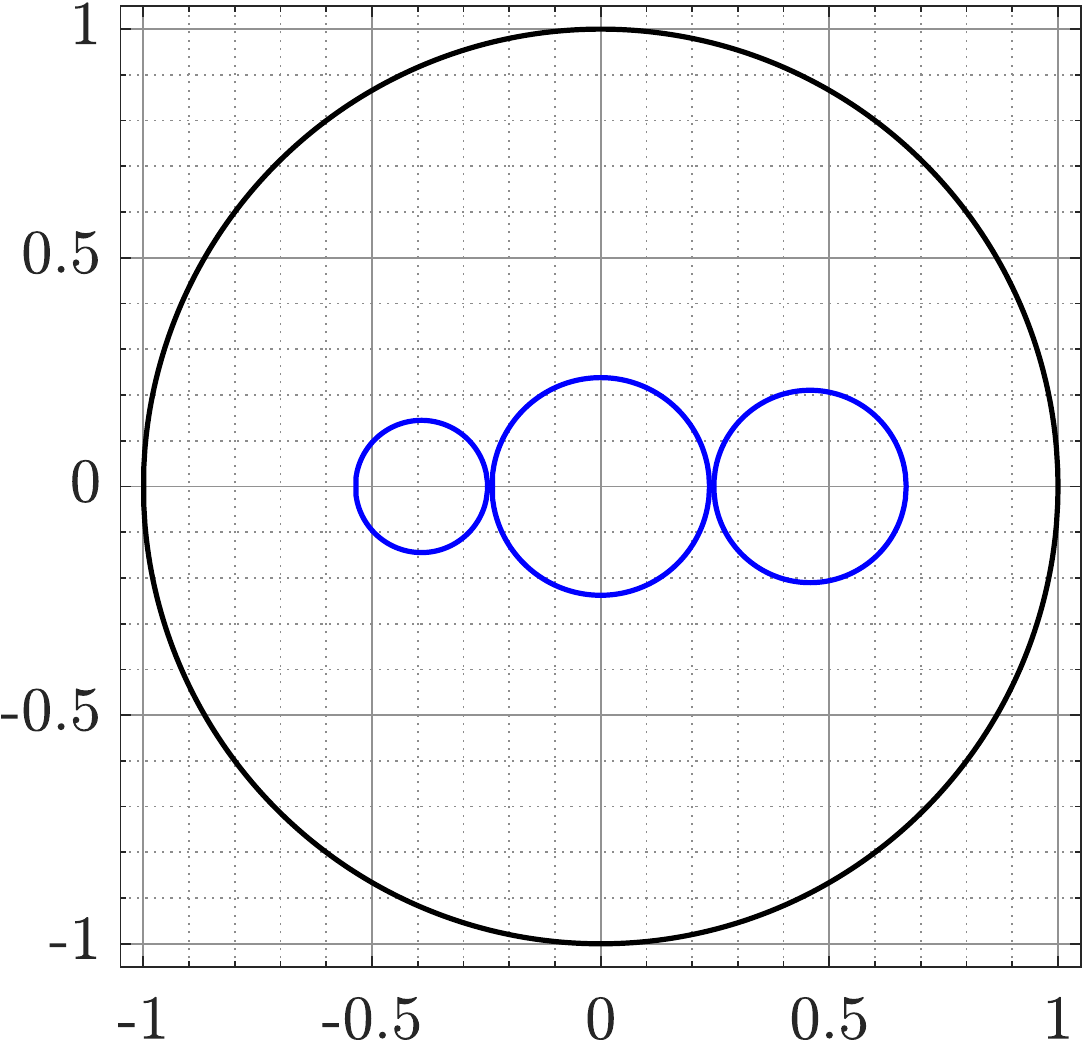}
		\includegraphics[width=0.275\textwidth]{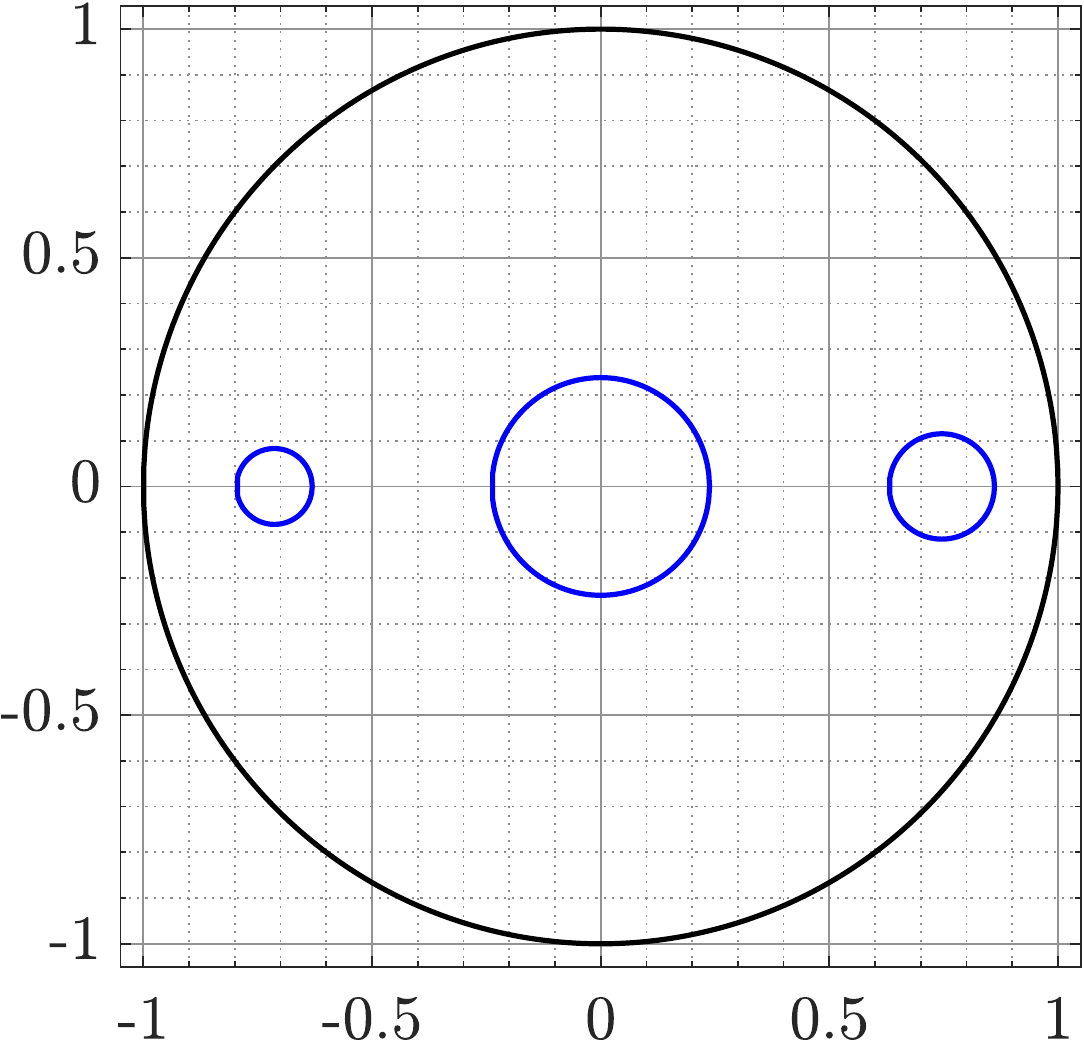}
		\includegraphics[width=0.2525\textwidth]{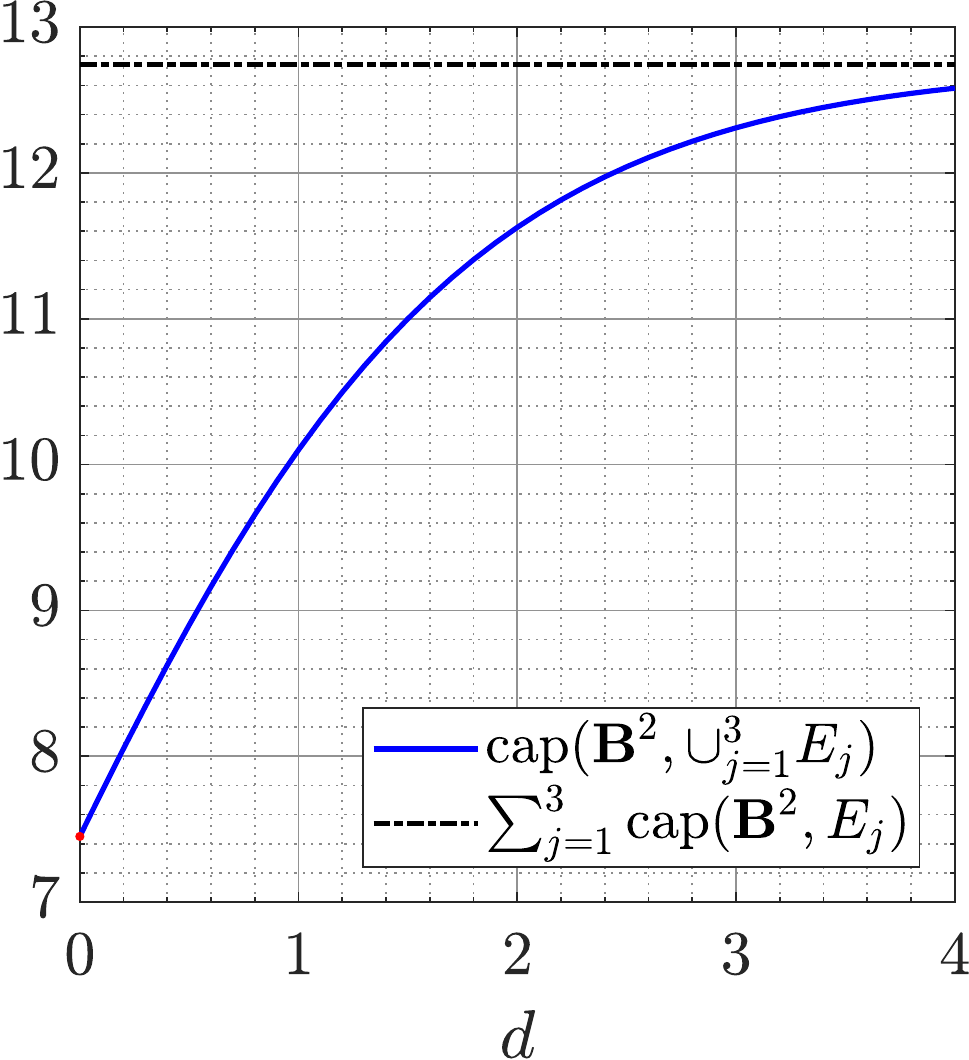}
	}\\
	\subfloat{
		\includegraphics[width=0.275\textwidth]{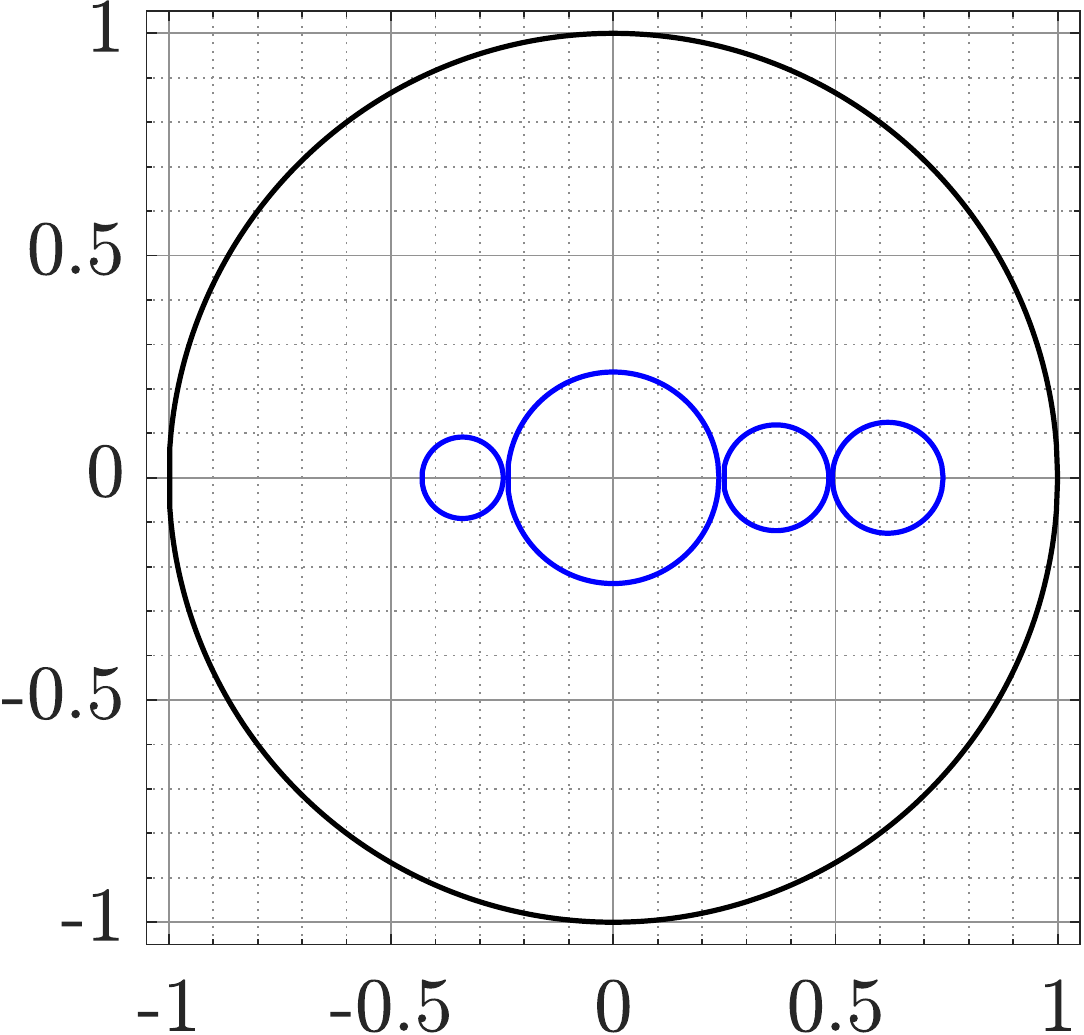}
		\includegraphics[width=0.275\textwidth]{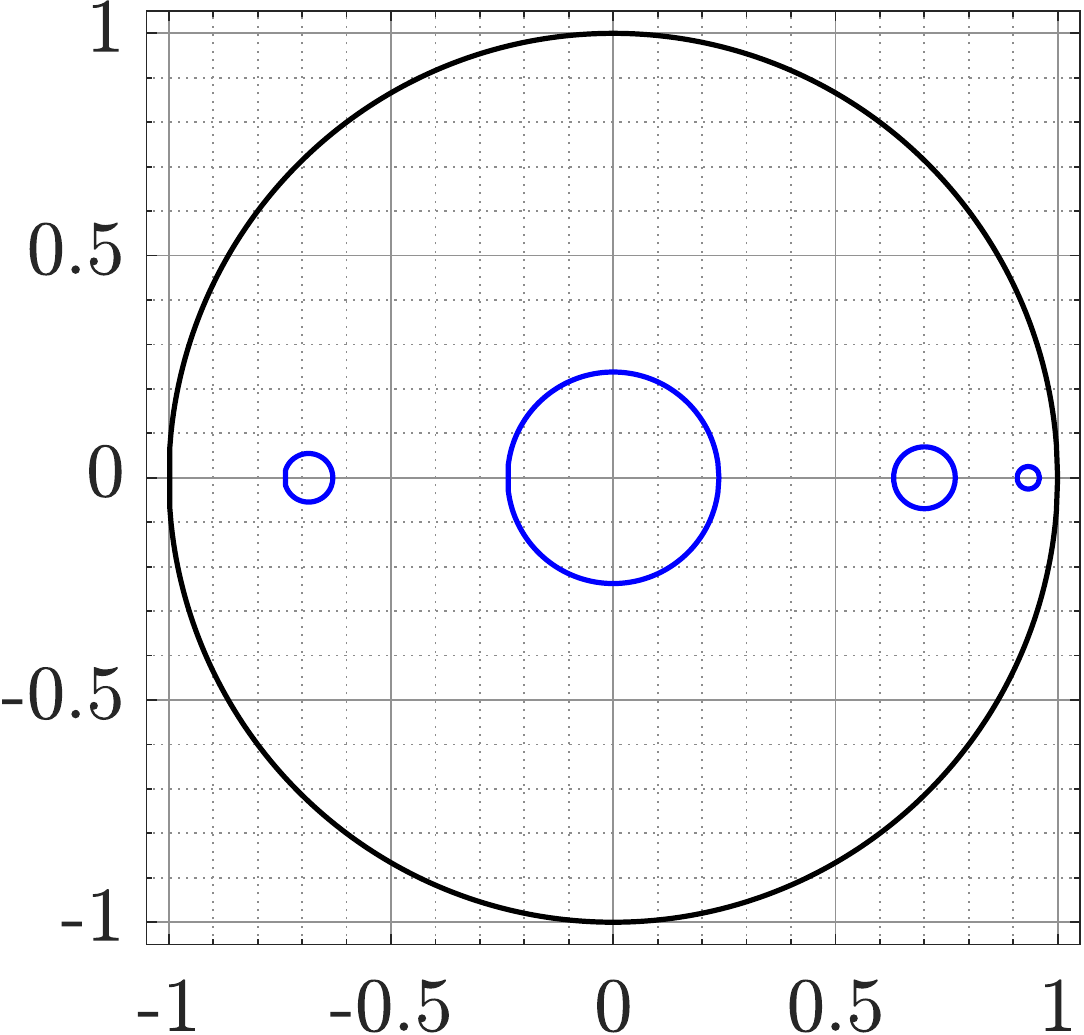}
		\includegraphics[width=0.2525\textwidth]{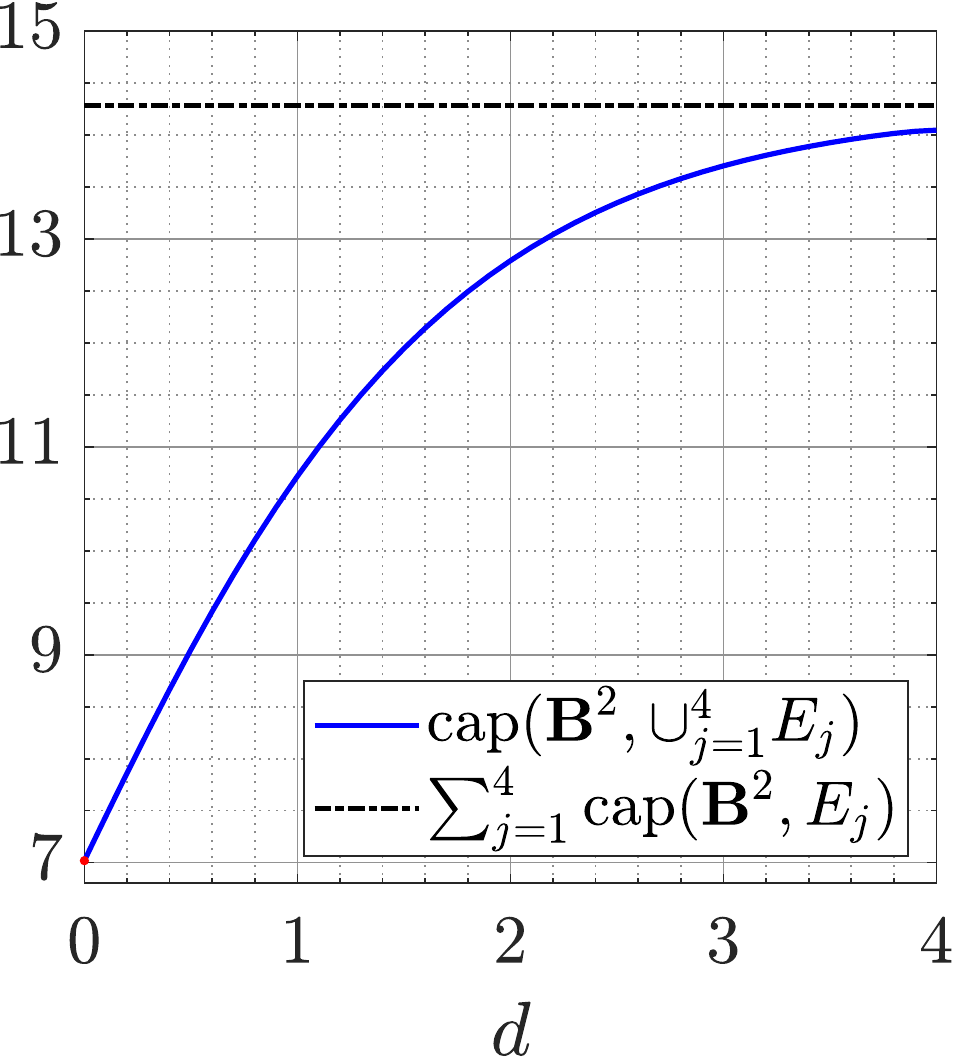}
	}
	\caption{The hyperbolic disks when the hyperbolic distance $d$ between them is $d=0.02$ (left) and $d=1$ (middle). 
		On the right, ${\rm cap}( \B, \cup_{j=1}^k E_j)$ as a function $d$.
		In the first row: $r_1=0.55d_1/2$ and $r_2=0.45d_1/2$ where $d_1=\rho_{\mathbb{B}^2}(-0.6,0.6)$.
		In the second row: $r_1=0.35d_1/2$, $r_2=0.25d_1/2$, and $r_3=0.40d_1/2$.
		In the third row: $r_1=0.35d_1/2$, $r_2=0.15d_1/2$,  $r_3=0.20d_1/2$, and $r_4=0.30d_1/2$.
	} 
	\label{fig:collineardisks}
\end{figure}

The cases ${\rm cap}(\B, \cup_{j=1}^m E_j)$ for $m=2,3,4$ over the range $0.02\le d\le 4$
are shown in Figure~\ref{fig:collineardisks}. The conjectured lower bound with $d=0$ is computed with $hp$-FEM (see the `red dot' in Figure~\ref{fig:collineardisks} (right)),
all other capacities are computed with BIE.
From Figure~\ref{fig:collineardisks} we also see that 
\[
{\rm cap}( \B,\cup_{j=1}^m E_j)\approx \sum_{j=1}^m{\rm cap}(\B,E_j),
\]
as the separation $d$ becomes large.

\begin{table}
	\centering
	\caption{Disks with collinear centers: $m$ hyperbolic disks $E_j$ with radii $r_j$ 
		and centers on the diameter $(-1,1)$ with $\sum_{j=1}^m 2 r_j=d_1=\rho_{\mathbb{B}^2}(-0.6,0.6)$.
		Conjectured lower and upper bounds of the capacity ${\rm cap}(\B, \cup E_j)$.}
	\label{tbl:collinearbounds}
	\begin{tabular}{lll}
		\hline
		$m$ & Lower & Upper \\ \hline
		$2$     & $8.515312094751020$ & $11.463763614692954$\\
		$3$     & $7.450131756754710$ & $12.744594178229441$\\
		$4$     & $7.017838565418236$ & $14.282099489357595$\\ \hline
	\end{tabular}
\end{table}

\subsubsection{Verification of results}
Let us consider the case with four disks and set $E = \cup_{j=1}^4 E_j$.
The initial position is when the disks
are contiguous, tangent to each other, and then the hyperbolic distance $d$ between 
the disks increases from $0$ to $0.3$. 
The conclusion is that the value $d=0$ yields the minimal value of the capacity of the constellation. 

\begin{figure}
	\centering
	\includegraphics[width=0.45\textwidth]{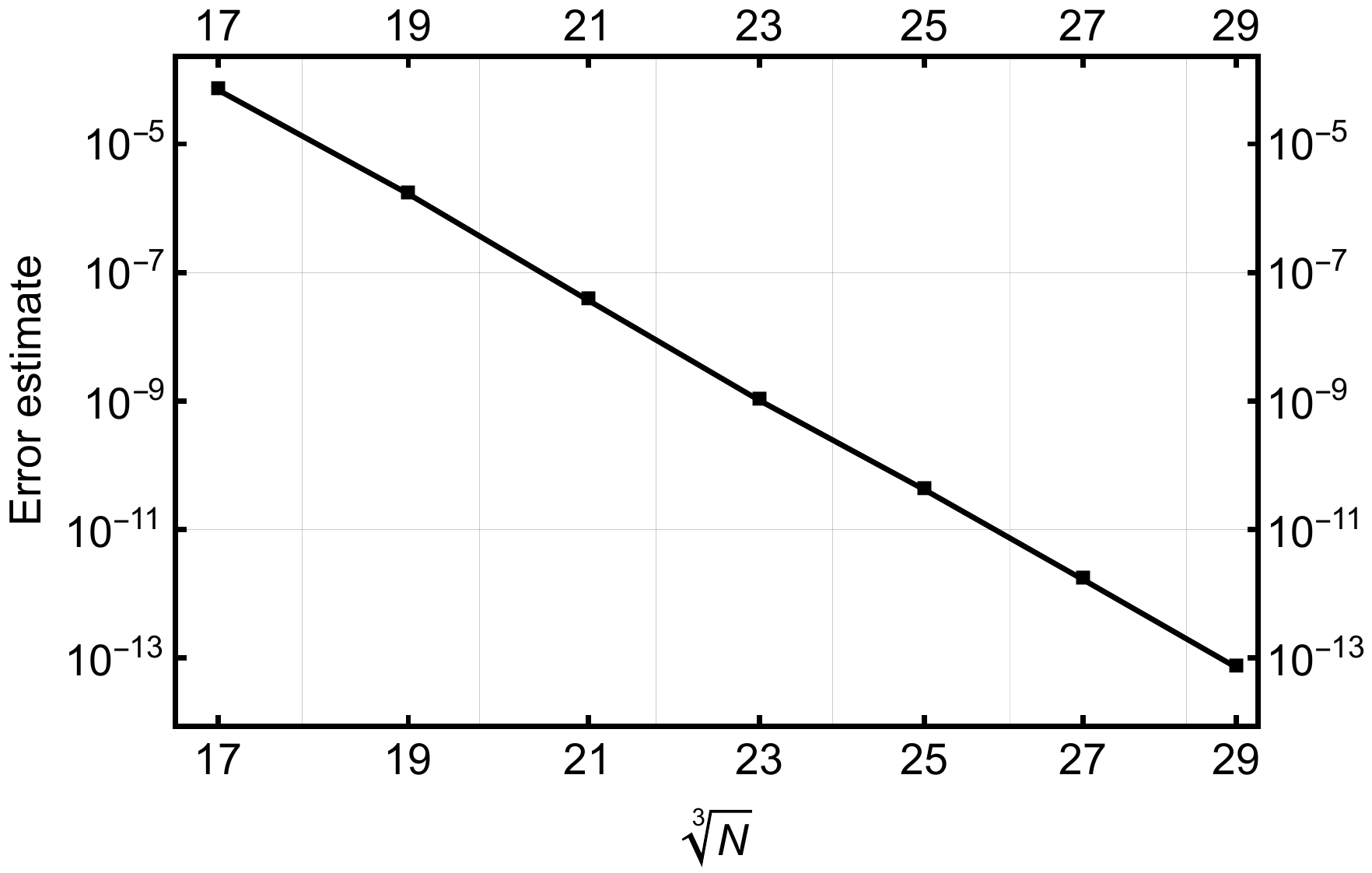}
	\hfill
	\includegraphics[width=0.4\textwidth]{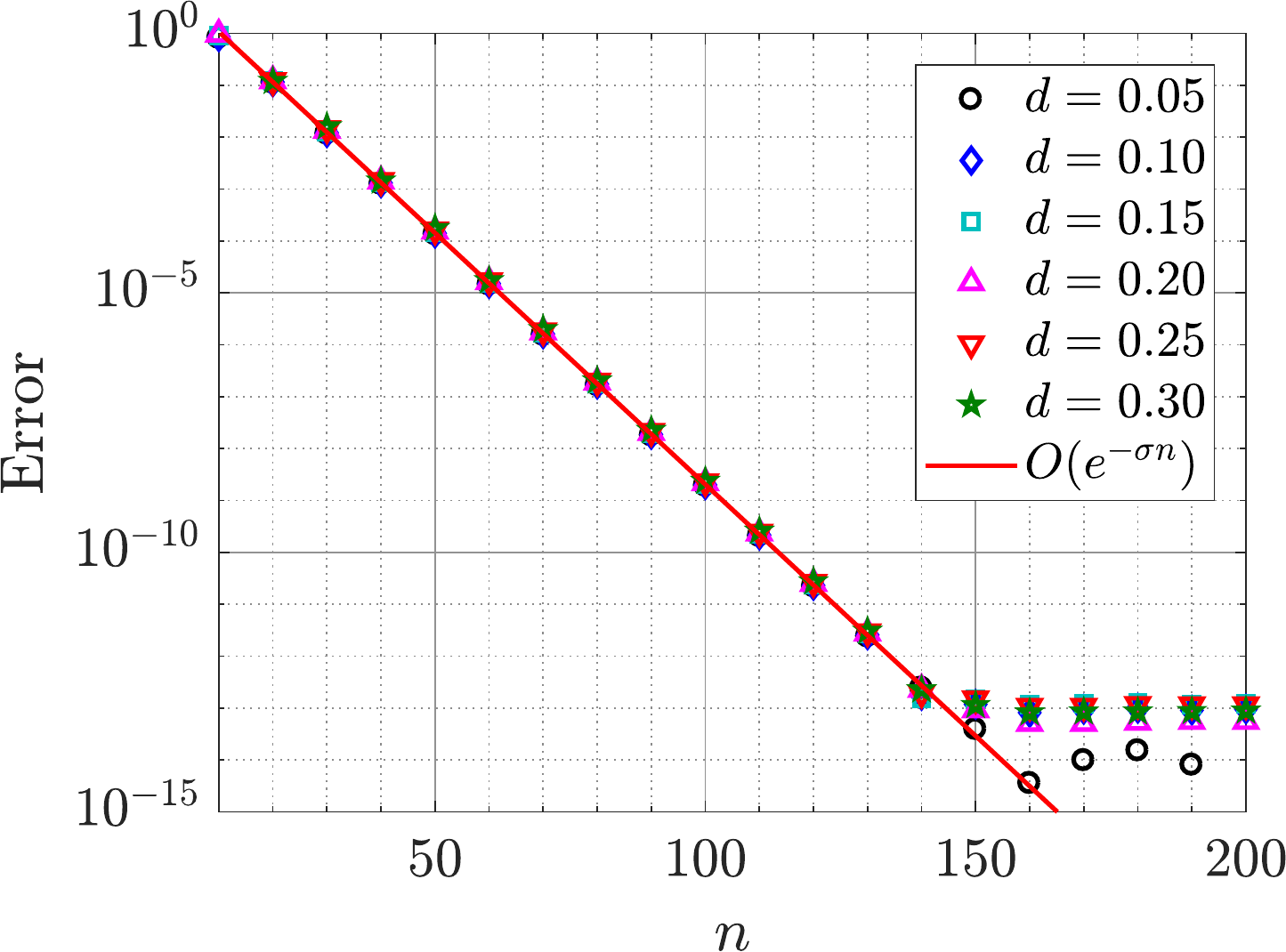}
	\caption{The error for the constellation of four disks in Table~\ref{horizTab}. 
		Left: The $hp$-FEM error estimate as a function of $\sqrt[3]{N}$, where $N$ is the number of d.o.f. (logplot) for four disk configuration with contacts ($d=0$).
		The observed constant or the slope of the graph $= 37.1$. 
		Right: The errors in the computed values of ${\rm cap}({\mathbb B}^2,E)$ using the BIE method as functions of $n$, for $\alpha=0.8\i$ where $\sigma=-\log|\alpha|\approx0.223$.
	}\label{fig:fourmesherror}
\end{figure}

The values of the capacity ${\rm cap}({\mathbb B}^2,E)$ in Table~\ref{horizTab} have been computed using both methods, the FEM and the BIE method. 
For the BIE, we use $n=2^7$ and $\alpha=0.8\i$.  
Table~\ref{horizTab} shows the absolute differences between the computed values which indicates a good agreement between the two methods. As in~\cite{hnv}, the values computed using the FEM will be considered as reference values and used to estimate the error in the values computed by the BIE method for several values of $n$. 
The BIE method cannot be used for $d=0$. The error for $d=0.05,0.1,\ldots,0.3$ is presented in Figure~\ref{fig:fourmesherror} (right) which illustrates the exponential convergence with order of convergence $O(e^{-\sigma n})$ where $\sigma=-\log|\alpha|\approx0.223$. 
Numerical experiments (not presented here) with other values of $\alpha$ indicate that the order of convergence depends on $\alpha$ as well as the centers $z_1,\ldots,z_m$ and the radii $r_1,\ldots,r_m$ of the inner circles. A detailed analysis of the order of convergence for the above BIE method is a subject of future work.

\begin{table}
	\centering
	\caption{Computed values of ${\rm cap} (\mathbb{B}^2,E)$ when $m=4$ for a constellation with disk radii 
		(from left to right) $r_1=0.15d_1/2$, $r_2=0.35d_1/2$, $r_3=0.20d_1/2$, and $r_4=d_1/2-(r_1+r_2+r_3)$ 
		where $d_1=\rho_{\mathbb{B}^2}(-0.6,0.6)$. 
		The centers on the diameter $(-1,1)$ as a function of the hyperbolic distance $d$ between disks, i.e., 
		$c_1=-\th((r_2+r_1+d)/2)$, $c_2=0$, $c_3= \th((r_2+r_3+d)/2)$, $c_4=-\th((r_2+2r_3+r_4+2d)/2)$. 
	}
	\label{horizTab}
	\begin{tabular}{llll}
		\hline
		$ d $ &  FEM                    & BIE          & Agreement\\
		\hline
		$0.00$ & $7.017838565413617$  & ---                  & ---\\
		$0.05$ & $7.230698262298420$  & $7.230698262298405$  & $1.51\times 10^{-14}$\\
		$0.10$ & $7.442082617728579$  & $7.442082617728490$  & $8.88\times 10^{-14}$\\ 
		$0.15$ & $7.651760366696882$  & $7.651760366696745$  & $1.37\times 10^{-13}$\\
		$0.20$ & $7.859490827905997$  & $7.859490827905935$  & $6.22\times 10^{-14}$\\
		$0.25$ & $8.064996233395842$  & $8.064996233395734$  & $1.08\times 10^{-13}$\\
		$0.30$ & $8.267972932727597$  & $8.267972932727497$  & $9.95\times 10^{-14}$\\
		\hline
	\end{tabular}
\end{table}

\subsection{{Four disks: Permutation of contiguous disks}} We consider next two cases where all
the disks of the constellation have fixed hyperbolic radii $A > B > C > D > 0$ but their relative
ordering is not constrained other than that  
each disk
is tangent to at least one other disk of the constellation and
their hyperbolic centers lie (a) either on the diameter $(-1,1)$ or 
(b) on the circle $\{z: |z|=1/2\}.$

Now the question is what is the effect of the permutation of the disks on the capacity.
There are 24 permutations with 12 different capacities due to symmetry.
For every realisation, the radii are denoted by $r_j$ from left to right and
the constellations are denoted by $E_D$ and $E_C$, respectively.
For $E_D$ we set $(A,B,C,D) = (1/2,2/5,1/4,1/5)$, and
for $E_C$ slightly perturbed $(A,B,C,D) = (1/2,1/3,1/4,1/5)$.
The results are collected in Table~\ref{tbl:permTable} and Figure~\ref{fig:appetizer} 
shows the observed extremal permutations.
Interestingly, the resulting capacities have exactly the same
dependence on the relative sizes of the radii.

\begin{table}
	\centering
	\caption{Permutations of contiguous constellations. $E_D$ with centers on the segment $(-1,1)$
		and $(A,B,C,D) = (1/2,2/5,1/4,1/5)$.
		$E_C$ with centers on the circle $\{z: |z|=1/2\}$
		and $(A,B,C,D) = (1/2,1/3,1/4,1/5)$.}\label{tbl:permTable}
	\begin{tabular}{lllllll}
		\hline
		Case  & $r_1$ &   $r_2$ & $r_3$ & $r_4$ & $\capa C(E_D)$ & $\capa C(E_C)$ \\ 
		\hline
		1  & D & B & A & C & 6.781488018927628 & 6.451424010111881 \\
		2  & D & A & B & C & 6.788910565780309 & 6.455800945561348 \\
		3  & D & C & A & B & 6.843774515059010 & 6.475070264106950 \\
		4  & C & D & A & B & 6.882473842468833 & 6.485425869048534 \\
		5  & A & B & C & D & 6.890544149275032 & 6.496389476635198 \\
		6  & B & C & A & D & 6.897202225461369 & 6.500210100051595 \\
		7  & C & A & D & B & 6.919626376828870 & 6.520197932005349 \\
		8  & A & B & D & C & 6.928074481413122 & 6.523073055329720 \\
		9  & A & C & B & D & 6.932436180755356 & 6.542555705939787 \\
		10 & C & B & D & A & 6.962814943144452 & 6.542981227003898\\
		11 & A & C & D & B & 7.053764008325471 & 6.575258877036491\\
		12 & A & D & C & B & 7.055565195334228 & 6.576332514877286\\
		\hline
	\end{tabular}
\end{table}

\begin{table}
	\centering
	\caption{Hyperbolic radii used in Figure~\ref{fig:moving}.}\label{tbl:hypRad}
	\begin{tabular}{lllll}\hline
		Case & $r_1$ & $r_2$ & $r_3$ & $r_4$ \\ \hline
		1 & 0.4 & 0.2 & 0.5 & 0.25 \\
		2 & 0.2 & 0.5 & 0.3 & 0.5 \\
		3 & 0.5 & 0.5 & 0.5 & 0.2 \\
		4 & 0.2 & 0.7 & 0.4 & 0.1 \\ \hline
	\end{tabular}
\end{table}

\begin{figure}
	\centering
	\subfloat{
		\includegraphics[width=0.22\textwidth]{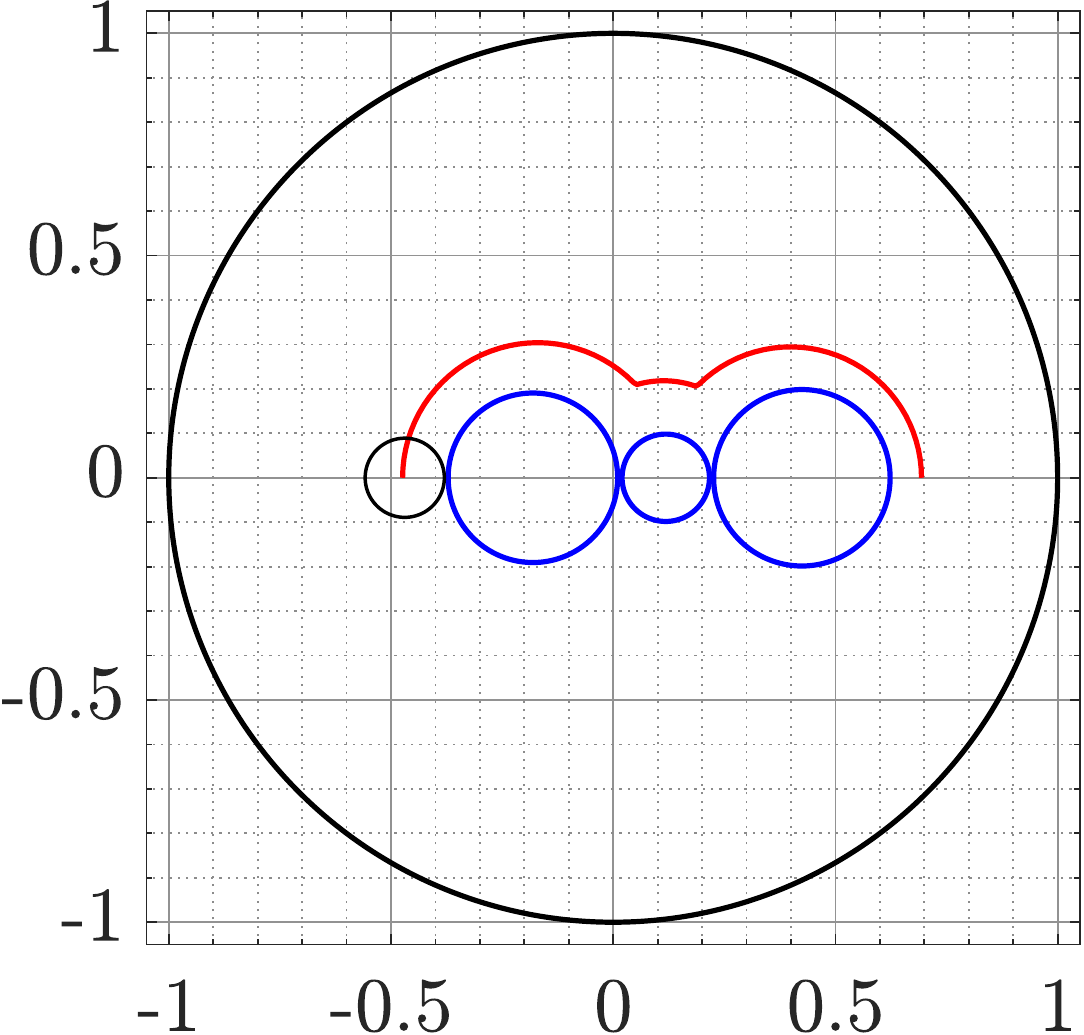}
	}
	\subfloat{
		\includegraphics[width=0.22\textwidth]{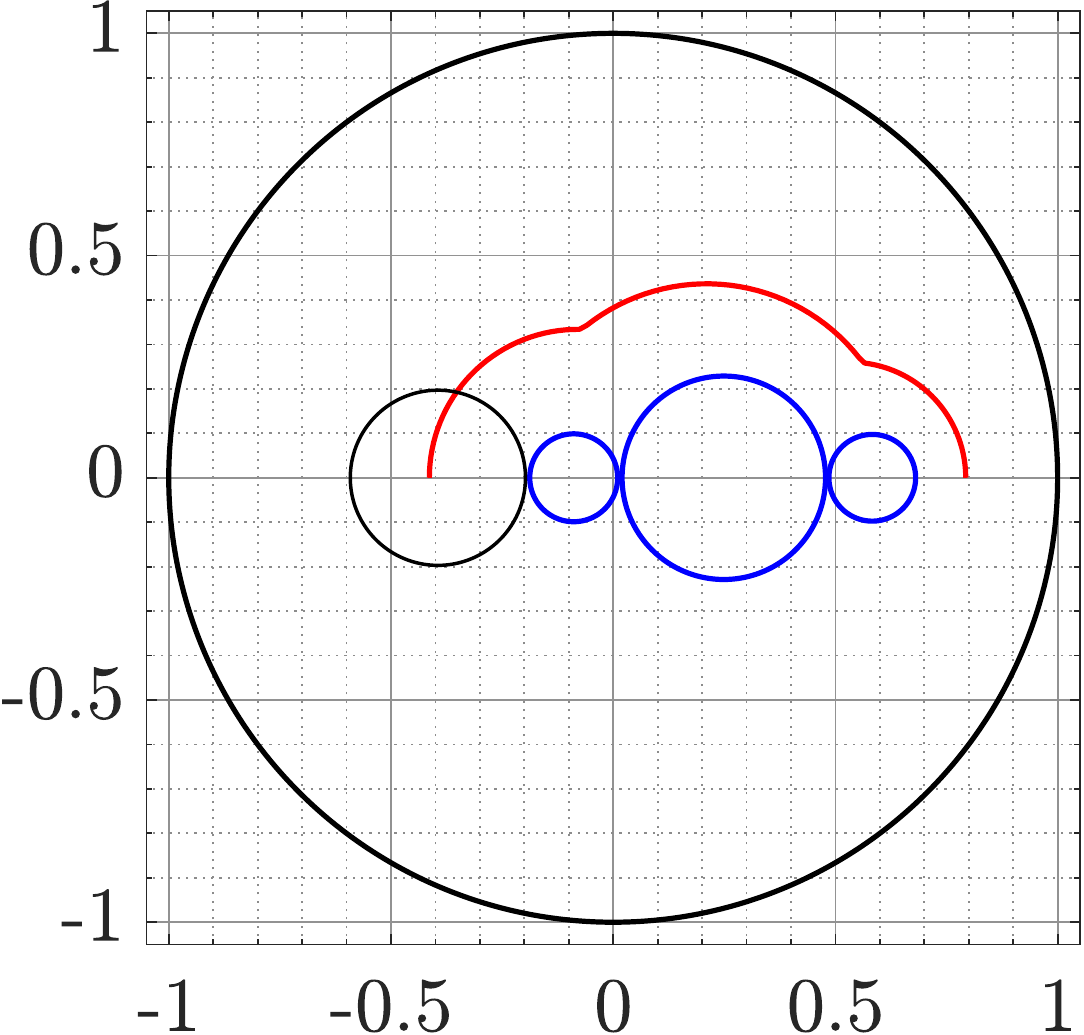}
	}
	\subfloat{
		\includegraphics[width=0.22\textwidth]{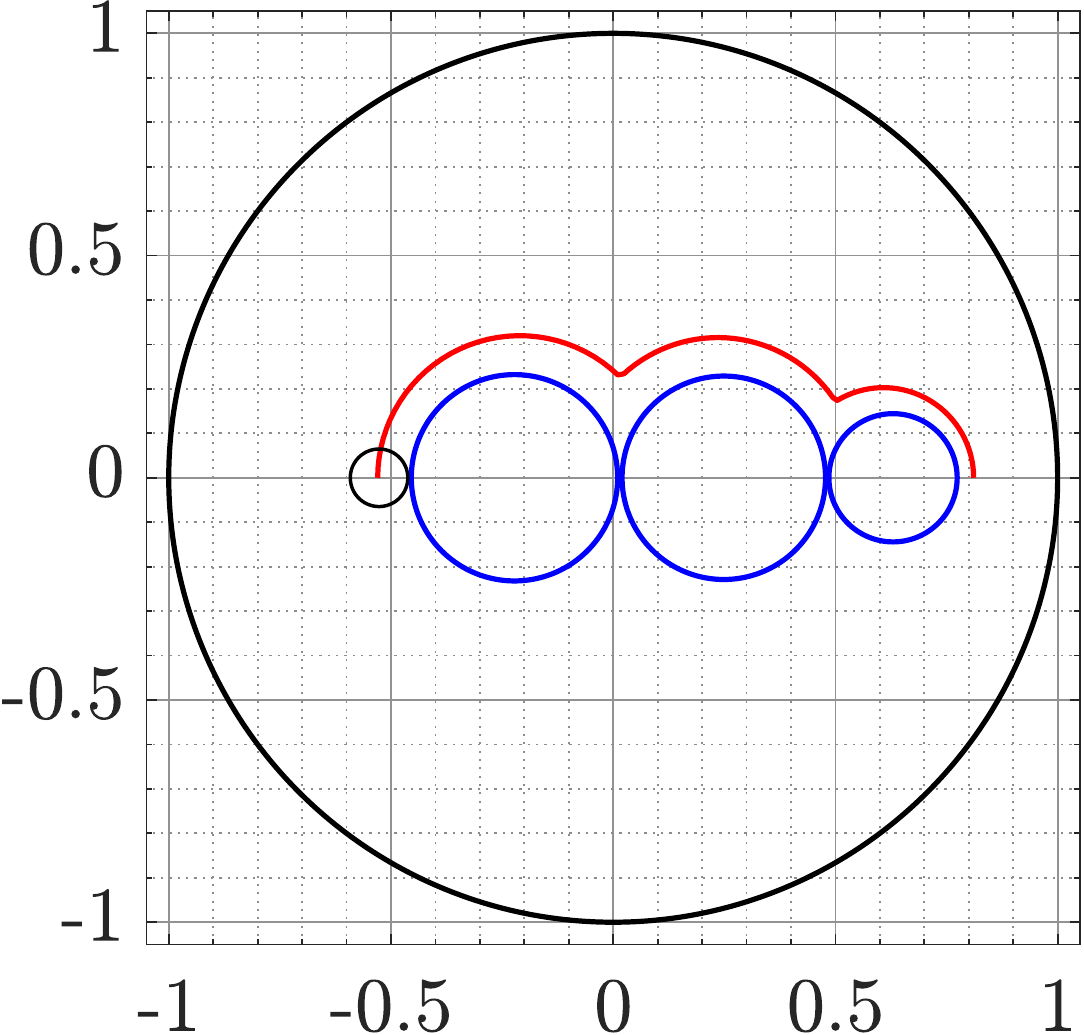}
	}
	\subfloat{
		\includegraphics[width=0.22\textwidth]{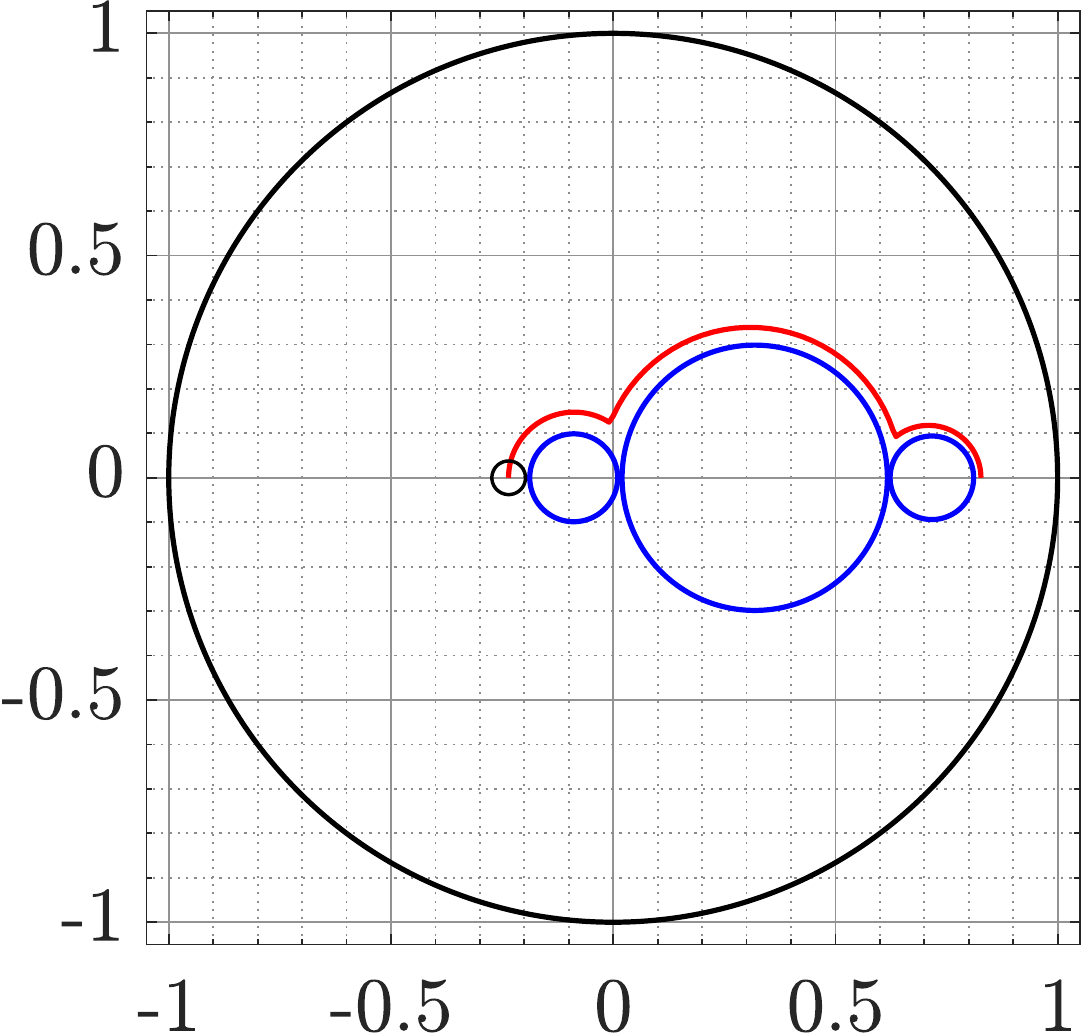}
	}\\
	\subfloat{
		\includegraphics[width=0.22\textwidth]{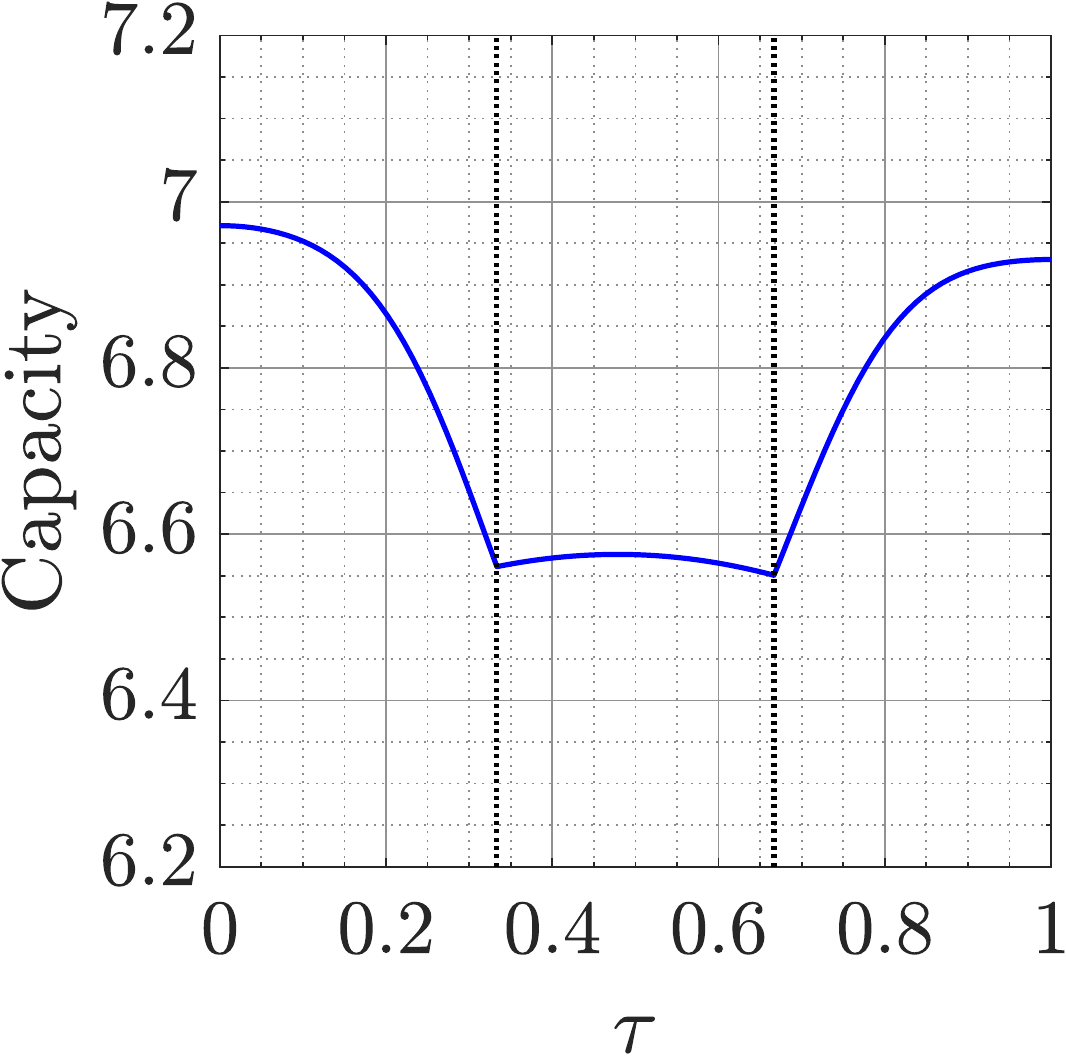}
	}
	\subfloat{
		\includegraphics[width=0.22\textwidth]{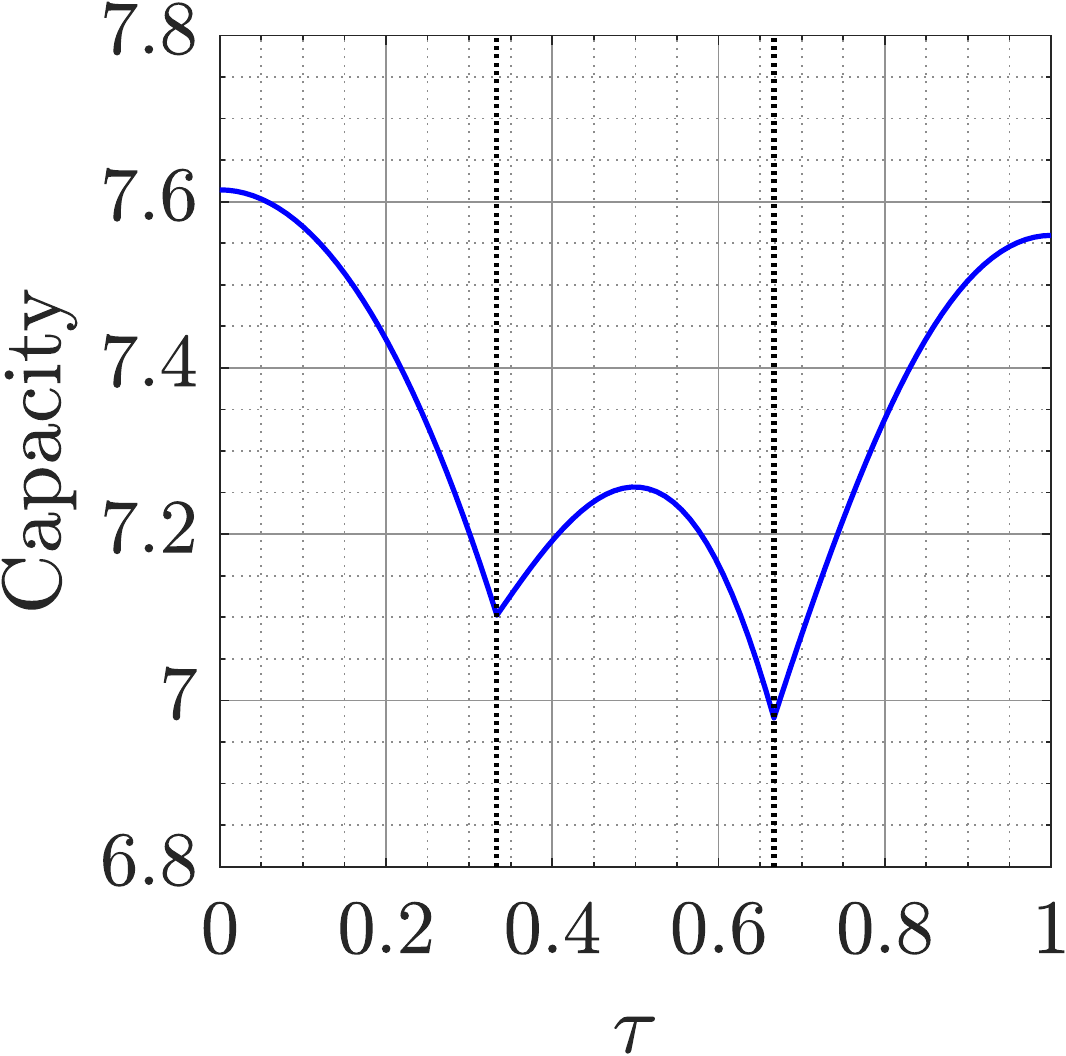}
	}
	\subfloat{
		\includegraphics[width=0.22\textwidth]{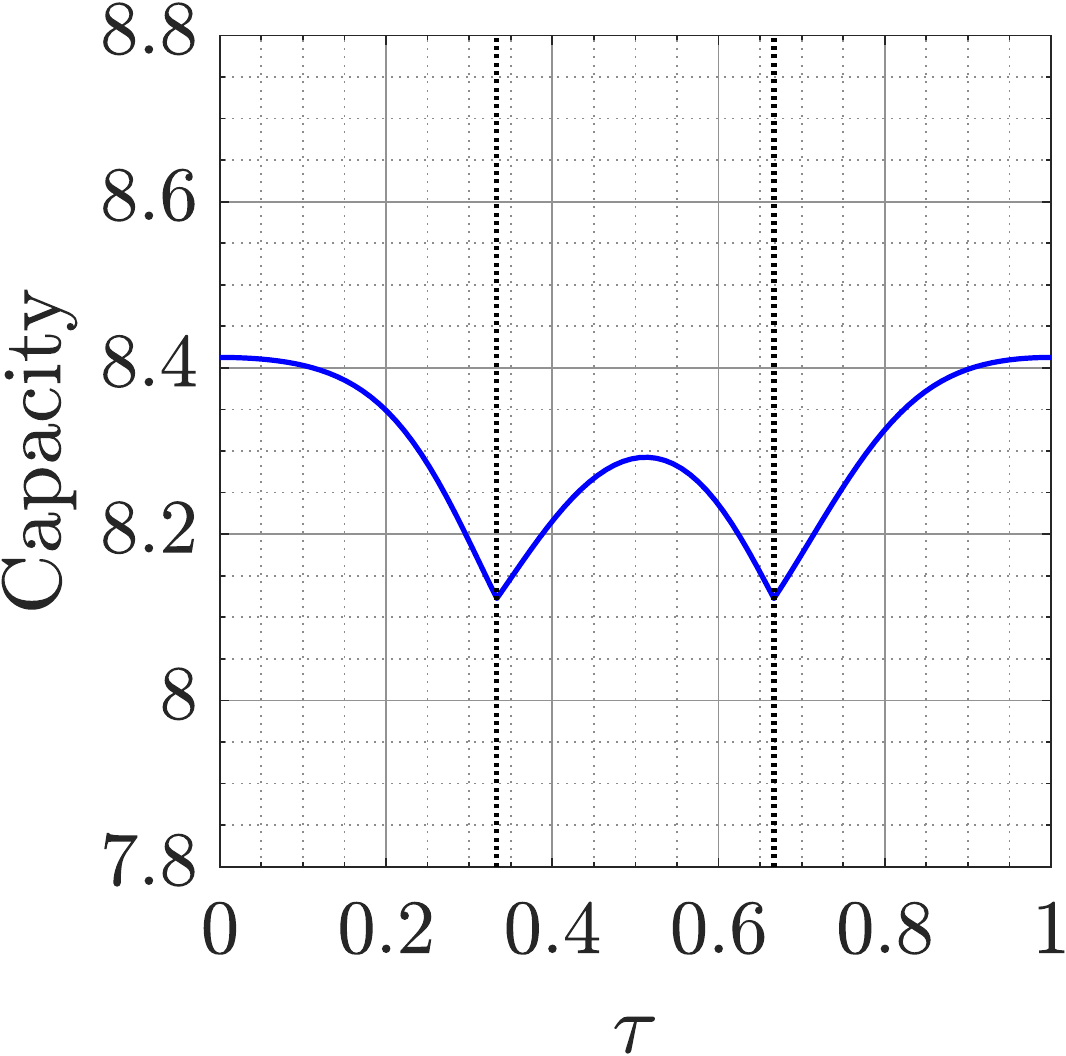}
	}
	\subfloat{
		\includegraphics[width=0.22\textwidth]{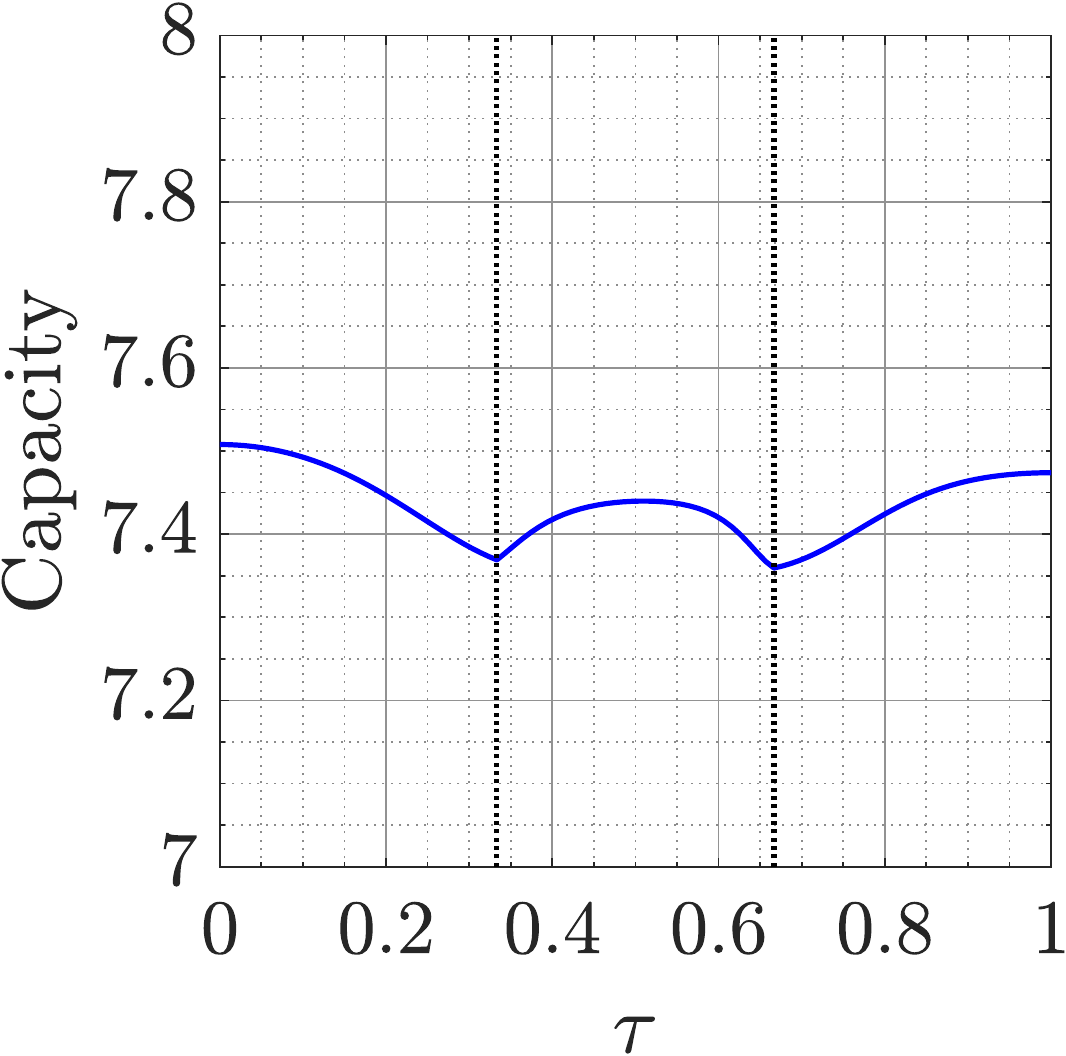}
	}
	\caption{Three immobile disks, one rolling disk. Cases 1 to 4 from left to right.
		Dependence of the capacity on the relative location of the rolling disk.
		The hyperbolic center $z_4$ of the moving disk is on the red curve shown in the figure.}\label{fig:moving}
\end{figure}

\subsection{{Three immobile disks, one rolling disk}} 
In the final experiment of the section we study the situation when one disk is free to roll on the remaining three contiguous immobile disks, centers on the diameter $(-1,1)$ and tangent to each other. 
The route of the mobile disk is parametrized with a parameter $\tau \in[0,1]$ where the values $0$ and $1$ are for the case when also the mobile disk has its center on  the diameter $(-1,1)$ and the values $1/3$ and $2/3$ correspond to the intermediate points on the route when the rolling disk is tangent to two immobile disks. 
Depending on the radii, it might also happen that there is only one such point. 
In Figure~\ref{fig:moving} below we see that for the values $1/3$ and $2/3$ the capacity of the constellation attains a local minimum.
The numerical results for this example are computed using the BIE method. So, instead of assuming that the disks are touching each other, we assume that the disks are close to each other such that the hyperbolic distance between them is $d=0.02$.
In all cases the hyperbolic centers of the three fixed disks are $z_1=-\th((r_1-d)/2)$, $z_2=\th((r_2+2d)/2)$, and $z_3=\th((r_3+2r_2+3d)/2)$. 
The hyperbolic center $z_4$ of the moving disk is on the red curve shown in the figure. 
The observed results are summarized in the second row of Figure~\ref{fig:moving}.

\section{Minimizing Capacity: Optimization under Free Mobility}\label{fed}
In this section we consider a series of experiments, where some disks are given
fixed positions but the others are free to move within constraints.
The constraints can restrict the admissible configurations to specific regions.
In the most general case, the only constraint is that the disks should not overlap.
In all simulations it is assumed that the disks have a minimal separation $\delta > 0$.
In those cases where the disks touch, that is, $\delta = 0$, only $hp$-FEM results 
are reported.

\subsection{Three fixed disks. One freely moving disk}\label{sc:3f1m} 

\begin{figure} %
	\centering
	\subfloat[$r=0.6$]{
		\includegraphics[width=0.31\textwidth]{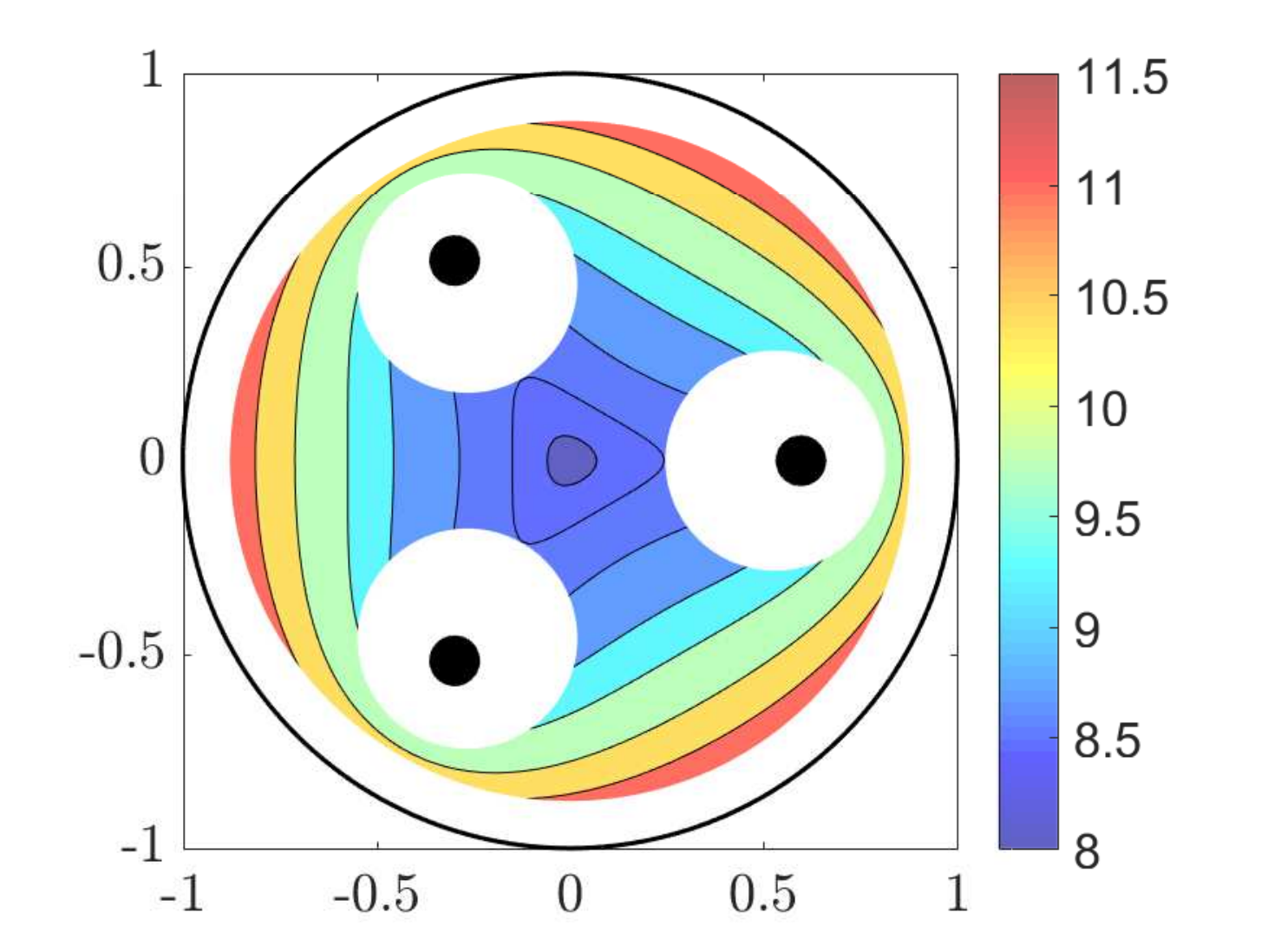}
	}
	\subfloat[$r=0.5$]{
		\includegraphics[width=0.31\textwidth]{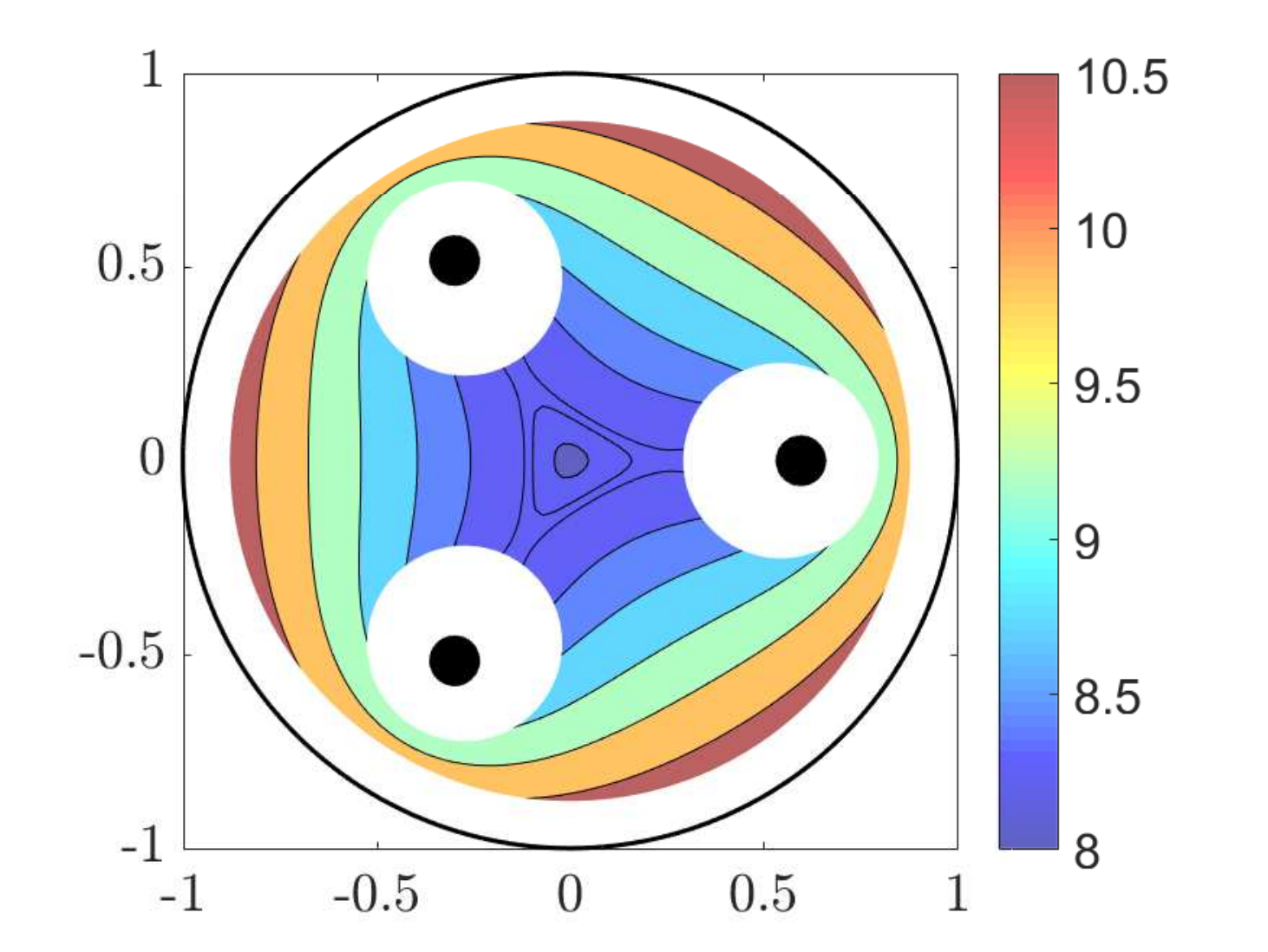}
	}
	\subfloat[$r=0.4$]{
		\includegraphics[width=0.31\textwidth]{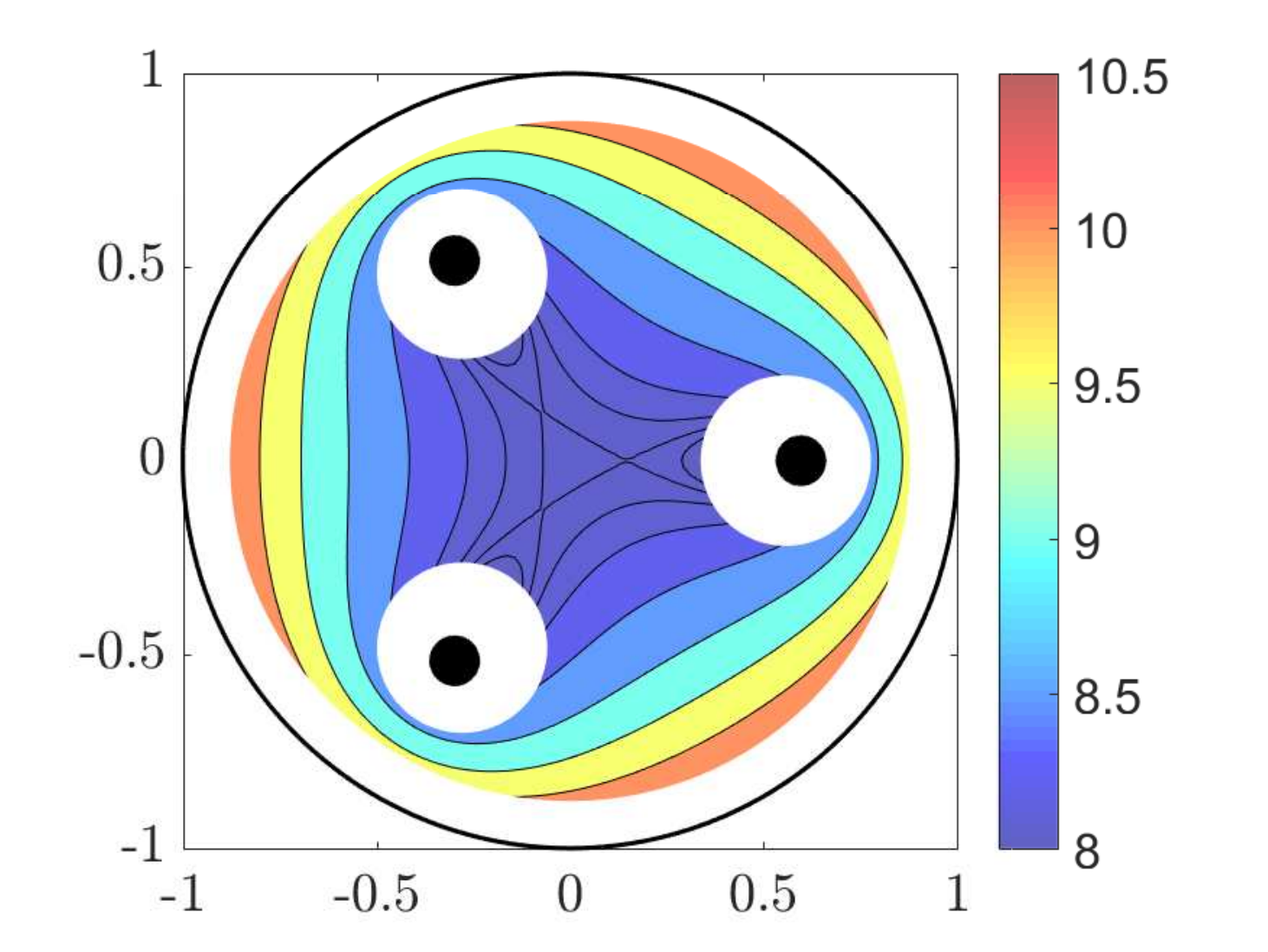}
	}\\
	\subfloat[$r=0.3$]{
		\includegraphics[width=0.31\textwidth]{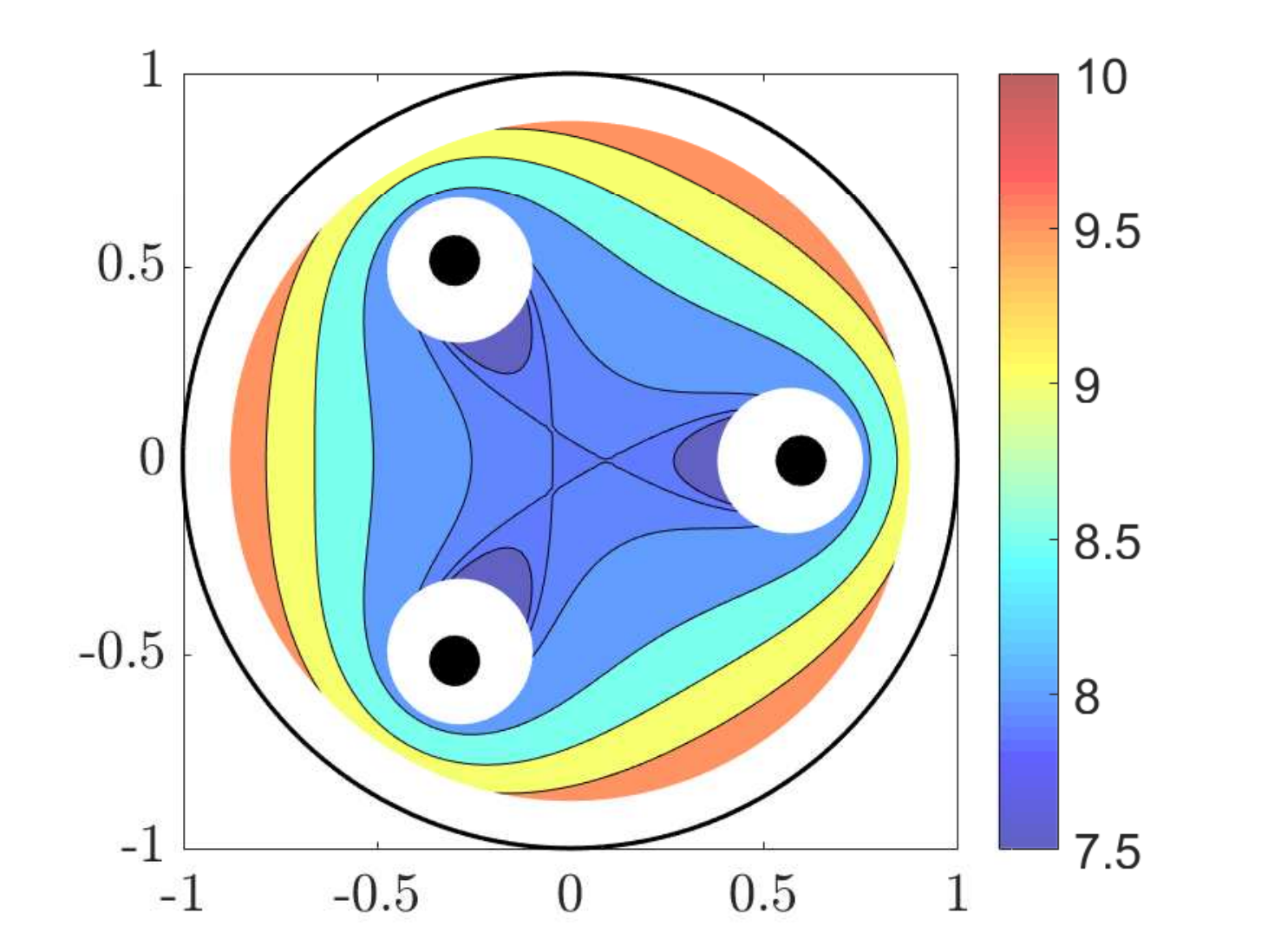}
	}
	\subfloat[$r=0.2$]{
		\includegraphics[width=0.31\textwidth]{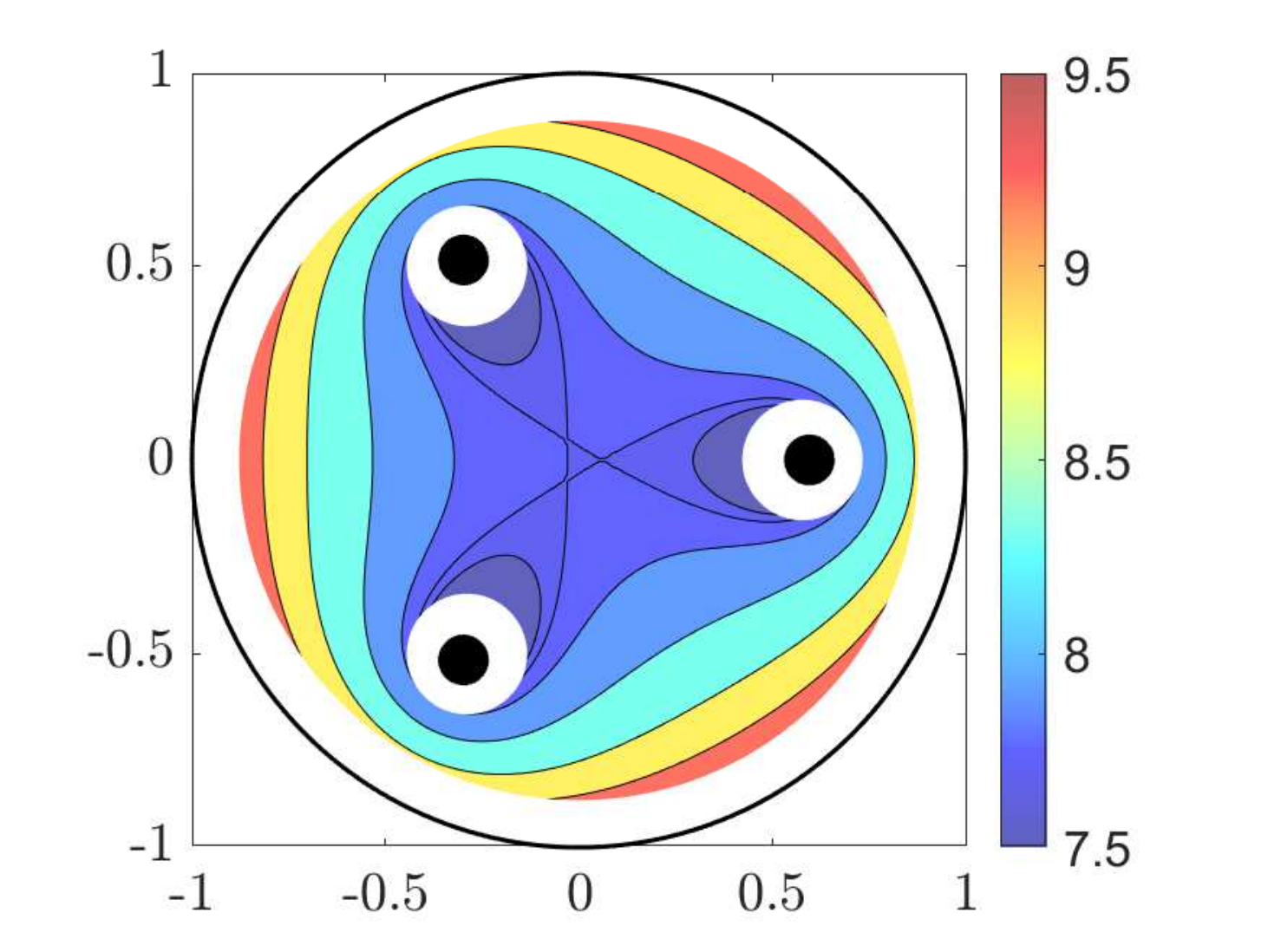}
	}
	\subfloat[$r=0.1$]{
		\includegraphics[width=0.31\textwidth]{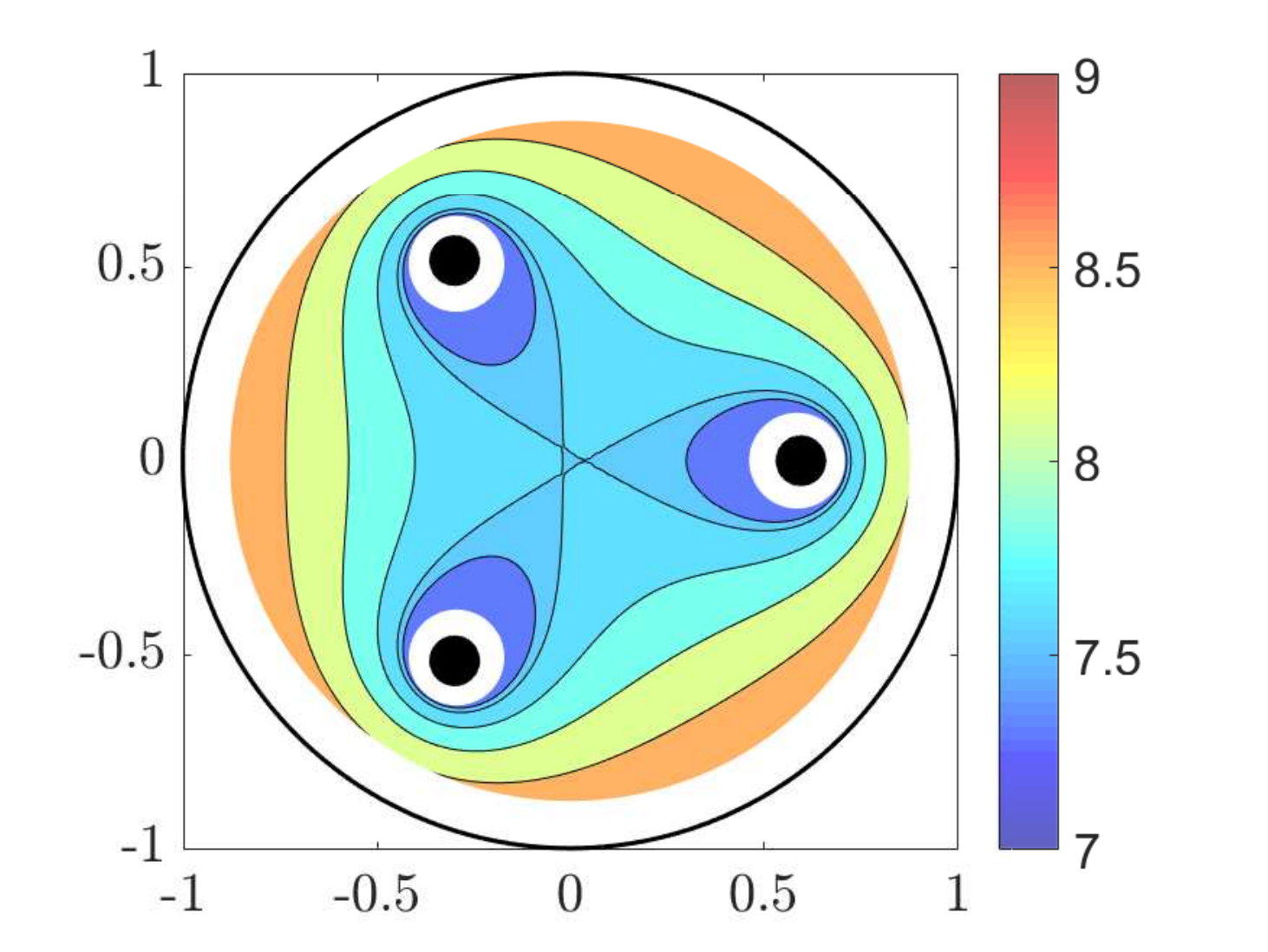}
	}
	\caption{The level curves of the capacity $u(x,y)$ of a constellation of four disks as a function of the center $z=x+\i y$ of the fourth disk.
		Three hyperbolic disks with equal hyperbolic radii $= 0.2$ are at fixed locations, whereas the fourth one with a given radius $r$ is free to move (the mobile fourth disk is not shown). The number of local minima depends on the radius of the fourth disk.
	} 
	\label{fig:4disk-3fixed}
\end{figure}

Consider three hyperbolic disks with equal hyperbolic radii $= 0.2$, and whose centers are at
$
0.5 e^{2(k-1)\pi\i/3}, \quad k=1,2,3.
$
We consider a fourth hyperbolic disk whose hyperbolic radius is $r$ and its hyperbolic center is $z=x+\i y$ such that the four disks are non-overlapping. 
Let a function $u(x,y)$ be defined by
$
u(x,y)= {\rm cap}( \B, E),
$
where $E$ is the union of the four disks. The level curves of the function $u(x,y)$ for six cases of $r$ are given in Figure~\ref{fig:4disk-3fixed}. 
Notice, that the locations of the local minima depend on the chosen radius $r$
of the free disk. Due to symmetry, there is a local minimum at the origin in every case.
The results suggest that there exists a critical radius $r_c$ such that
the global minimum is found at the origin for all sufficiently large $r$, that is, $r > r_c$,
but next to one of the fixed disks for $r < r_c$.
The interior-point method is guaranteed to converge to one of the local minima, and
therefore for all $r$ a local minimum may be attained when the mobile disk is centered at the origin.
%

\subsection{One fixed disk. Two moving disks on a circle} 
\begin{figure} %
	\centering
	\subfloat{
		\includegraphics[width=0.275\textwidth]{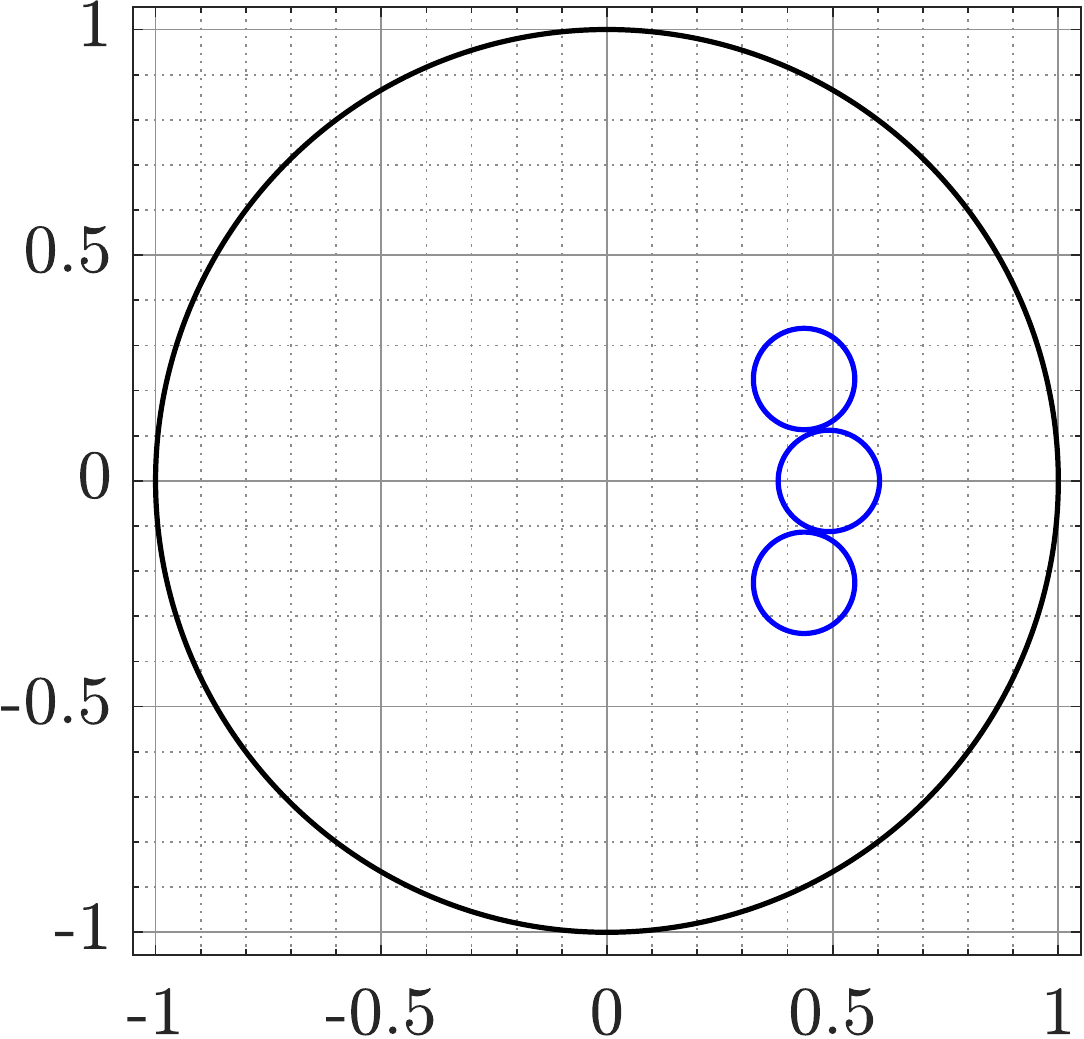}
	}
	\subfloat{
		\includegraphics[width=0.275\textwidth]{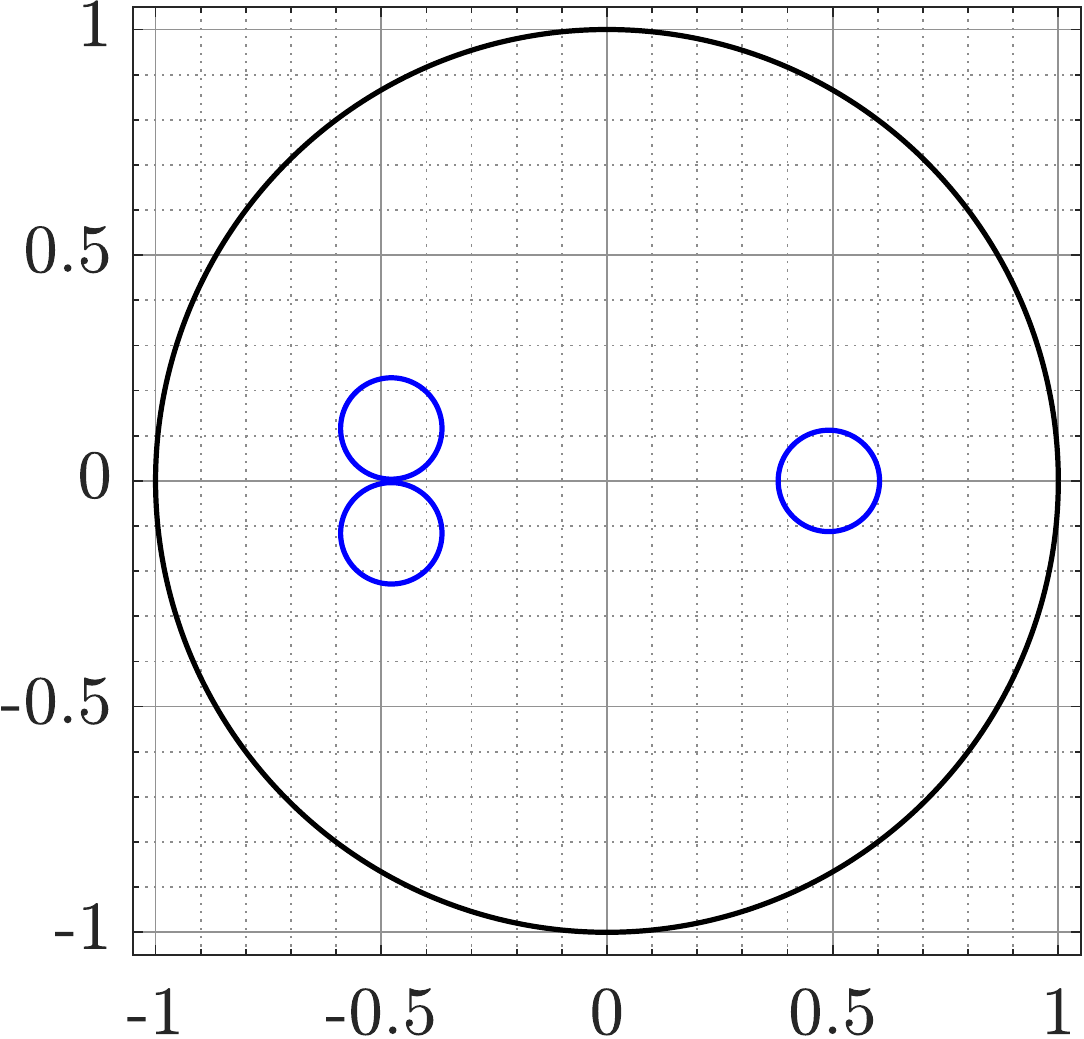}
	}
	\subfloat{
		\includegraphics[width=0.275\textwidth]{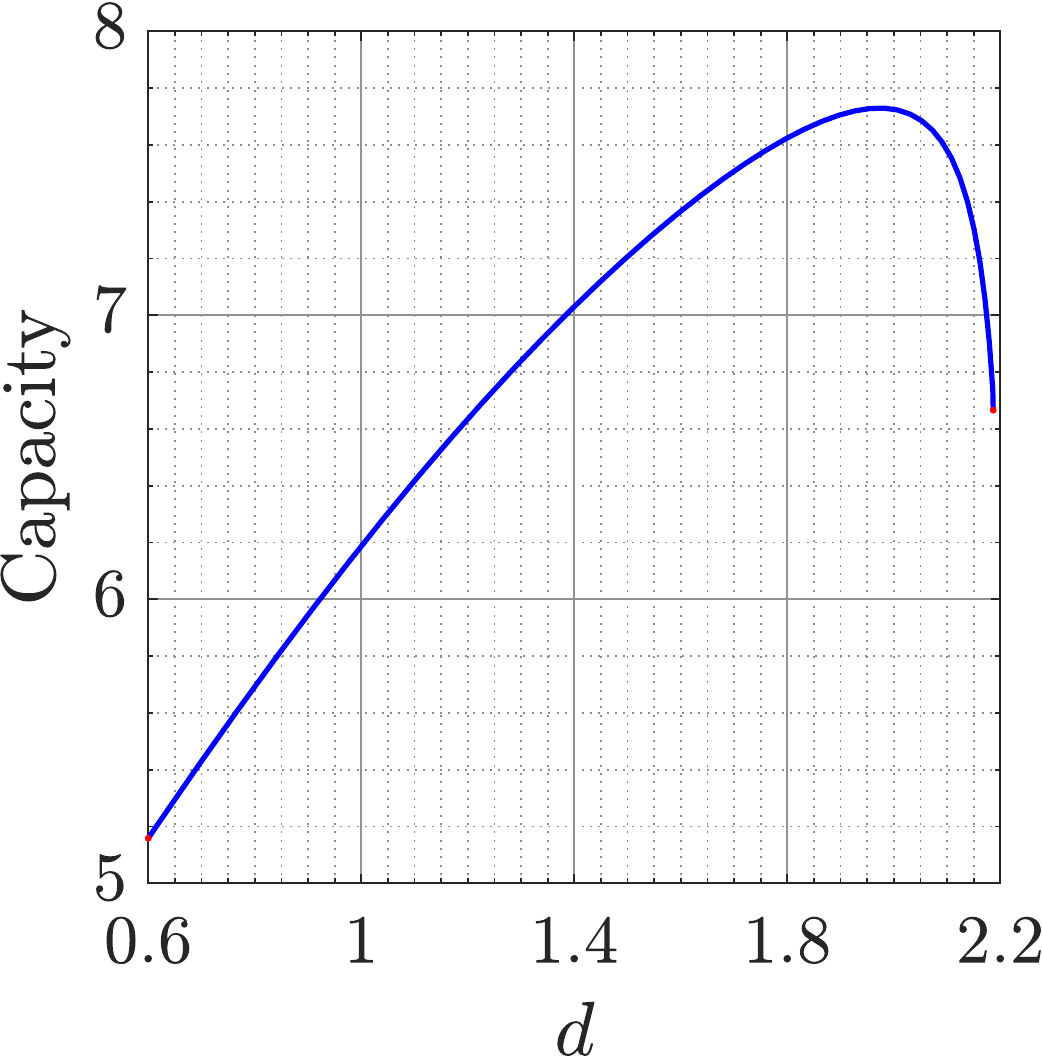}
	}
	\caption{Three disks with equal hyperbolic radius $= 0.3$ on the circle $|z|=0.5$.
		One disk is fixed on the positive real line and the other two move symmetrically
		on upper- and lower-half planes, respectively. 
		The left and middle figures illustrate the minimal $d_{\rm min}$ and maximal $d_{\rm max}$ values of the hyperbolic distance $d$ between the hyperbolic centers of the disk on the real line and the disk on the upper-half plane. The right figure
		shows the capacity for the range between these extreme valued  $d_{\rm min}\le d \le d_{\rm max}$.} 
	\label{fig:3diskslimit}
\end{figure}

Let us next consider three disks $D_1,D_2,D_3$ with equal hyperbolic radii
$r=0.3$. The centers of these three disks are placed on the circle $|z|=0.5$.
We assume that the disk $D_1$ is fixed with center on the positive real
line, $D_2$ is in the upper-half plane and $D_3$ is in the lower-half
plane. Starting when the three disks are touching each others (see
Figure~\ref{fig:3diskslimit} (left)), these disks start moving
away from each other such that the hyperbolic distance $d$ between the
hyperbolic centers of $D_1$ and $D_2$ is the same as for $D_1$ and
$D_3$. When all these disks are touching each other, $d=2r$. The maximum
value $d_{\rm max}$ of $d$ is obtained when the the disks $D_2$ and $D_3$ are
touching each other (see Figure~\ref{fig:3diskslimit} (middle)). The
values of the capacity as a function of $d$ are shown in
Figure~\ref{fig:3diskslimit} (right) where the values of the capacity for
$2r<d<d_{\rm max}$ are computed by the BIE method and; for
$d=2r$ and $d=d_{\rm max}$ by the FEM. 
The minimal capacity is found when $d=2r$ and the maximal when the centers of the three disks form an equilateral triangle.

\subsection{One fixed disk. Three moving disks on a circle} 
\begin{figure} %
	\centering
	\subfloat{
		\includegraphics[width=0.275\textwidth]{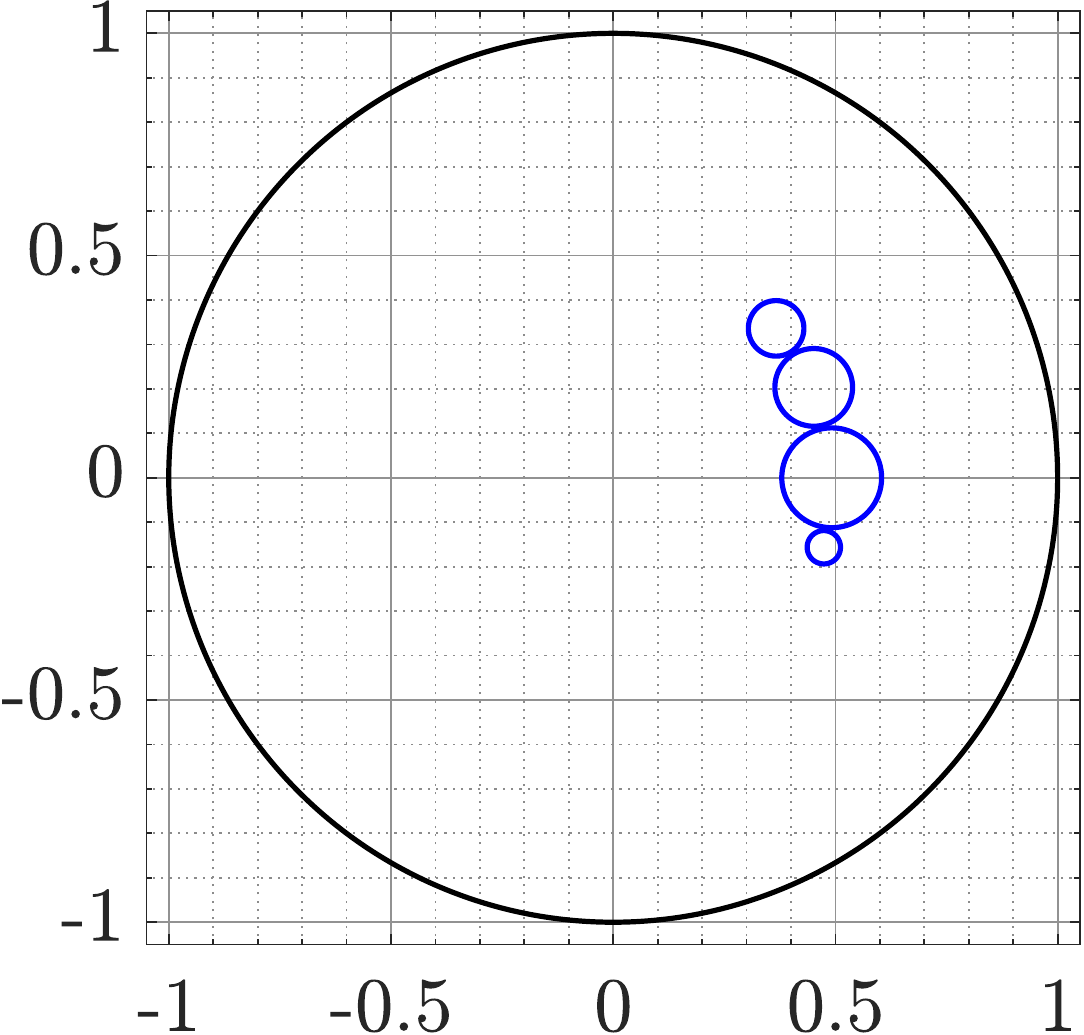}
	}
	\subfloat{
		\includegraphics[width=0.275\textwidth]{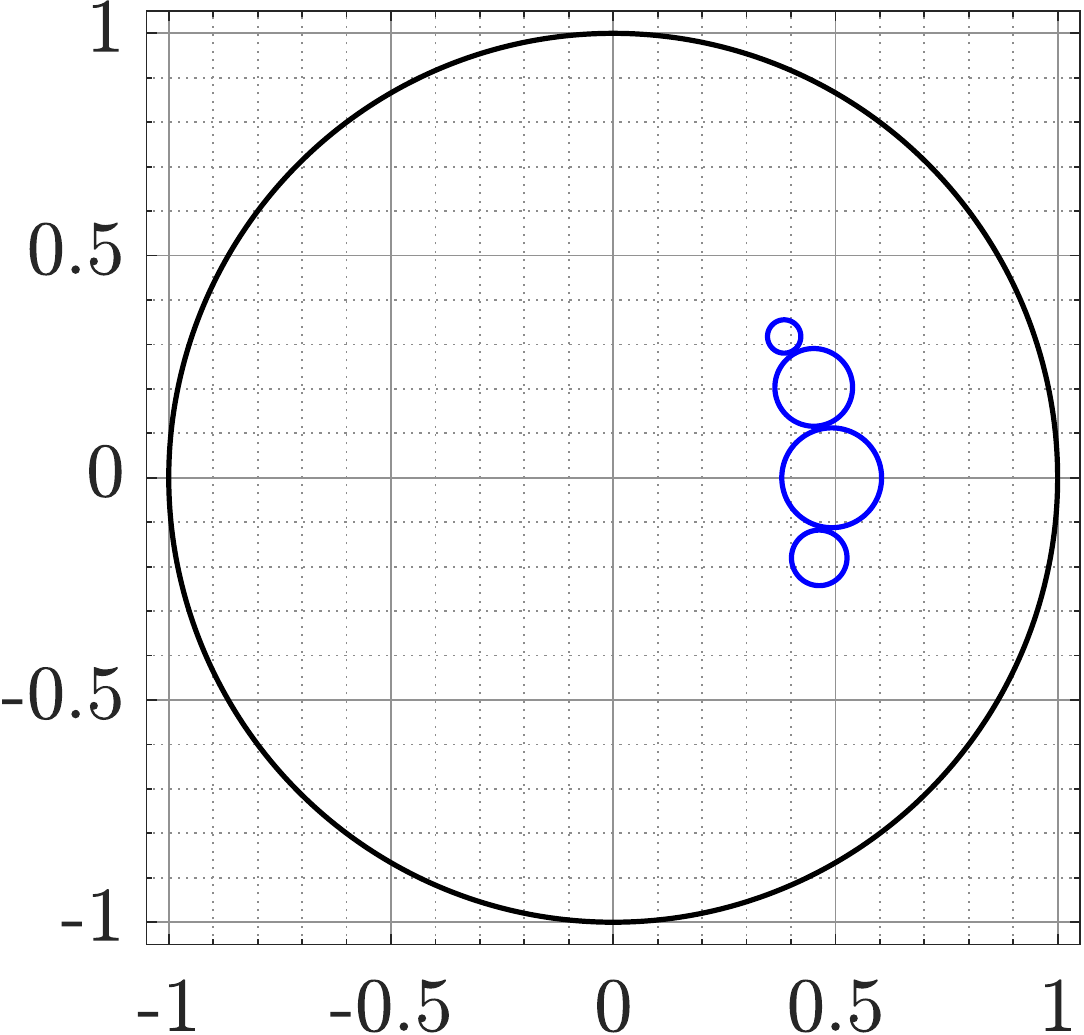}
	}
	\subfloat{
		\includegraphics[width=0.275\textwidth]{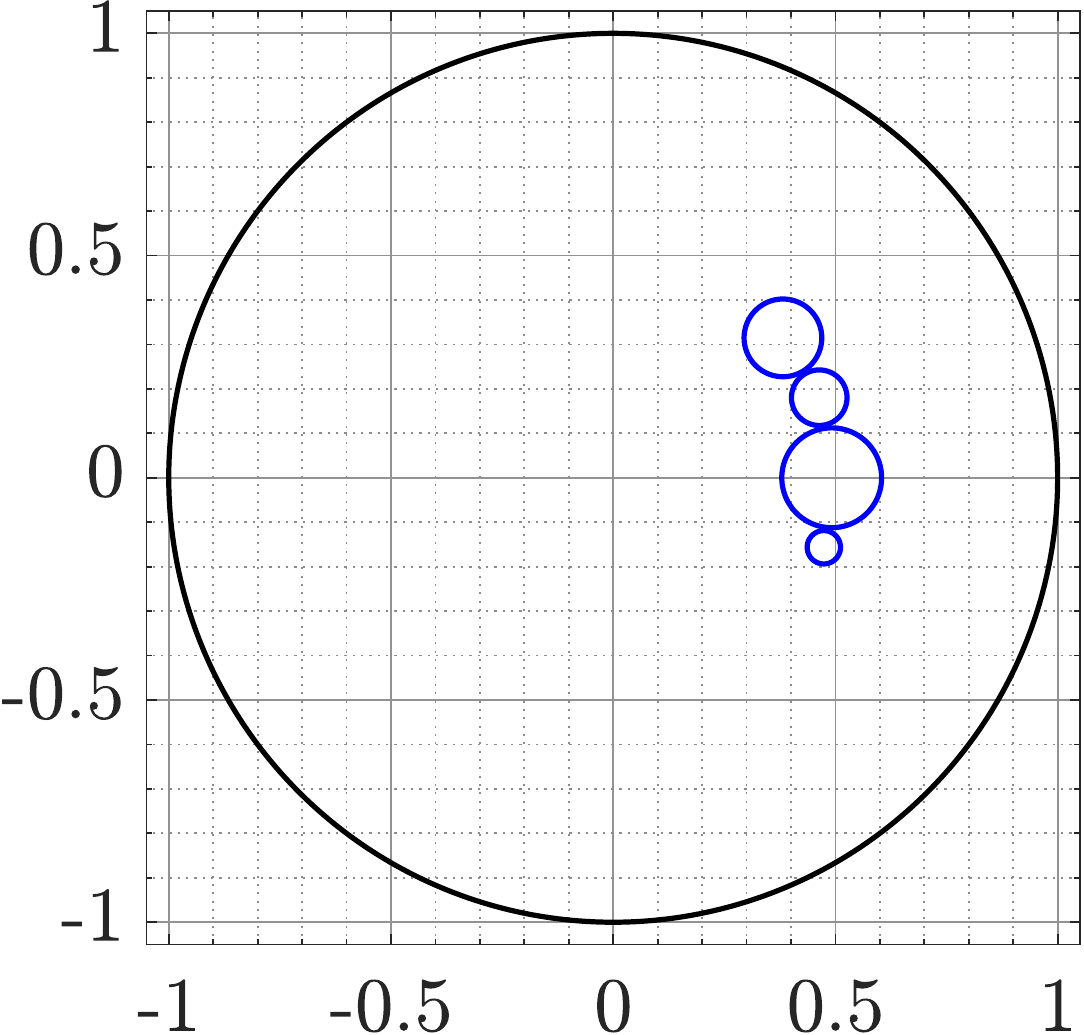}
	}
	\caption{Four disks with hyperbolic radii $3/30$, $5/30$, $7/30$, and $9/30$, and with centers on the circle $|z|=0.5$.
		Representative configurations of the optimised cases.} 
	\label{fig:opt-4disk-f}
\end{figure}

Staying on the circle $|z|=0.5$ we consider four disks with 
centers on the circle and hyperbolic radii 
$3/30$, $5/30$, $7/30$, and
$9/30$. 
Without any loss of generality, we will assume that the disk with hyperbolic
radius $9/30$ is fixed with its center on the positive real line at the point $0.5$. Then, we
search for the positions of the other three disks that minimize the capacity.
The initial positions of these three disks are assumed to be
$0.5e^{2k\pi\i/4}$ for $k=1,2,3$. For the optimized positions, we have
obtained six positions, with three
different values of the capacity due to symmetry (see Figure~\ref{fig:opt-4disk-f}). For the disks in the first column in
Figure~\ref{fig:opt-4disk-f}, the capacity is $4.6269$. The capacity is
$4.6193$ for the second column and $4.6621$ for the third column.

\subsection{One fixed disk. Three moving disks} 
\begin{figure} %
	\centering
	\subfloat{
		\includegraphics[width=0.275\textwidth]{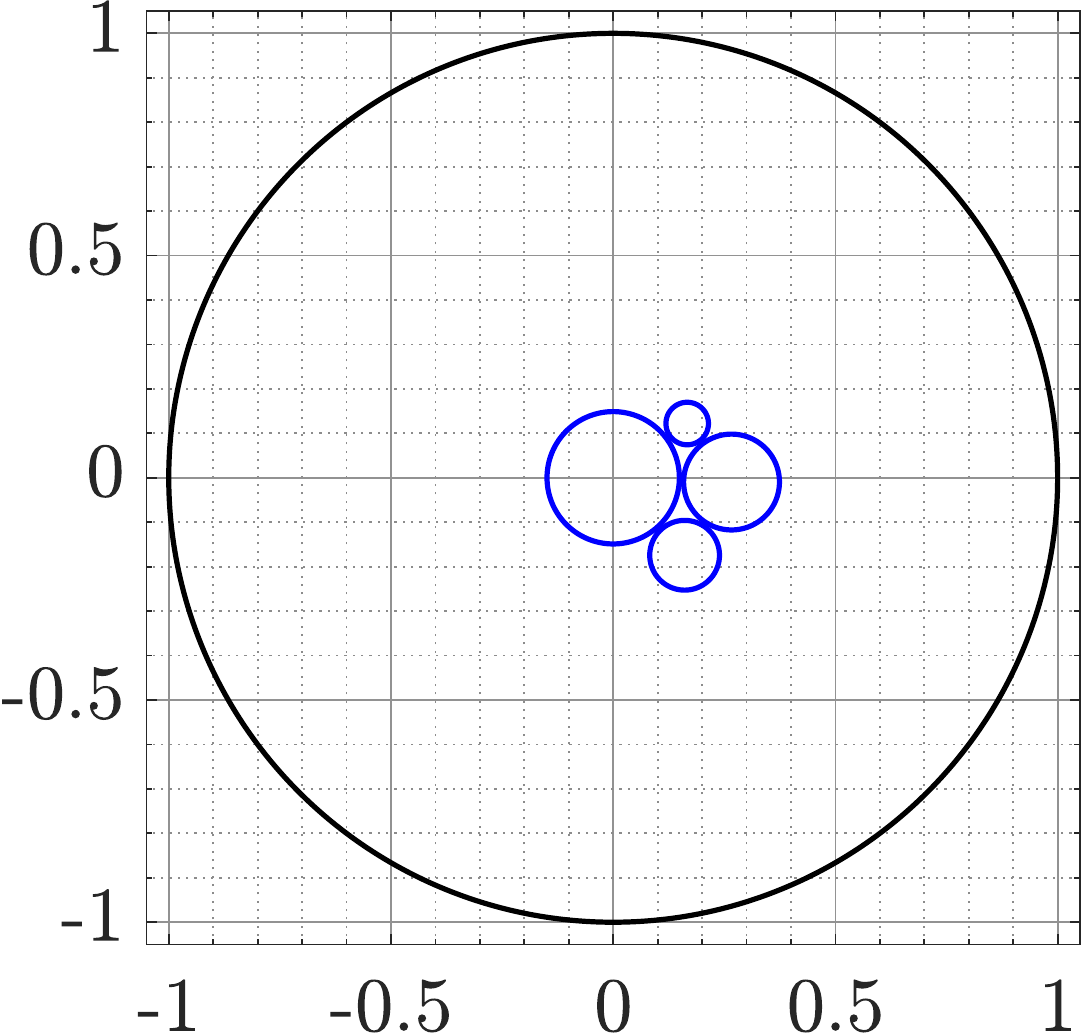}
	}
	\subfloat{
		\includegraphics[width=0.275\textwidth]{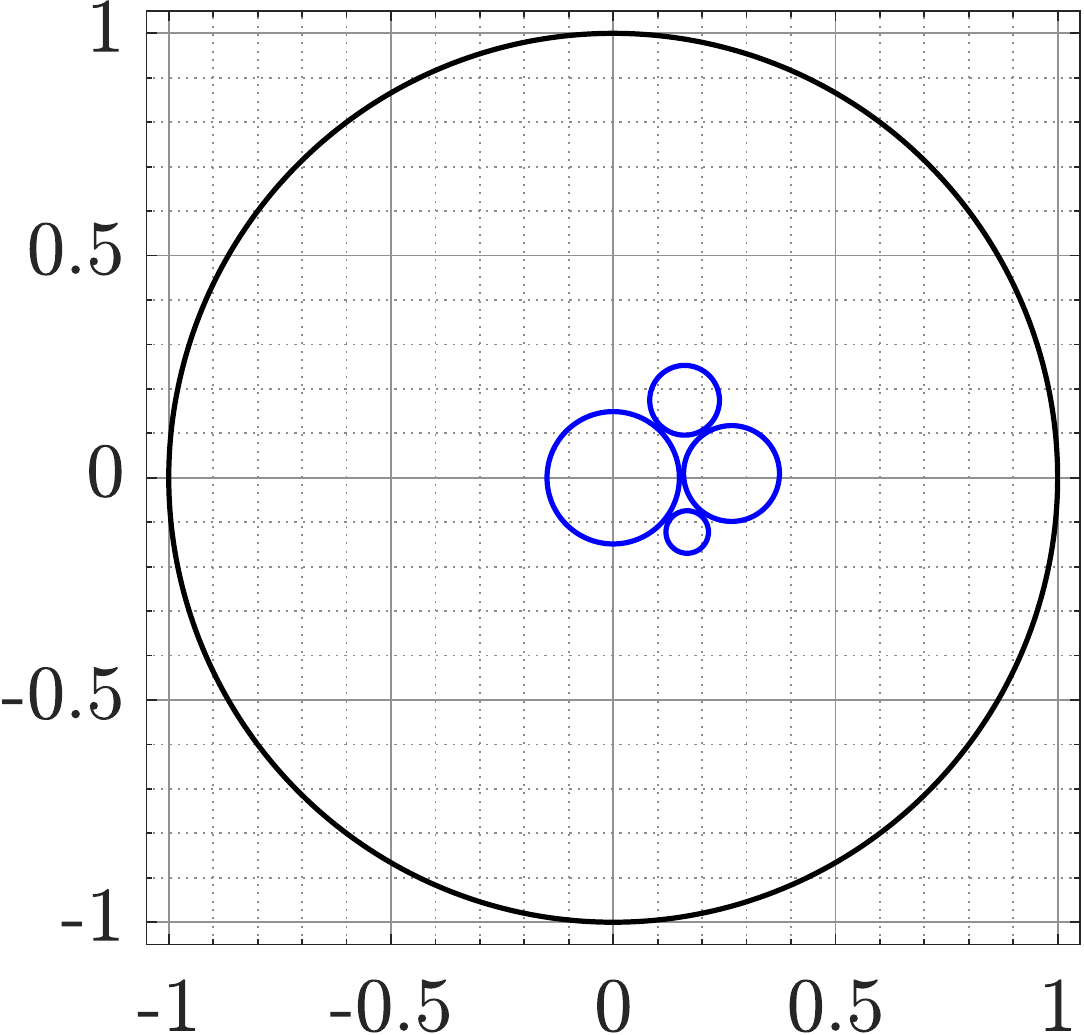}
	}
	\subfloat{
		\includegraphics[height=1.4in]{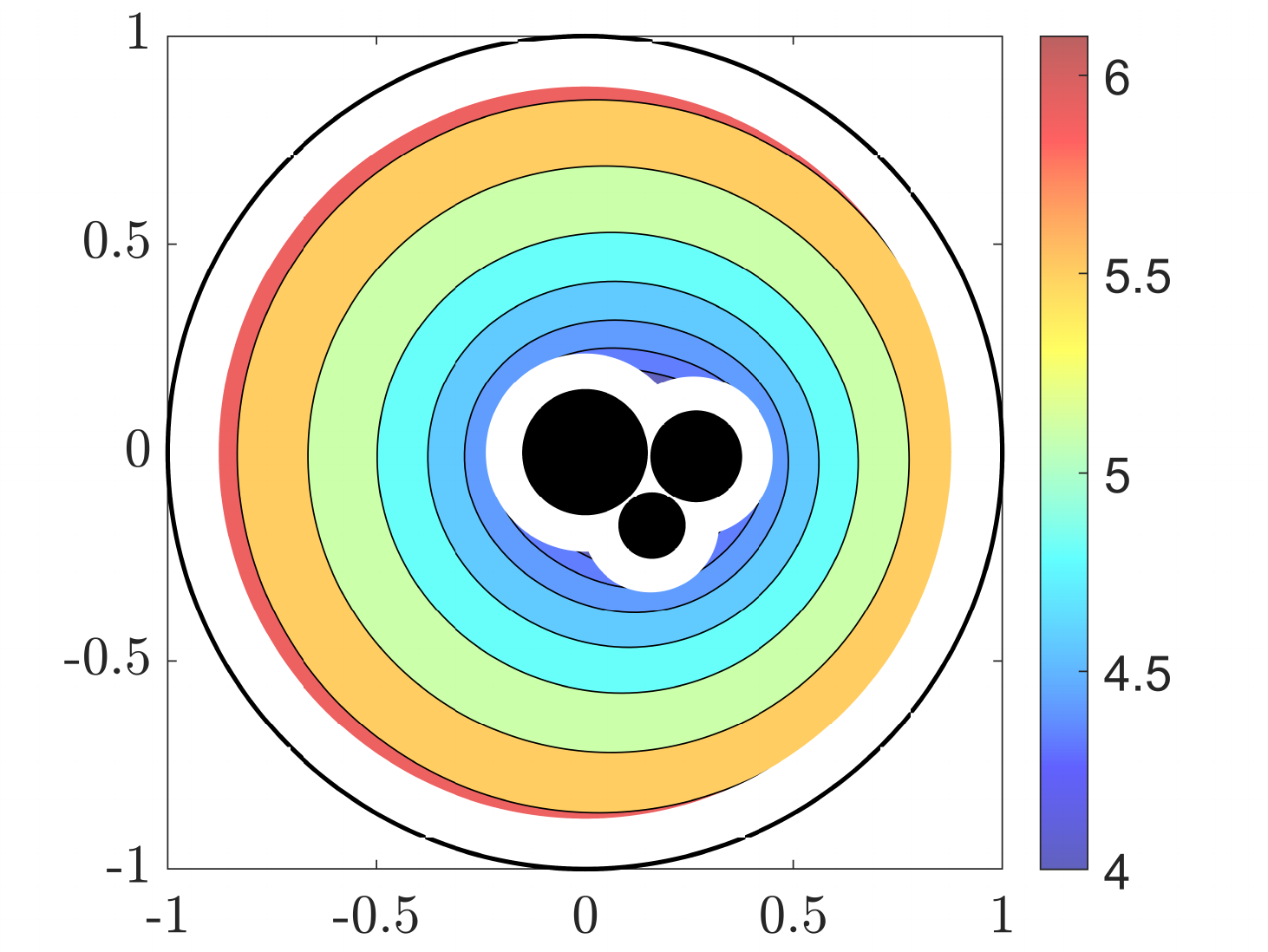}
	}
	\caption{Four disks with hyperbolic radii $3/30$, $5/30$, $7/30$, and $9/30$.
		The three largest disks have fixed positions, and the smallest one, centered at $z=x+i y$, is free to move. The level curves of
		$
		u(x,y)= {\rm cap}( \B, E),
		$
		where $E$ is the union of the four disks, indicate three local minima. 
		The two configurations on the left have converged to the global minimum.
	} 
	\label{fig:opt-4disk-0}
\end{figure}

Finally, we consider four disks with hyperbolic radii $3/30$, $5/30$, $7/30$, and $9/30$. 
This time, we will assume that the disk with hyperbolic radius $9/30$ is fixed with its center at the origin. 
The task is to find the positions of the three free disks that minimize the capacity where 
the initial positions of the three disks are assumed to be $0.5e^{2k\pi\i/3}$ for $k=0,1,2$. 
For the optimized positions, we have obtained two configurations, as shown in Figure~\ref{fig:opt-4disk-0}, 
with the capacity $4.2322$ which is the global minimum.

If we assume that the three disks with hyperbolic radii $5/30$, $7/30$, and $9/30$ have fixed positions as in Figure~\ref{fig:opt-4disk-0} (left), and the small disk with hyperbolic radius $3/30$ is moving. Assume that the center of the small disk is $z=x+\i y$ such that the four disks are non-overlapping. Let a function $u(x,y)$ be defined by
$
u(x,y)= {\rm cap}( \B, E),
$
where $E$ is the union of the four disks. 
The level curves of the function $u(x,y)$ are given in Figure~\ref{fig:opt-4disk-0} (right). 
As we can see from the figure, the capacity has three local minima and the capacity 
for the position in Figure~\ref{fig:opt-4disk-0} (left) is the global minimum. 
This experiment has been repeated multiple times with different initial starting positions
for the free disks and every one one of the local minima has been observed.

\subsubsection{On Computational Costs}
Naturally, the optimisation problems are the computationally most expensive ones
of all our numerical experiments.
In Table~\ref{tbl:timing} performance data on
the four disks free mobility problem is presented.
Comparison of the two methods is only qualitative, since both underlying
hardware and the interior-point implementations are different.
However, some conclusions can be derived.
In all cases the interior-point tolerance is the same, $\epsilon = 10^{-6}$, and within the
$hp$-FEM simulations, meshing is performed with the same discretization
control in every evaluation.
For optimal performance, the individual solutions must be accurate
enough so that the error induced by numerical approximation
of the gradients and Hessians is balanced with other sources of error.
For the $hp$-FEM it appears that the same mesh with $p=4$ is not adequate
in comparison with the one at $p=8$. Even though the time spent in one individual
iteration step is doubled, the overall time for $p=8$ is
significantly lower. 
At every evaluation the number of degrees of freedom is roughly 13000
(initial configuration: 13542, and final: 12589).
Similarly, for BIE the performance at $n=2^7$ is
superior to that at $n=2^4$.

\begin{table}
	\centering
	\caption{Solution times for the minimization process when
		one disk is fixed and three disks are mobile.
		Number of steps is number of iterations in the interior-point algorithm.
		Number of evaluations is the total number of solves performed during the
		minimization.
	}\label{tbl:timing}
	\begin{tabular}{llrrr}
		Method & Discretization & Time & Number of steps & Number of evaluations \\\hline
		BIE & $n=2^4$ & 472.9 & 151 & 1455\\
		& $n=2^7$     & 85.6  & 24  & 192 \\
		& $n=2^{10}$  & 150.7 & 24  & 192 \\\hline
		$hp$-FEM & $p=4$ & 39600& 202& 39568\\
		& $p=6$ & 11100 & 37 & 7494 \\
		& $p=8$ & 9100 & 20 & 4150 \\\hline
	\end{tabular}
\end{table}

The two implementations have very different requirements per iteration step.
Observe that the number of iteration steps is comparable, yet the
number of evaluations is not. The average time for one evaluation in BIE
is four to five times faster than one evaluation in $hp$-FEM.
Matlab and Mathematica results have been computed on modern Intel and
Apple Silicon computers, respectively.

\subsection{Hyperbolic area lower bound} 

\begin{figure} %
	\centering
	\subfloat{
		\includegraphics[width=0.23\textwidth]{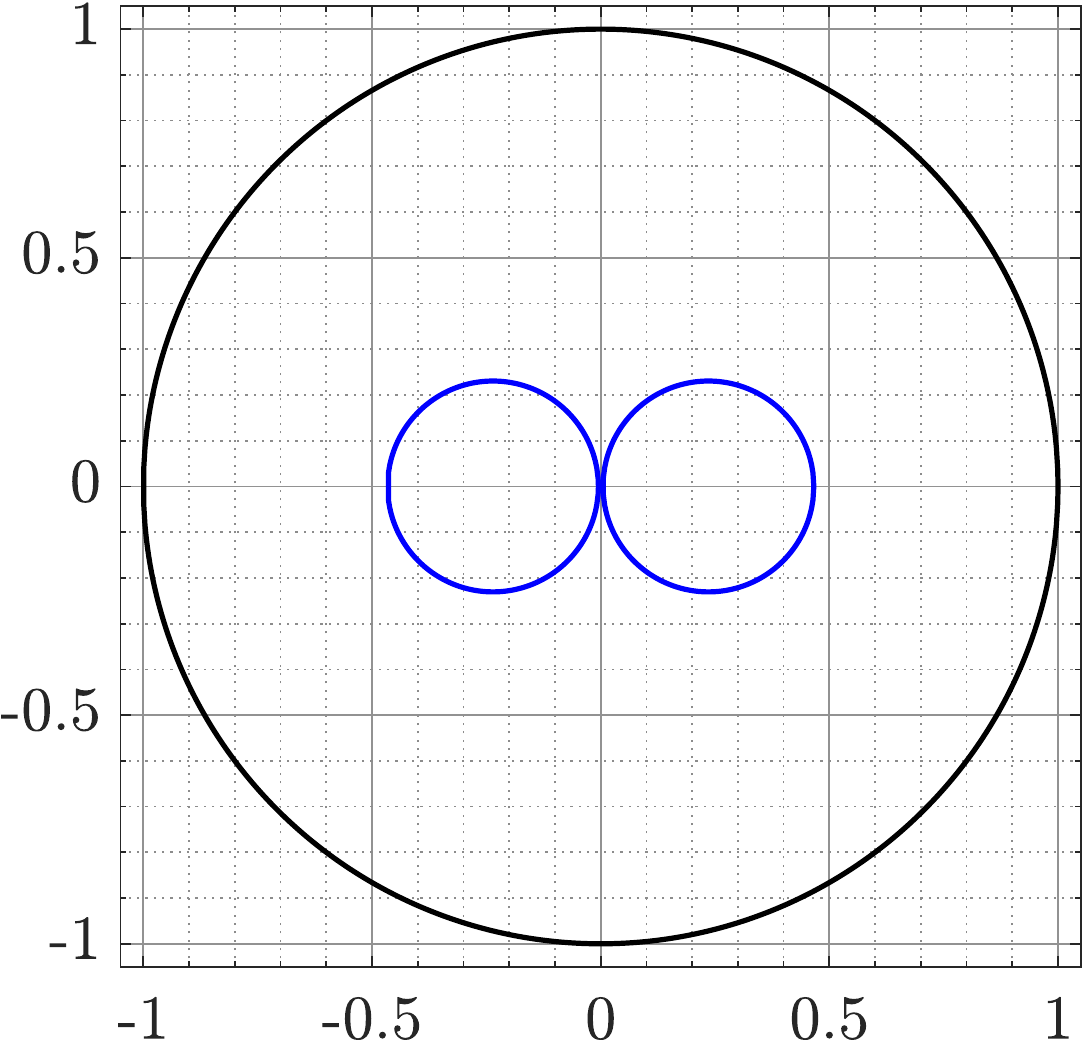}
	}
	\subfloat{
		\includegraphics[width=0.23\textwidth]{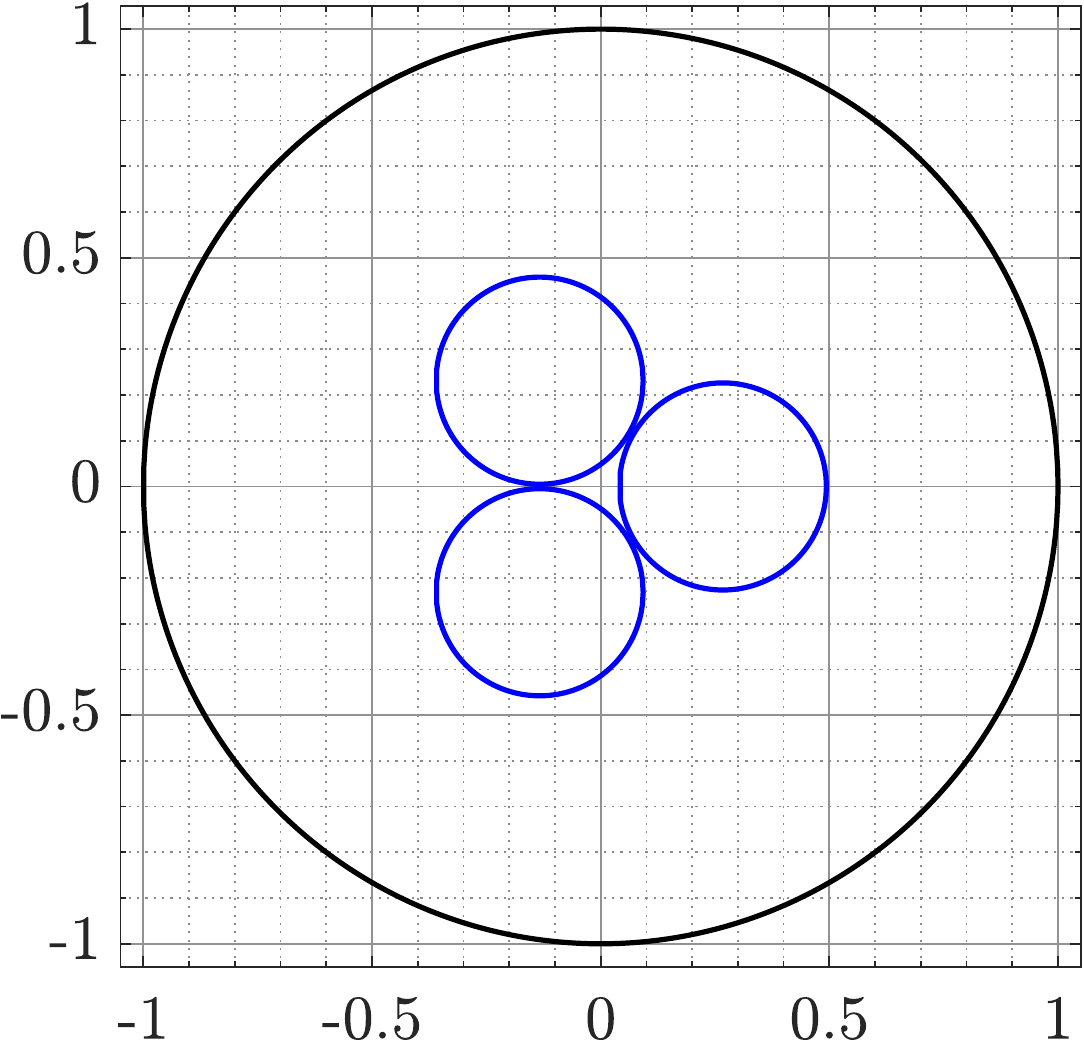}
	}
	\subfloat{
		\includegraphics[width=0.23\textwidth]{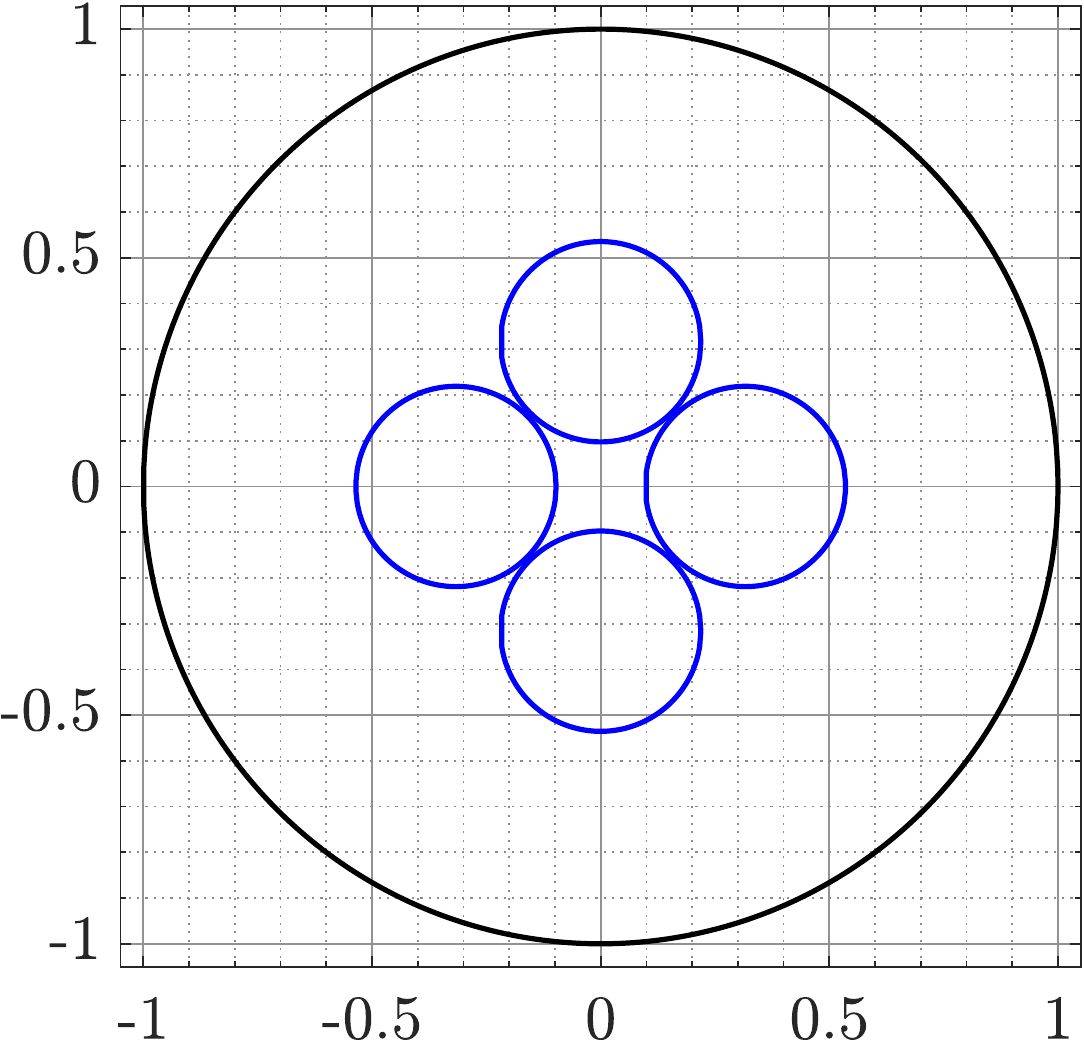}
	}
	\subfloat{
		\includegraphics[width=0.23\textwidth]{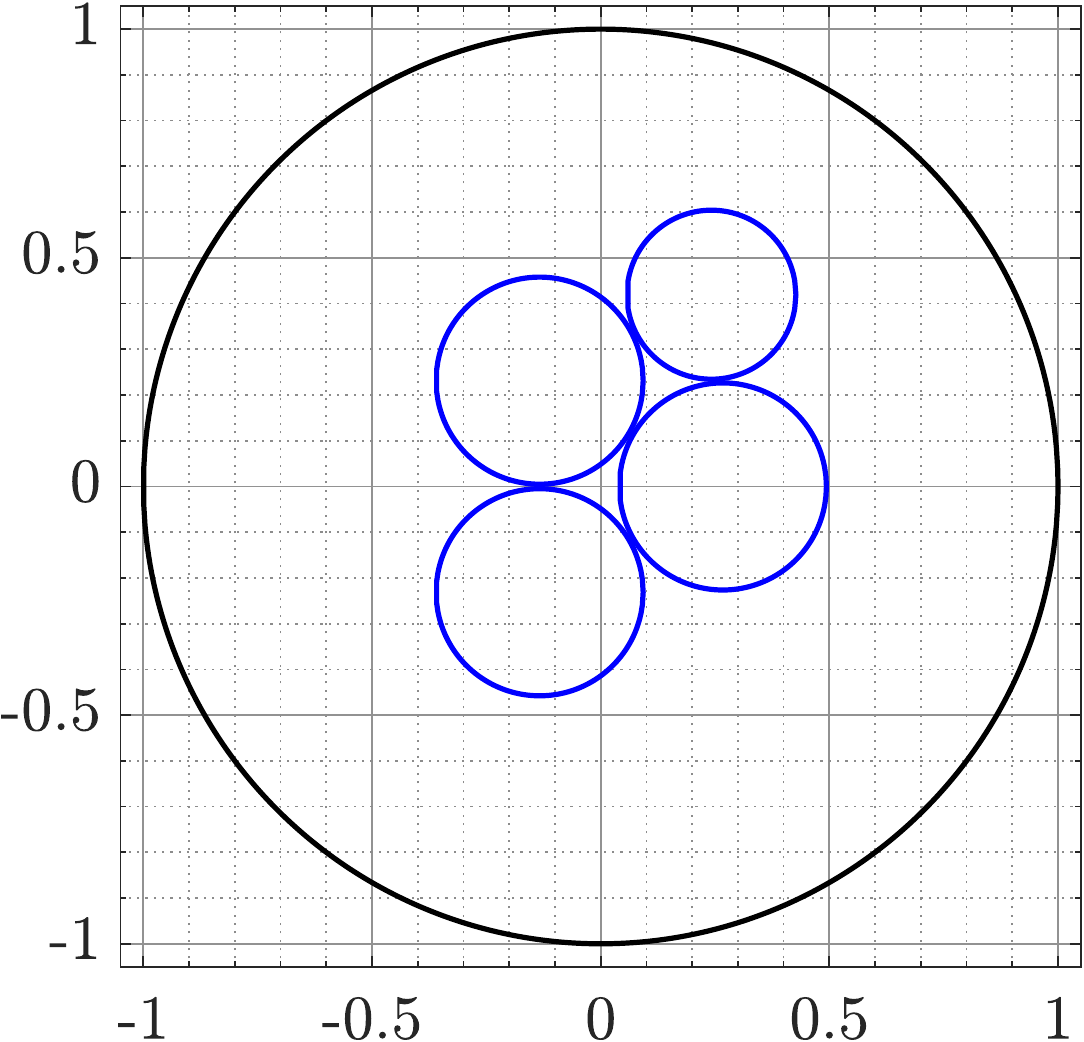}
	}
	\caption{The four types of condensers $(\B,E_r)$ for the hyperbolic radius $r=0.5$.
		From left: $m=2$, $m=3$, $m=4$ (Case I), $m=4$ (Case II).} 
	\label{fig:lowerbound}
\end{figure}

\begin{figure} %
	\centering
	\subfloat{
		\includegraphics[width=0.4\textwidth]{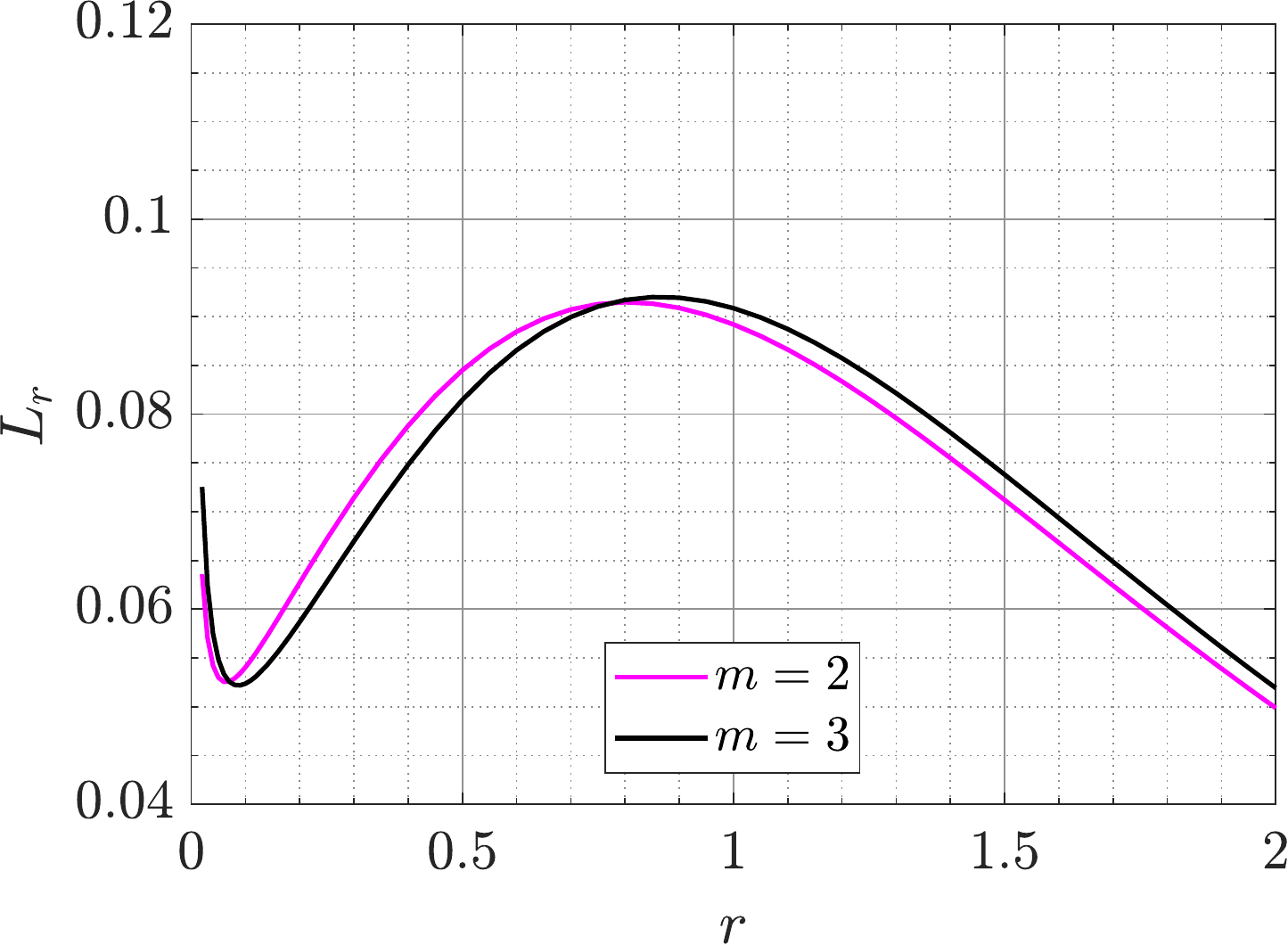}
	}	
	\subfloat{
		\includegraphics[width=0.4\textwidth]{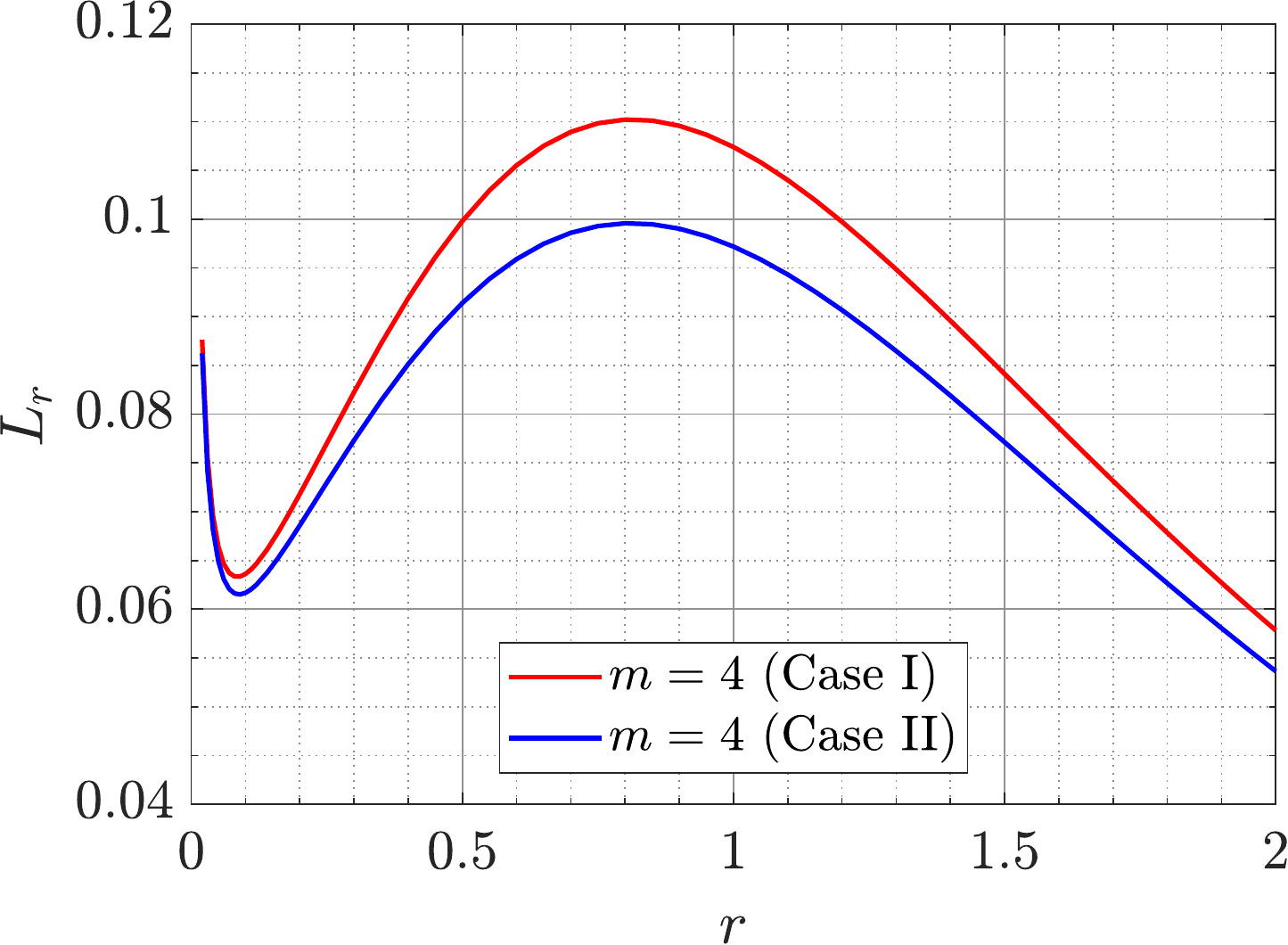}
	}
	\caption{The ratio $L_r$ for the four types of condensers $(\B,E_r)$.
		As $r \to \infty$ the improvement relative to the lower bound $L(r)$ tends to zero
		as expected. Left: For $m=2$, $m=3$, the improvements are very similar. Right: 
		In line with our experiments above, the Case II is indeed optimal, and gives us
		an improved lower bound.
	} 
	\label{fig:lowerboundcap}
\end{figure}

Finally, we compute the capacity of a constellation of disjoint hyperbolic disks and compare the computed values with the Hyperbolic area lower bound~\cite{G}. 
Let $E_r$ be the union of $m$ disjoint hyperbolic disks with equal hyperbolic radii $r$ such the hyperbolic distance between any two disks is $0.02$ (see Figure~\ref{fig:lowerbound} for $r=0.5$ and $m=2,3,4$). For $m=4$, we consider two cases (as shown in Figure~\ref{fig:lowerbound}) where the centers of the disks in Case I are on the real and imaginary axes. In Case II, the centers are on the rays $e^{\i\theta}$ for $\theta=0,\pi/3,2\pi/3,4\pi/3$. The hyperbolic area of these $m$ disks is $4m\pi\sh^2(r/2)$. Consider the hyperbolic disk $B_\rho(0,M)$ whose hyperbolic area is the same as the hyperbolic area of the $m$ disks, then  
\[
M=2\arsh\left(\sqrt{m}\sh(r/2)\right).
\] 
Then $L(r) = \capa(\B,B_\rho(0,M))$ is the hyperbolic area lower bound of $\capa(\B,E_r)$. In view of~\eqref{hypDat0}, we have
\[
L(r) = \capa(\B,B_\rho(0,M)) = \capa(\B,B^2(0,\th (M/2))) = \frac{2\pi}{\log\cth(M/2)}.
\]
The BIE method is then used to compute $\capa(\B,E_r)$ for several values of $r$ with $0.02\le r\le 2$. 
Our computed minimum value of the capacity can be considered a
lower bound of the capacity of the constellation of $m$ disjoint
hyperbolic disks. We compare the computed value with the hyperbolic area lower bound by defining
\[
L_r = \frac{\capa(\B,E_r)-L(r)}{L(r)}.
\]
The graph of $L_r$ is shown in Figure~\ref{fig:lowerboundcap} for $0.02\le r\le 2$ and $m=2,3,4$. 
As $r \to \infty$ it appears that the improvement tends to zero.
This is a consequence of the nature of hyperbolic geometry. With one disk fixed in the centre
the other three will have ever smaller contributions to the capacity since their
Euclidean areas tend to zero as in Figure~\ref{fig:hypgeom} (right).
It is an indication of the complexity of the problem that the 
graphs in Figure~\ref{fig:lowerboundcap} do not reveal any simple connection between
the number of the disks and the minimal capacity.

	\section{Conclusions}
We study  lower bounds for the conformal capacity of a  constellation of disjoint hyperbolic disks $E_j \subset \mathbb{B}^2$, $j=1,...,m$,
using a novel idea: instead of using a symmetrization transformation,
which usually leads to fusion of the disjoint disks, we are looking for a lower
bound in terms of another constellation which yields a minimal value. The traditional
symmetrization transformation  \cite{PS}, \cite{B},
\cite{Du}, is now replaced by free mobility of individual
disks with the constraint that the hyperbolic radii of the disks are 
invariant and the disks are non-overlapping. In this process, due to the 
conformal invariance, the conformal capacity of each disk stays 
invariant, whereas the capacity of the whole constellation may significantly
vary. Moreover, the hyperbolic area of the constellation is also constant.

The optimization methods we used produced (locally) minimal constellations
such that the disks group together, as closely as possible. This coalescing is reminiscent of the behavior of some animal colonies in cold
weather conditions for the purpose of heat flow minimization. 
Mathematical methods are not available for analytic treatment of the problems, 
but we are convinced that there is a strong connection with combinatorial geometry, topics like
packing and covering problems. Such problems often have many local minima
\cite[p. 157]{ber}.

We carried out numerical simulations using two different methods, the BIE and $hp$-FEM
methods and the close agreement of the two computational methods confirmed
the  results. Because of the complexity
of the problem we studied various subproblems where disk centers satisfied
constraints such that the centers are on the interval $(-1,1)$ or  at the
same distance from the origin. In both cases we observed the grouping phenomenon
(cf. Figure~\ref{fig:appetizer}) and, moreover, noticed that permutation of disks has influence 
on the capacity if the radii are different. Because the hyperbolic area 
of a constellation is a constant, it is now clear that the hyperbolic 
area alone does not define the constellation capacity.

This observation led us to compare our computed lower bound to Gehring's
sharp lower bound given in terms of hyperbolic area. The conclusion was that we 
obtained in some cases approximately $10\%$ improvement when $m=4$.

The numerical agreement of the BIE and $hp$-FEM methods was very good, typically
ten decimal places or better, and the expected exponential convergence
was observed, see Figure~\ref{fig:fourmesherror}. The performance of the BIE method was significantly
faster than the $hp$-FEM method when it comes to computational time and flexilibity to 
modify the code to new situations. This is probably due to the heavy data 
structure of the $hp$-FEM method due to hierarchial triangulation refinement 
process of the method.

A vast territory of open problems remains. First, it would be interesting to
study whether some kind heuristic methods would lead to "close to extremal"
constellations, to be used as initial steps of the minimization. 
Such a method could be based on some computationally cheaper
object function than the capacity itself: for instance, first, the maximization of the
number of the mutual contact points of the constellation. 
Second, the case of $m>5$ disks of equal radii seems to be completely open. 
Perhaps in this case the number of locally minimal constellations grows exponentially as a function of $m.$ 
Third, one could study constellations of other types of geometric figures like hyperbolic triangles.

\bibliographystyle{amsplain}

\end{document}